\documentclass[sn-mathphys,Numbered]{sn-jnl}

\usepackage{bm}
\usepackage{graphicx}
\usepackage{epstopdf}
\usepackage{subfig}
\usepackage{multirow}%
\usepackage{amsmath,amssymb,amsfonts}%
\usepackage{amsthm}%
\usepackage{mathrsfs}%
\usepackage[title]{appendix}%
\usepackage{xcolor}%
\usepackage{textcomp}%
\usepackage{manyfoot}%
\usepackage{booktabs}%
\usepackage{algorithm}%
\usepackage{algorithmicx}%
\usepackage{algpseudocode}%
\usepackage{listings}%
\usepackage{rotating}
\usepackage{natbib}

\theoremstyle{thmstyleone}%
\theoremstyle{thmstyletwo}%
\theoremstyle{thmstylethree}%
\makeatletter
\newcommand{\biggg}{\bBigg@{3}}

\newcommand{\Biggg}{\bBigg@{3.5}}

\newcommand{\bigggg}{\bBigg@{4}}

\newcommand{\Bigggg}{\bBigg@{4.5}}

\makeatother
\begin{document}

\title[Article Title]{The macroscopic finite-difference scheme and modified equations of the general propagation multiple-relaxation-time lattice Boltzmann model}


\author[1,2,3]{\fnm{} \sur{Ying Chen}}

\author[1,2,3]{\fnm{} \sur{Xi Liu}}

\author*[1,2,3]{\fnm{} \sur{Zhenhua Chai}}\email{hustczh@hust.edu.cn}

\author[1,2,3]{\fnm{} \sur{Baochang Shi}}

\affil[1]{\orgdiv{School of Mathematics and Statistics}, \orgname{Huazhong
		University of Science and Technology}, \orgaddress{ \city{ Wuhan} \postcode{430074}, \state{Hubei}, \country{China}}}

\affil[2]{\orgdiv{Institute of Interdisciplinary Research for Mathematics and Applied Science}, \orgname{Organization}, \orgname{Huazhong
		University of Science and Technology}, \orgaddress{ \city{ Wuhan} \postcode{430074}, \state{Hubei}, \country{China}}}
	
\affil[3]{\orgdiv{Hubei Key Laboratory of Engineering Modeling and
		Scientific Computing}, \orgname{Organization}, \orgname{Huazhong
		University of Science and Technology}, \orgaddress{ \city{ Wuhan} \postcode{430074}, \state{Hubei}, \country{China}}}


\abstract{In this paper, we first present the general propagation multiple-relaxation-time lattice Boltzmann (GPMRT-LB) model and obtain the corresponding macroscopic finite-difference (GPMFD) scheme on conservative moments.  Then based on the Maxwell iteration method, we conduct the analysis on the truncation errors and modified equations (MEs) of the GPMRT-LB model and GPMFD scheme at both diffusive and acoustic scalings.     For the nonlinear anisotropic convection-diffusion equation (NACDE) and Navier-Stokes equations (NSEs), we also derive the  first- and second-order MEs of the GPMRT-LB model and GPMFD scheme.   In particular, for the one-dimensional convection-diffusion equation (CDE) with the constant velocity and diffusion coefficient, we can develop a fourth-order GPMRT-LB (F-GPMRT-LB) model and the corresponding fourth-order GPMFD (F-GPMFD) scheme at the diffusive scaling.  Finally, two benchmark problems, Gauss hill problem and Poiseuille flow in two-dimensional space, are used to test the GPMRT-LB model and GPMFD scheme, and it is found that the numerical results are not only in good agreement with corresponding analytical solutions, but also have a second-order convergence rate in space.     Additionally, a numerical study on one-dimensional CDE also demonstrates that the F-GPMRT-LB model and F-GPMFD scheme can achieve a fourth-order accuracy in space, which is consistent with our theoretical analysis.}




\maketitle

\section{Introduction}

The kinetic theory based lattice Boltzmann (LB) method, as a highly efficient numerical approach at the mesoscopic level, has been widely used to study the fluid flow problems  (e.g., the multiphase flows \cite{Kruger2017,Wang2019}, fluid flows in porous media \cite{Succi2001}) governed by the Navier-Stokes equations (NSEs) \cite{Guo2013}  for its a second-order accuracy in space \cite{Junk2005,Dellar2013} and advantange in treating complex boundary conditions. On the other hand, the LB method has also been extended to solve some special kinds of  partial differential equations (PDEs), including diffusion equations \cite{Ancona1994, Rasin2005, Huber2010,Suga2010,Lin2022,Silva2023,Chen2023-diffusion}, convection-diffusion equations (CDEs) \cite{Gin2005, Chopard2009,Yoshida2010,Gin2013,Chai2013,Jettestuen2016,Li2019,Straka2020,Michelert2022,Chen2023},  Burgers equations \cite{Li2012,Li2018,Rong2023}, general real and complex nonlinear convection-diffusion equations \cite{ Shi2009}, and nonlinear anisotropic convection-diffusion equations (NACDEs) \cite{Chai2016}.

Usually, the LB method suffers from numerical instability when the relaxation parameter is close to 2. To solve the problem, there are two possible approaches that can be adopted. The first one is to introduce the multiple-relaxation-time (MRT) collision operator \cite{DdHum2002, Chai2020, Chai-Shi-2023} with some adjustable free relaxation parameters \cite{Lallemand2000, Pan2006, Cui2016, Luo2011},  which is more general and more stable than the single- and two-relaxation-time lattice Boltzmann (SRT-LB \cite{Qian1992} and TRT-LB \cite{Gin2005}) models.   The second one is to use the general propagation LB model within the framework of the time-splitting method where two free parameters are introduced into the propagation step for NSEs \cite{Guo2001} or NCDE \cite{Guo2018} based on the Lax-Wendroff (LW) scheme \cite{Mcnamara1995, Zhang2001} and fractional propagation (FP) scheme \cite{Qian1997}. This model is more stable, and the popular standard LB model, LW  and FP schemes can be viewed as its special cases. Considering the advantages of the MRT-LB model and general propagation LB model  in the numerical stability, in this work, we will consider the more general propagation MRT-LB (GPMRT-LB) model.

In the framework of LB method, several asymptotic analysis approaches have been used to derive macroscopic PDEs, including the Chapman-Enskog analysis \cite{Chapman1990}, Maxwell iteration \cite{Ikenberry1956, Yong2016}, direct Taylor expansion \cite{Chai2020}, recurrence equations method \cite{Gin2012, GinzburgD2009}, and equivalent equations method \cite{Dubois2008, Dubois2009, Dubois2019}, while they cannot be applied to clarify the relation between the LB model and the macroscopic PDE based numerical scheme (the so-called macroscopic numerical scheme). Over the past years, some significant contributions have been made to bridge the gap between the LB model and macroscopic numerical scheme for a specified PDE, and these efforts aim to provide clear and rigorous consistency, accuracy, and the derivation of modified equation (ME). 
Most of the existing works, however, are limited to the LB models and the  macroscopic finite-difference schemes for diffusion equations and CDEs. On the one hand, Ancona \cite{Ancona1994} first presented the SRT-LB model with D1Q2 lattice structure for the one-dimensional diffusion equation, and found that it is consistent with a macroscopic three-level second-order DuFort-Frankel scheme \cite{Fort1953}. Then Suga \cite{Suga2010} demonstrated that the SRT-LB model with D1Q3 lattice structure for the one-dimensional diffusion equation is equivalent to a macroscopic four-level fourth-order finite-difference scheme. Lin et al. \cite{Lin2022} extended this work to consider a more general MRT-LB model, and derived a macroscopic four-level sixth-order finite-difference scheme. Silva \cite{Silva2023} focused on the TRT-LB model for the diffusion equation with a linear source term, and obtained a macroscopic fourth-order finite-difference scheme. However, all of the aforementioned works are only limited to one-dimensional problems. In a recent work, Chen \cite{Chen2023-diffusion} considered  the MRT-LB model with the D2Q5 lattice structure for the two-dimensional diffusion equation, and obtained a five-level fourth-order finite-difference scheme. On the other hand, Dellacherie \cite{Dellacherie2014} analyzed the SRT-LB model for the one-dimensional CDE with D1Q2 lattice structure, and illustrated that this LB model is equivalent to a three-level finite-difference scheme called LFCCDF (Leap-Frog difference for the temporal derivative, central difference for the convective term, and DuFort-Frankel approximation for the diffusive term) scheme \cite{Kwok1993}. Cui et al. \cite{Cui2016} showed that for the one-dimensional steady CDE, the MRT-LB model can be written as a macroscopic second-order central-difference scheme. Following the similar idea,  Wu et al. \cite{Wu2020} derived a macroscopic finite-difference scheme of the MRT-LB model composed of natural moments, and further performed a more general analysis on the discrete effects  of some boundary schemes  for the CDE. Recently, Chen et al. \cite{Chen2023} obtained  a macroscopic four-level fourth-order finite-difference scheme from the MRT-LB model with D1Q3 lattice structure for the CDE. In addition, Li et al. \cite{Li2012} demonstrated that for the one-dimensional Burgers equation, the SRT-LB model with the D1Q2 lattice structure can be expressed as a macroscopic three-level second-order finite-difference scheme.  Junk \cite{Junk} and Inamuro \cite{Inamuro}
found that the SRT-LB model is equivalent to a macroscopic two-level finite-difference
scheme if the relaxation parameter is equal to one, and at the diffusive scaling, it has a second-order convergence rate for the incompressible NSEs \cite{Junk}.  d$'$Humi$\grave{{\rm e}}$re and Ginzburg \cite{GinzburgD2009} conducted a theoretical analysis on the TRT-LB model with recurrence equations, and illustrated that when the magic parameter is fixed as $\Lambda^{eo}= 1/4$, the model can be written as a macroscopic three-level finite-difference scheme with a second-order accuracy in space. Chai et al. \cite{Chai2020} further presented a more general analysis on the TRT-LB model, and also derived the three-level finite-difference schemes for steady and unsteady problems.  

It is worth noting that the works mentioned above are limited to some specific problems and/or lattice structures. To obtain the macroscopic finite-difference scheme from a given LB model with the D$d$Q$q$ ($q$ discrete velocities in $d$-dimensional space) lattice structure, Fu\v{c}\'{i}k et al. \cite{Fucik2021} developed a general computational tool \cite{Fucik2023}, while the origin of this algorithm remains unclear. In contrast, Bellotti et al. \cite{Bellotti2022} conducted a precise algebraic characterization of the LB model, and investigated the relationship between the MRT-LB model and macroscopic numerical scheme. They found that the LB model can be expressed as a macroscopic multiple-level  finite-difference scheme  with a second-order accuracy. Furthermore, they also carried out some analysis on the truncation errors and MEs at both  diffusive and acoustic scalings \cite{Bellotti-1-2022}, which are consistent with the results based on the asymptotic analysis methods \cite{Chai2020}. However, it should be noted that they  only considered the MRT-LB model  with a diagonal relaxation matrix and the first-order ME at the diffusive scaling \cite{Bellotti-1-2022}.  In this work, we will first extend the previous works \cite{Bellotti2022,Bellotti-1-2022} to consider the more general GPMRT-LB models that are developed for NACDE and NSEs, and  derive the macroscopic finite-difference (GPMFD) schemes. Then we will conduct some detailed analysis on the truncation errors, the first- and second-order MEs of the GPMRT-LB models and GPMFD schemes at both diffusive and acoustic scalings.

The remainder of this paper is organized as follows.  In Section \ref{Derivation-GPMFD}, we present the details on how to derive  the GPMFD scheme on conservative moments from the GPMRT-LB model. In Section \ref{Truncation-Error-Modified-Equation}, the truncation errors at both diffusive and acoustic scalings are derived through the  Maxwell iteration method, followed by the first- and second-order MEs of the GPMRT-LB model and GPMFD scheme. In section \ref{Fourth-Order-for-CDE}, we develop a fourth-order GPMRT-LB (F-GPMRT-LB) model and GPMFD (F-GPMFD) scheme at the diffusive scaling for the one-dimensional CDE with the constant velocity and diffusion coefficient. In Section \ref{Numerical-Simulation}, some simulations of the  Gauss hill problem, Poiseuille flow, and   one-dimensional CDE are carried out to test  the proposed GPMRT-LB model and GPMFD scheme. Finally, some conclusions are given in Section \ref{Conclusion}.

	\section{The GPMFD scheme of the GPMRT-LB model}\label{Derivation-GPMFD}
\subsection{Preparation}\label{preparation}
To begin our analysis, we first discretize the  problems in $d$ $(d=1,2,3)$ dimensional space without considering the boundary conditions. In the LB method \cite{Guo2013},  the space  is discretized by $\mathcal{L}:=\Delta x\mathcal{Z}^d$ with a constant lattice spacing $\Delta x>0$ in all directions, and the  more general  rectangular lattice structure \cite{Chai-Shi-2023} is not considered here. The time is uniformly discretized by $\mathcal{T}:=\Delta t\mathcal{N}$ with $t_n:=n\Delta t,n\in\mathcal{N}$,  $\Delta t$ is the time step. Additionally, we also introduce the so-called lattice velocity, defined by $\lambda:=\Delta x/\Delta t$. It should be noted that the discretizations of the spatial and temporal domains are completely independent on the scaling between $\Delta x$ and $\Delta t$, this means that it does not have an influence on the derivation of the macroscopic finite-difference schemes of the LB models. However, the scaling has a significant effect on the consistency analysis, i.e., the truncation errors and MEs (see Section \ref{Truncation-Error-Modified-Equation} for details). 

It is known that in the LB method, the evolution process can be split into the collision and propagation steps, and to simplify the following analysis on the derivation of the macroscopic finite-difference scheme,  it is necessary to introduce the time and the space operators associated with the discrete velocities. \\
\textbf{Definition 1} Let $\bm{z}\in\mathcal{Z}^d$ and $z\in \mathcal{Z}$, the space shift operator on the space lattice $\mathcal{L}$ denoted by $T_{\Delta x}^{\bm{z}}$ and the time shift operator on the time lattice $\mathcal{T}$ represented  by $T^{z}_{\Delta t}$, are defined as follows:
\begin{subequations}\label{de-time-space}
	\begin{align}
		&\big[T_{\Delta x}^{\bm{z}}h\big](\bm{x},t)=h(\bm{x}+\bm{z}\Delta x,t),\bm{x}\in\mathcal{L},t\in\mathcal{T},\label{Space-oper}\\
		&\big[T^{z}_{\Delta t}h\big](\bm{x},t)=h(\bm{x},t+z\Delta t),\bm{x}\in\mathcal{L},t\in\mathcal{T},\label{Time-oper}
	\end{align}
\end{subequations}
where $h(\bm{x},t)$ is an arbitrary function defined on the space lattice $\mathcal{L}$ and time lattice $\mathcal{T}$. According to Eq. (\ref{de-time-space}), the space and time shift operators can also be expressed in the series forms:
\begin{subequations}
	\begin{align}
		&T_{\Delta x}^{\bm{z}}=\sum_{k=0}^{+\infty}\frac{\Delta x^k\big(\bm{z}\cdot\nabla\big)^k}{k!},\label{series-space-oper}\\
		&T^z_{\Delta t}=\sum_{k=0}^{+\infty}\frac{(z\Delta t)^k\partial_t^k}{k!},\label{series-time-oper}
	\end{align}
\end{subequations}
where the gradient operator $\nabla=\Big(\partial_{x_1},\partial_{x_2},\ldots,\partial_{x_d}\Big)$.

\subsection{GPMRT-LB  model}
In the LB method, there are some discrete procedures that can be used to solve the discrete velocity Boltzmann equation \cite{Guo2001,Bao2006},

\begin{align}\label{DVBE}
	\frac{\partial f_i}{\partial t}+\bm{c}_i\cdot\nabla f_i=\Omega_i+F_i,
\end{align}
where the function $f_i$ represents the particle distribution at position $\bm{x}$ and time $t$, $\{\bm{c}_i=c\bm{e}_i,i=1,2,\ldots,q\}$ denotes the set of discrete velocities in D$d$Q$q$ lattice structure, $\Omega_i$ is the general collision operator, and $F_i$ is the discrete source or force term. Based on the previous works \cite{Guo2001,Guo2018}, we can develop a GPMRT-LB model  for the NACDE and NSEs.

With the time-splitting method, Eq. (\ref{DVBE}) can be separated into two steps:
\begin{subequations}
	\begin{align}
		&\frac{\partial f_i}{\partial t}=\Omega_i+F_i,\label{collision}\\
		&\frac{\partial f_i}{\partial t}+\bm{c}_i\cdot\nabla f_i=0,\label{propagation}
	\end{align}
\end{subequations}
which are the collision and propagation steps, respectively.  

Following the approach presented in Ref. \cite{Chai2020}, the collision step in 
 Eq. (\ref{collision})  can be reformulated as 
\begin{align}\label{collison2}
	&f^{\star}_i(\bm{x},t)=f_i(\bm{x},t)-\bm{\Lambda}_{ik}f_k^{ne}(\bm{x},t)+\Delta t\bigg[G_i+F_i+\frac{\Delta t}{2}\overline{D}_iF_i\bigg](\bm{x},t),i=1,2,\ldots,q.
\end{align}
Here $f^{\star}(\bm{x},t)$ denotes the post-collision distribution function, $\bm{\Lambda}=\big(\bm{\Lambda}\big)_{ik}$ is a $q\times q$ invertible collision matrix and can be defined as $\bm{\Lambda}:=\bm{M}^{-1}\bm{SM}$, $\bm{M}$ and $\bm{S}$ are the invertible transform and relaxation matrices, respectively. $f^{ne}_i(\bm{x},t)=f_i(\bm{x},t)-f_i^{eq}(\bm{x},t)$ represents the nonequilibrium distribution function, $G_i$ is the auxiliary source distribution function and can be used to remove additional terms. $\overline{D}_i=\partial_t+\gamma\bm{c}_i\cdot\nabla$ with $\gamma\in\{0,1\}$ being a parameter to be determined \cite{Chai2020}, and in this work, we consider $\gamma=0$ for the NACDE and $\gamma=1$ for the NSEs.

To discretize the propagation step (\ref{propagation}), we adopt the explicit two-level, three-point scheme \cite{Guo2001,Guo2018}:
\begin{equation}\label{discre-propagation-step}
	\begin{aligned}
		f_i(\bm{x},t+\Delta t)=p_0f^{\star}(\bm{x},t)+p_{-1}f^{\star}_i(\bm{x}-\bm{\lambda}_i\Delta t,t)+p_1f^{\star}_i(\bm{x}+\bm{\lambda}_i\Delta t,t),i=1,2,\ldots,q,\\
	\end{aligned}	
\end{equation}
where the free parameters $p_0,p_{-1}$ and $p_1$ should satisfy the following conditions:
\begin{align}\label{p}
	&p_0=1-b,p_{-1}=\frac{a+b}{2},p_1=\frac{b-a}{2},
\end{align}
with \begin{align}\label{c-a-relation}
	\bm{\lambda}_i=\lambda\bm{e}_i,	\lambda=\frac{\Delta x}{\Delta t}, a=\frac{|\bm{c}_i|}{|\bm{\lambda}_i|},c=a\lambda,(0<a\leq 1).
\end{align}
Substituting Eqs. (\ref{p}) and (\ref{c-a-relation}) into then Eq. (\ref{discre-propagation-step}) yields
\begin{align}\label{discre-propagation-step2}
	f_i(\bm{x},t+\Delta t)=&f_i^{\star}(\bm{x},t) -\frac{a}{2}\big[f^{\star}_i(\bm{x}+\bm{\lambda}_i\Delta t,t)-f^{\star}_i(\bm{x}-\bm{\lambda}_i\Delta t,t)\big]\notag\\
	& +\frac{b}{2}\big[f^{\star}_i(\bm{x}+\bm{\lambda}_i\Delta t,t)-2f^{\star}(\bm{x},t)+f^{\star}_i(\bm{x}-\bm{\lambda}_i\Delta t,t)\big],
\end{align} 
where $a$ and $b$ are considered as two free parameters. Here we would like to point out that Eq. (\ref{discre-propagation-step2}) can reduce to the propagation step of the standard LB model \cite{Guo2013} when $a=b=1$, and additionally, based on the stability structure analysis \cite{Zhao2020}, the two parameters $a$ and $b$ should satisfy the following condition:
\begin{equation}
	a^2\leq b\leq1.
\end{equation}

To simplify analysis and for the sake of brevity, the GPMRT-LB model composed of Eqs. (\ref{collison2}) and (\ref{discre-propagation-step2}) can also be expressed in a  matrix form,
\begin{subequations}
	\begin{align}
		\bm{m}^{\star,n}=&\big(\bm{I}_q-\bm{S}\big)\bm{m}^n+\bm{S}\bm{m}^{eq,n}+\Delta t\bm{\tilde{F}}^n,\label{matrix-vector-collision}\\
		\bm{m}^{n+1}(\bm{x})=&\bm{M}
		\Big(p_0\bm{M}^{-1}\bm{m}^{\star,n}(\bm{x})+p_{-1}\bm{M}^{-1}\bm{m}^{\star,n}(\bm{x}-\bm{\lambda}_i\Delta t)\notag\\
		&+p_1\bm{M}^{-1}\bm{m}^{\star,n}(\bm{x}+\bm{\lambda}_i\Delta t)\Big),\label{matrix-vector-propagation}
	\end{align}
\end{subequations}
where $\bm{I}_q\in R^{q\times q}$ is the identify matrix and 
\begin{subequations}
	\begin{align}
		&\Big\{\bm{m}^{\star,n};\bm{m}^{n};\bm{m}^{eq,n};\bm{\tilde{F}}^{n}\Big\}=\bm{M}\Big\{\bm{f}^{\star};\bm{f};\bm{f}^{eq};\bm{G}+\bm{F}+\frac{\Delta t}{2}\bm{\overline{D}}\bm{F}\Big\}(\bm{x},t_n),\\
		&\bm{m}^{\star,n}(\bm{x}\pm\bm{\lambda}_i\Delta t)=\bm{M}\Big(f_1(\bm{x}\pm\bm{\lambda}_1\Delta t,t_n),\ldots,f_q(\bm{x}\pm\bm{\lambda}_q\Delta t,t_n)\Big)^T, 
	\end{align} 
	\end{subequations}
here $\bm{\psi}:=\Big(\psi_1,\psi_2,\ldots,\psi_q\Big)^T$ with $\bm{\psi}$ representing $\{\bm{f}^{\star},\bm{f},\bm{f}^{eq},\bm{G},\bm{F}\}$ and $\bm{\overline{D}}:=$\textbf{diag}$\big(\overline{D}_1,\overline{D}_2,\ldots,\overline{D}_q \big)$. 

In the following, we assume that the matrices $\bm{M}$ and $\bm{S}$ are independent on the space and time. Then according to the space and time shift operators defined by Eqs.  (\ref{Space-oper}) and (\ref{Time-oper}), the GPMRT-LB model, i.e., Eqs. (\ref{matrix-vector-collision}) and (\ref{matrix-vector-propagation}),  can be rewritten as
\begin{equation}\label{GPMRT-LB-model}
	\big[T^1_{\Delta t}\bm{I}_q-\bm{A}\big]\bm{m}^{n}=\bm{B}\bm{m}^{eq,n}+\Delta t\bm{W}\bm{\tilde{F}}^n,
\end{equation}
where 
\begin{align} 
	&\bm{A}=\bm{W}\Big(\bm{I}_q-\bm{S}\Big),\bm{B}=\bm{WS},  
\end{align}
with 
\begin{subequations}\label{FD-scheme-quan}
	\begin{align}   
		&\bm{W}=\bm{M}\bm{\overline{T}}\bm{M}^{-1},\bm{\overline{T}}=\Big(p_0\bm{T}_0+p_{-1}\bm{T}_{-1}+p_1\bm{T}_1\Big),\\
		&\bm{T}_0=\bm{I}_q,\bm{T}_{-1}=\makebox{\textbf{diag}} \Big(T^{-\bm{e}_{1}}_{\Delta x},T^{-\bm{e}_{2}}_{\Delta x},\ldots,T^{-\bm{e}_{q}}_{\Delta x}\Big),\bm{T}_{1}=\makebox{\textbf{diag}}\Big(T^{\bm{e}_1}_{\Delta x},T^{\bm{e}_2}_{\Delta x},\ldots,T^{\bm{e}_{q}}_{\Delta x}\Big).
	\end{align}
\end{subequations}
Now we present a remark on the GPMRT-LB model, i.e., 
Eq. (\ref{GPMRT-LB-model}).\\
\textbf{Remark 1} The term $\overline{D}_iF_i$ in the collision step (\ref{collison2}) with $\gamma=1$ can be discretized by   an implicit difference scheme \cite{Chai2016},
\begin{equation}\label{implicit}
	\overline{D}_iF_i=\frac{F_i(\bm{x}+\bm{\lambda}_i\Delta t,t+\Delta t)-F_i(\bm{x},t)}{\Delta t},i=1,2,\ldots,q.
\end{equation} 
If we substitute Eq. (\ref{implicit}) into the GPMRT-LB model (\ref{GPMRT-LB-model}), one can obtain
\begin{equation}\label{GPMRT-LB-overD-D}
	\big[T^1_{\Delta t}\bm{I}_q-\bm{A}\big]\bm{\overline{m}}^{n}=\bm{B}\Big[\bm{m}^{eq,n}-\frac{\Delta t}{2}\bm{MF}^n\Big]+\Delta t\bm{W}\bm{\overline{F}}^n,
\end{equation}
which can be considered as a modified GPMRT-LB model with $\bm{\overline{m}}^n=\bm{m}^n-\frac{\Delta t}{2}\bm{MF}^n$ and $\bm{\overline{F}}^n=\bm{M}\big(\bm{G}^n+\bm{F}^n\big)$.
\subsection{Derivation of  the GPMFD scheme}\label{part-derivation-GPMFD}
In this part, we will provide some  details on how to derive the corresponding GPMFD scheme from the GPMRT-LB model (\ref{GPMRT-LB-model}). Without loss of generality, we assume that the first $N$ rows in $\bm{M}$ correspond to the $N$ conservative moments, and denote $i = 1, 2,\ldots, J$ as $i\in\{1\sim J\}$ for brevity. Now, we focus on two cases with $N = 1$  and $N > 1$, which are corresponding to the NACDE and NSEs  considered in this work.\\
\textbf{Proposition 1} For the case of $N=1$, the GPMRT-LB model (\ref{GPMRT-LB-model}) can be written as a multiple-level finite-difference scheme on the conservative moment $m_1$,
\begin{align}\label{GPMFD-N-eq1}
	m^{n+1}_{1}
	=-\sum_{k=1}^{q}\gamma_k \bm{m}^{n+k-q} _{1}
	+\sum_{k=1}^{q}\bigg[\sum_{l=1}^{k}
	\gamma_{q+1+l-k}\bm{A}^{l-1}\Big(\bm{B}\bm{m}^{eq|n-k+1}
	+\Delta t\bm{W\tilde{F}}^{n-k+1}\Big)\bigg]_{1},
\end{align} 
or 
\begin{align}\label{det-GPMFD-N-eq1}
	&\det{\big(T^1_{\Delta t}\bm{I}_q-\bm{A}\big)}m^n_1=\bigg[{\rm adj}\big(T^1_{\Delta t}\bm{I}_q-\bm{A}\big)\big(\bm{B}\bm{m}^{eq}+\Delta t\bm{W}\bm{\tilde{F}}\big)\bigg]_1,
\end{align}
where $\big(\gamma_k\big)_{k=1}^{q+1}$ are the coefficients of the monic characteristic polynomial $p_{\bm{A}}(x)=\sum_{k=1}^{q+1}\gamma_kx^{k-1}$ of matrix $\bm{A}$,  and ${\rm adj}(\cdot)$ represents the adjugate matrix.
\begin{proof} Let $n\in \mathcal{N}$, then for any $k\in \mathcal{N}$, applying Eq. (\ref{GPMRT-LB-model}) recursively gives
	\begin{align}\label{recurr-equation}
		T^1_{\Delta t}\bm{m}^{n}=\bm{m}^{n+1} &=\bm{A}\bm{m}^n +\bm{Bm}^{eq,n} +\Delta t\bm{W}\bm{\tilde{F}}^n\notag\\
		&
		=\bm{A}^2\bm{m}^{n-1}+\bm{AB}\bm{m}^{eq|n-1}+\bm{Bm}^{eq,n}+\Delta t\big(\bm{AW\tilde{F}}^{n-1}+\bm{W\tilde{F}}^n\big) \notag\\
		&=\ldots=\bm{A}^{k}\bm{m}^{n+1-k} +\sum_{l=0}^{k-1}
		\bm{A}^{l}\big(\bm{B}\bm{m}^{eq,n-l}+\Delta t\bm{W\tilde{F}}^{n-l}\big).  
	\end{align} 
	In order to fix the first term on the right-hand side of Eq. (\ref{recurr-equation}) regardless
	of the value of $k$, especially for $k\in\{1\sim (q+1)\}$, one can introduce $\tilde{n}:=n+1-k$ and obtain 
	\begin{align}\label{proof-1-eq1}
		&\bm{m}^{\tilde{n}+k}=\bm{A}^{k}\bm{m}^{\tilde{n}}
		+\sum_{l=0}^{k-1}
		\bm{A}^{l}\Big(\bm{B}\bm{m}^{eq,\tilde{n}+k-1-l}
		+\Delta t\bm{W\tilde{F}}^{\tilde{n}+k-1-l}\Big). 
	\end{align} 
	From Eq. (\ref{proof-1-eq1}), we also have
	\begin{align}
		&\sum_{k=1}^{q+1}\gamma_k\bm{m}^{\tilde{n}+k}
		= \bigg(\sum_{k=1}^{q+1}\gamma_k\bm{A}^{k}\bigg) 
		\bm{m}^{\tilde{n}}
		+\sum_{k=1}^{q+1}\gamma_k\bigg[\sum_{l=0}^{k-1}
		\bm{A}^{l}\Big(\bm{B}\bm{m}^{eq,\tilde{n}+k-1-l}
		+\Delta t\bm{W\tilde{F}}^{\tilde{n}+k-1-l}\Big)\bigg].
	\end{align} 
		Based on the characteristic polynomial $p_{\bm{A}}(x)=x^{-1}\sum_{k=1}^{q+1}\gamma_kx^k$ and the Cayley-Hamilton Theorem \cite{Bellotti2022}, we can obtain $\bm{A}\Big(\sum_{k=1}^{q+1}\gamma_k\bm{A}^{k-1}\Big)=\bm{0}$  and the following equation:
	\begin{align}\label{double-sum}
		&\bm{m}^{n+1}
		=-\sum_{k=1}^{q}\gamma_k\bm{m}^{n+k-q}
		+ \sum_{k=1}^{q+1}\gamma_k\bigg[\sum_{l=0}^{k-1}
		\bm{A}^{l}\Big(\bm{B}\bm{m}^{eq,n-q-1+k-l}
		+\Delta t\bm{W\tilde{F}}^{n-q-1+k-l}\Big)\bigg] ,
	\end{align}
	where $\tilde{n}+q=n$ and $\gamma_{q+1}=1$ have been used. Changing the indices of the last double summation  term on the right-hand side of Eq. (\ref{double-sum}) yields
	\begin{align}
		&\bm{m}^{n+1}
		=-\sum_{k=1}^{q}\gamma_k\bm{m}^{n+k-q}
		+\sum_{k=1}^{q}\bigg[\sum_{l=1}^{k}
		\gamma_{q+1+l-k}\bm{A}^{l-1}\Big(\bm{B}\bm{m}^{eq,n-k+1}
		+\Delta t\bm{W\tilde{F}}^{n-k+1}\Big)\bigg], 
	\end{align} 
	where the first row of the above equation is the macroscopic finite-difference scheme of the GPMRT-LB model (\ref{GPMRT-LB-model}) on the conservative moment $m_1$.
	
	To simplify the derivation of truncation errors and MEs in Section \ref{Truncation-Error-Modified-Equation}, an alternative form of the finite-difference scheme (\ref{GPMFD-N-eq1}) can be given by Eq. (\ref{det-GPMFD-N-eq1}). The reason is as follows. 
	
	Based on the relation between adjugate matrix and determinant, we have
	\begin{equation}\label{relation-det-adj}
		\big(	T^1_{\Delta t}\bm{I}_q-\bm{A}\big){\rm adj}\big(	T^1_{\Delta t}\bm{I}_q-\bm{A}\big)={\rm adj}\big(	T^1_{\Delta t}\bm{I}_q-\bm{A}\big)\big(	T^1_{\Delta t}\bm{I}_q-\bm{A}\big)=\det{\big(	T^1_{\Delta t}\bm{I}_q-\bm{A}\big)}\bm{I}_q,
	\end{equation}
	multiplying (\ref{GPMRT-LB-model}) by ${\rm adj}\big(	T^1_{\Delta t}\bm{I}_q-\bm{A}\big)$ yields  
	\begin{align}\label{det-GPMFD}
		&\det{\big(T^1_{\Delta t}\bm{I}_q-\bm{A}\big)}\bm{m}^n = {\rm adj}\big(T^1_{\Delta t}\bm{I}_q-\bm{A}\big)\big(\bm{B}\bm{m}^{eq,n}+\Delta t\bm{W}\bm{\tilde{F}}^n\big),
	\end{align}
	then selecting the first row of above Eq. (\ref{det-GPMFD}) gives Eq. (\ref{det-GPMFD-N-eq1}). In addition, according to the algebraic expression of the adjugate matrix and characteristic polynomial, i.e.,
	\begin{subequations}
		\begin{align}
			{\rm adj}\big(T^1_{\Delta t}\bm{I}_q-\bm{A}\big)&=\sum_{k=1}^{q}\bigg(\sum_{l=1}^{q-k+1}\gamma_{k+l}\bm{A}^{l-1}\bigg)T_{\Delta t}^{k-1}\bm{I}_q\notag\\
			&=T_{\Delta t}^{q-1}\sum_{k=1}^{q}\bigg(\sum_{l=1}^{q-k+1}\gamma_{k+l}\bm{A}^{l-1}\bigg)T_{\Delta t}^{k-q}\bm{I}_q,\label{adj-expression}\\
				 \det{\big(T^1_{\Delta t}\bm{I}_q-\bm{A}\big)}&=\sum_{k=1}^{q+1}\gamma_kT_{\Delta t}^{k-1}=T_{\Delta t}^{q-1}\sum_{k=1}^{q+1}\gamma_kT_{\Delta t}^{k-q},\label{det-expression}
		\end{align} 
		\end{subequations} 
	if we introduce $\tilde{k}:=q-k+1$, Eq. (\ref{adj-expression}) can be expressed as
	\begin{align}\label{adj-expression2}
		&{\rm adj}\big(T^1_{\Delta t}\bm{I}_q-\bm{A}\big)= T_{\Delta t}^{q-1}\sum_{k=1}^{q}\Big(\sum_{l=1}^{\tilde{k}}\gamma_{q+1+l-\tilde{k}}\bm{A}^{l-1}\Big)T_{\Delta t}^{1-\tilde{k}}\bm{I}_q.
	\end{align} 
	From Eqs. (\ref{GPMFD-N-eq1}), (\ref{det-expression}), and (\ref{adj-expression2}), it is easy to show that the two finite-difference schemes, i.e., Eqs. (\ref{GPMFD-N-eq1}) and (\ref{det-GPMFD-N-eq1}), on the conservative moment $m_1$ are identical.   \end{proof}

It should be noted that whether the first $N>1$ rows of Eq. (\ref{det-GPMFD}) can be considered as the finite-difference schemes on  $N>1$ conservative moments  is unclear. To this end, we first present a corollary of the Proposition 1. \\
\textbf{Corollary 1} Assuming that the relaxation matrix $\bm{S}$ of the GPMRT-LB model with $N\geq1$ conservative moments is a diagonal one  with diagonal elements $s_1,s_2,\ldots,s_j,\ldots,s_q$ located in the range $(0,2)$, then the $j_{th}$ $\big(j\in\{1\sim N\}\big)$ row of Eq. (\ref{det-GPMFD}) on the $j_{th}$ conservative moment does not depend on the relaxation parameter $s_j$.
\begin{proof} Assume that the force and source terms are independent of the relaxation matrix $\bm{S}$ [here we do not discretize $\overline{D}_jF_j^n$  in Eq. (\ref{det-GPMFD})], now we consider whether the two terms $\det{\big(T^1_{\Delta t}\bm{I}_q-\bm{A}\big)}m_j^n$ and $\Big[{\rm adj}\big(T^1_{\Delta t}\bm{I}_q-\bm{A}\big)\bm{B}\bm{m}^n\Big]_j$ in Eq. (\ref{det-GPMFD}) are independent on $s_j$  or not.  First, for any $j\in\{1\sim N\}$, the matrix $\bm{B}$ can be decomposed as
	\begin{equation}
		\bm{B}=\bm{B}|_{s_j=0}+\Big( \bm{B}-\bm{B}|_{s_j=0} \Big)=\bm{B}|_{s_j=0}+ \bm{B}_jI_j^T,
	\end{equation} where $\bm{B}_j$ and $I_j$ represent the $j_{th}$ column of matrix $\bm{B}$ and the identity matrix $\bm{I}_q$, respectively. Then based on the following relations:
	\begin{align}
		&\Big[{\rm adj}\big(T^1_{\Delta t}\bm{I}_q-\bm{A}\big)\bm{B}_jI_j^T\bm{m}^{eq,n}\Big]_{j}=I_j^T{\rm adj}\big(T^1_{\Delta t}\bm{I}_q-\bm{A}\big)\bm{B}_jm^{eq,n}_j; m^{eq,n}_j=m_j^n,
	\end{align}   the $j_{th}$ row of Eq. (\ref{det-GPMFD}) can be rewritten as 
	\begin{align}\label{above0}
		&\Big[\det{\big(T^1_{\Delta t}\bm{I}_q-\bm{A}\big)}- I_j^T{\rm adj}\big(T^1_{\Delta t}\bm{I}_q-\bm{A}\big)\bm{B}_j\Big]m_j^n=\Big[{\rm adj}\big(T^1_{\Delta t}\bm{I}_q-\bm{A}\big) \bm{B}|_{s_j=0}\bm{m}^{eq,n} \Big]_j,
	\end{align} 
	thanks to the fact  $\det{\big(T^1_{\Delta t}\bm{I}_q-\bm{A}\big)}=\det{\big(T^1_{\Delta t}\bm{I}_q-\bm{A}-\bm{B}_jI^T_j\big)}
	+I_j^T{\rm adj}\big(T^1_{\Delta t}\bm{I}_q-\bm{A}\big)\bm{B}_j$, Eq. (\ref{above0}) can be reformulated as
	\begin{align}\label{above1}
		&\bigg[\det{\Big(T^1_{\Delta t}\bm{I}_q-\big(\bm{A}+\bm{B}_jI^T_j\big)\Big)}\bigg]m_j^n=\Big[{\rm adj}\big(T^1_{\Delta t}\bm{I}_q-\bm{A}\big) \bm{B}|_{s_j=0}\bm{m}^{eq,n} \Big]_j.
	\end{align} 
	From the term on the right-hand side of Eq. (\ref{above1}), one can find  that the $j_{th}$ column of matrix $\bm{A}$ depends solely on $s_j$, thus the $j_{th}$ row of the adjugate matrix ${\rm adj}\big(T^1_{\Delta t}\bm{I}_q-\bm{A}\big)$ is independent on $s_j$. Similarly, from the term on the left-hand side of Eq. (\ref{above1}), one can derive the relation $\bm{A}+\bm{B}_jI^T_j=\bm{A}\Big|_{s_j=0}$, which does not depend on $s_j$. Therefore, the relaxation parameter $s_j$ 	does not influence the $j_{th}$ row of Eq. (\ref{det-GPMFD}) corresponding to the $j_{th}$ conservative moment.\end{proof} 
Regarding the conclusion of Corollary 1, we give a remark.\\
\textbf{Remark 2} It is clear that the $j_{th}$ $\big(j\in\{1\sim N\}\big)$ row of Eq. (\ref{det-GPMFD-N-eq1}) can only be independent of $s_j$, but one cannot show that it is  independent on the relaxation parameters associated with other conservative moments. Therefore, for $N>1$ conservative moments, the point that considering the first $N$ rows of Eq. (\ref{det-GPMFD}) as the finite-difference schemes on the $N$ conservative moments is inconsistent with the GPMRT-LB model. In the following, we will do some treatments on the matrix $\bm{A}$ to derive the macroscopic  finite-difference scheme of the GPMRT-LB model with $N>1$ conservative moments, as outlined in the Proposition 2. Additionally, we also consider a more general block-lower-triangular relaxation matrix $\bm{S}$ in the following  Corollary 2.\\
\textbf{Proposition 2} For $N\geq 1$ conservative moments, the GPMRT-LB model (\ref{GPMRT-LB-model}) corresponds to a multiple-level finite-difference scheme on the $j_{th}$ conservative moment $m_j$ $\big(j\in\{1\sim N\}\big)$,

\begin{align}\label{GPMFD-N-geq1}
	m_j^{n+1}
	=&-\sum_{k=1}^{q+1-N}\gamma_{j,k}m_j^{n+N+k-1-q}+\sum_{k=1}^{q+1-N}\Big[\sum_{l=1}^{k}\gamma_{j,q+2-N+l-k}
	\tilde{\bm{A}}_j^{l-1}\overline{\bm{A}}_j\bm{m}^{n-k+1}
	\Big]_j\notag\\
	&+\sum_{k=1}^{q+1-N}\bigg[\sum_{l=1}^{k}\gamma_{j,q+2-N+l-k}
	\tilde{\bm{A}}_j^{l-1}\Big(\bm{B}\bm{m}^{eq|n-k+1}
	+\Delta t\bm{W\tilde{F}}^{n-k+1}\Big)\bigg]_j,
\end{align} 
or 
\begin{align}\label{det-GPMFD-N-geq1}
	\det{\big(T^1_{\Delta t}\bm{I}_q-\bm{\tilde{A}}_j\big)}m^n_j=&\Big[{\rm adj}\big(T^1_{\Delta t}\bm{I}_q-\bm{\tilde{A}}_j\big)\bm{\overline{A}}_j\bm{m}^n\Big]_j\notag\\
	&+\Big[{\rm adj}\big(T^1_{\Delta t}\bm{I}_q-\bm{\tilde{A}}_j\big)\big(\bm{B}\bm{m}^{eq,n}+\Delta t\bm{W}\bm{\tilde{F}}^n\big)\Big]_j,
\end{align} 
where 
\begin{subequations}
	\begin{align}
		&\bm{\tilde{A}}_j=\bm{A}\bm{P}_j, \bm{P}_j:=\sum_{l=j,N+1,\ldots,q}I_lI_l^T, \label{A*P}\\
		&\bm{\overline{A}}_j=\bm{A}-\bm{\tilde{A}}_j\label{A*(I-P)},
	\end{align}
\end{subequations}
	and  $\big(\gamma_{j,k}\big)_{k=1}^{q+2-N}$ are the coefficients of the monic characteristic polynomial $p_{\bm{\tilde{A}}_j}(x)=\sum_{k=1}^{q+2-N}\gamma_{j,k}x^{k+N-2}$ of matrix $\bm{\tilde{A}}_j$.  The proof is similar to that of Proposition 1, and the details are not presented here. In addition, we would also like to point out that the finite-difference scheme (\ref{det-GPMFD-N-geq1}) has the following Corollary 2.\\
\textbf{Corollary 2} If the relaxation matrix  $\bm{S}$ of the GPMRT-LB model with $N\geq 1$ conservative moments is a block-lower-triangular form with the diagonal elements located in range $(0,2)$ [see Eq. (\ref{S})], then  for any $j\in\{1\sim N\}$, the finite-difference scheme (\ref{det-GPMFD-N-geq1}) on the $j_{th}$ conservative moment is independent on the lower triangular relaxation parameters $s_{il}$ $\big(l\in\{1\sim N\};i\in\{l\sim q\}\big)$,\\ 
\begin{equation}\label{S}
	\begin{aligned}
		&\left(\begin{matrix}
			s_{11}&0&\ldots&0&\ldots&\ldots&0\\
			s_{21}&s_{22}&0&0&\ldots&\ldots&0\\
			s_{31}&s_{32}&s_{33}&0&\ldots&\ldots&0\\
			\vdots&\ddots&\vdots&\ddots&\vdots&\ddots&\vdots\\
			s_{N1}&s_{N2}&s_{N3}&\ldots&s_{NN}&\ldots&0\\
			\vdots&\ddots&\vdots&\ddots&\vdots&\multicolumn{2}{c}{\multirow{2}{*}{$ \bm{S_r}$}}\\
			s_{q1}&s_{q2}&s_{q3}&\ldots &s_{qN}&&\\
		\end{matrix}\right),
	\end{aligned}
\end{equation}
where $\bm{S}_r\in R^{(q-N)\times (q-N)}$ is a block-lower-triangular relaxation matrix.
\begin{proof} Similar to the proof of Corollary 1, for any fixed $j\in\{1\sim N\}$, the matrix $\bm{B}$ can be decomposed as
	\begin{equation}
		\begin{aligned}
			&\bm{B}=\bm{B}\Big|_{s_{kj}=0,k\geq j}+
			\bm{B}_jI^T_j,
		\end{aligned}
	\end{equation}
	then with the help of the following relations:
	\begin{equation}
		\begin{aligned}
			&	\Big[{\rm adj}\big(T^1_{\Delta t}\bm{I}_q-\bm{\tilde{A}}_j\big)\bm{B}_jI_j^T\bm{m}^{eq,n}\Big]_j=I_j^T{\rm adj}\big(T^1_{\Delta t}\bm{I}_q-\bm{\tilde{A}}_j\big)\bm{B}_jm^{eq,n}_j; m^n_j=m^{eq,n}_j,
		\end{aligned}
	\end{equation} the finite-difference scheme (\ref{det-GPMFD-N-geq1}) can be rewritten as
	\begin{equation} \label{above2}
		\begin{aligned}
			&\Big[\det{\big(T^1_{\Delta t}\bm{I}_q-\bm{\tilde{A}}_j\big)}- I_j^T{\rm adj}\big(T^1_{\Delta t}\bm{I}_q-\bm{\tilde{A}}_j\big)\bm{B}_j\Big]m^n_j\\
			&\qquad\qquad=\Big[{\rm adj}\big(T^1_{\Delta t}\bm{I}_q-\bm{\tilde{A}}_j\big) \bm{\overline{A}}_j\bm{m}^n \Big]_j
			+\Big[{\rm adj}\big(T^1_{\Delta t}\bm{I}_q-\bm{\tilde{A}}_j\big) \bm{B}\Big|_{s_{kj}=0,k\geq j}\bm{m}^{eq,n} \Big]_j.
		\end{aligned}
	\end{equation}
	In the following, the proof is divided into four steps.\\
	\textbf{Step 1:} Regarding the term on the left-hand side of Eq. (\ref{above2}), the definition of $\bm{\tilde{A}}_j$ (\ref{A*P}) gives
\begin{subequations}
	\begin{align}
		&\Big[\bm{\tilde{A}}_j\Big]_{il}=\Big[\bm{A}\Big]_{il},i\in\{1\sim q\};l=\{j\}\cup\{(N+1)\sim q\} , \\	&\Big[\bm{\tilde{A}}_j\Big]_{il}=0,i\in\{1\sim q\};l\in\{1\sim N\}/\{j\},
	\end{align}
\end{subequations}	
	from which one can find  that  $\bm{\tilde{A}}_j$ and $\bm{B}_j$ are independent on $s_{il}\big(l\in\{1\sim N\}/\{j\};i\in\{l\sim q\}\big)$, this also means that the term on the left-hand side of Eq. (\ref{above2}) does not depend on the $l_{th}$ $\big(l\in\{1\sim N\}/\{j\}\big)$  column  of matrix $\bm{S}$.\\
	\textbf{Step 2:} For the terms on the right-hand side of Eq. (\ref{above2}), the definitions of $\bm{\overline{A}}_j$ (\ref{A*(I-P)}) and $\bm{B}\Big|_{s_{kj}=0,k\geq j}$ show that $\Big[\bm{\overline{A}}_j\Big]_{pr}$ and $\Big[\bm{B}\Big|_{s_{kj}=0,k\geq j}\Big]_{pr}$ $\big(p\in\{1\sim q\};r\in\{N+1\sim q\}\big)$, are independent on $s_{il}\big(l\in\{1\sim N\}/\{j\};i\in\{l\sim q\}\big)$. Due to $m^n_l=m^{eq,n}_l$ for $l\in\{1\sim N\}/\{j\}$,  $\Big[\bm{\overline{A}}_j\Big]_{il}$ and $\Big[\bm{B}\Big|_{s_{kj}=0,k\geq j}\Big]_{il}$ can be considered  together for $i\in\{1\sim q\}$ and $l\in\{1\sim N\}/\{j\}$, 
	\begin{align}
		&\Big[\bm{\overline{A}}_j\Big]_{il}+\Big[\bm{B}\Big|_{s_{kj}=0,k\geq j}\Big]_{il}\notag\\ &=\Big[\bm{\overline{A}}_j+\bm{B}\Big|_{s_{kj}=0,k\geq j}\Big]_{il}=
		\Big[\bm{A}+\bm{B}\Big]_{il}=\Big[\bm{W}\Big]_{il}=\Big[\bm{A}\Big|_{s_{kl}=0,k\geq l}\Big]_{il}=\Big[\bm{\overline{A}}_j\Big|_{s_{kl}=0,k\geq l}\Big]_{il},
	\end{align} 
	thus  the terms on the right-hand side of Eq. (\ref{above2}) are independent on the $l_{th}$ $\big(l\in\{1\sim N\}/\{j\}\big)$ column of $\bm{S}$.\\
	\textbf{Step 3:} Now we consider whether Eq. (\ref{above2}) depends on the $j_{th}$ column of matrix $\bm{S}$ or not. After some algebraic manipulations, Eq. (\ref{above2}) can also be   reformulated as 
	\begin{align}\label{above3}
		&\bigg[\det{\Big(T^1_{\Delta t}\bm{I}_q-\big(\bm{\tilde{A}}_j+ \bm{B}_j^T I_j\big)\Big)}\bigg]m^n_j\notag\\
		&\qquad=\Big[{\rm adj}\big(T^1_{\Delta t}\bm{I}_q-\bm{\tilde{A}}_j\big) \bm{\overline{A}}_j\bm{m}^n \Big]_j
		+\Big[{\rm adj}\big(T^1_{\Delta t}\bm{I}_q-\bm{\tilde{A}}_j\big) \bm{B}\Big|_{s_{kj}=0,k\geq j}\bm{m}^{eq,n} \Big]_j.
	\end{align} 
	From the terms on the the right-hand side of   Eq. (\ref{above3}), one can see that the $j_{th}$ column of $\bm{\tilde{A}}_j$ only depends on the $j_{th}$ column of matrix $\bm{S}$. Therefore, the $j_{th}$ row of ${\rm adj}\big(T^1_{\Delta t}\bm{I}_q-\bm{\tilde{A}}_j\big)$ does not depend on the $j_{th}$ column of matrix $\bm{S}$.\\
	\textbf{Step 4:} From the term on the left-hand side of   Eq. (\ref{above3}), for any $i\in\{1\sim q\}$, one can obtain 
	\begin{equation}
	\Big[\bm{\tilde{A}}_j+\bm{B}_jI_j^T\Big]_{il}= \begin{cases}
		\Big[\bm{A} \Big]_{il},&l\in\{(N+1)\sim q\}, \\ 
		0, & l\in\{1\sim N\}/\{j\}, \\
	\Big[\bm{A} +\bm{B}\Big]_{il}=\Big[\bm{W}\Big]_{il}=\bigg[\bm{\tilde{A}}_j\Big|_{s_{kj}=0,k\geq j}\bigg]_{il},&l=j, \\		
		\end{cases}
	\end{equation}
	thus the term on the left-hand side of Eq. (\ref{above3})  does depend on the $j_{th}$ column of matrix $\bm{S}$. Due to the equivalence of finite-difference scheme (\ref{det-GPMFD-N-geq1}), Eqs. (\ref{above2}) and (\ref{above3}), one can prove the Corollary 2.\end{proof}

We now conduct some discussion on the Corollary 2.\\
\textbf{Remark 3} According to the Remark 1 and Proposition 2, if the parameter $\gamma$ in $\overline{\bm{D}}$ is equal to one, the macroscopic finite-difference scheme on the $j_{th}$ $\big(j\in\{1\sim N\}\big)$ conservative moment of the GPMRT-LB model (\ref{GPMRT-LB-overD-D}) is given by 
\begin{align}\label{GPMFD-overD-D}
	&\det{\Big(T^1_{\Delta t}\bm{I}_q-\bm{\tilde{A}}_j\Big)}\bm{\overline{m}}_j^n=\bigg[{\rm adj}\Big(T^1_{\Delta t}\bm{I}_q-\bm{\tilde{A}}\Big)\bm{\overline{A}}_j\bm{\overline{m}}^n\bigg]_j\notag\\
	&\qquad+\bigg[{\rm adj}\Big(T^1_{\Delta t}\bm{I}_q-\bm{\tilde{A}}_j\Big)\bm{B}\Big(\bm{m}^{eq,n}-\frac{\Delta t}{2}\bm{MF}\Big)\bigg]_j+\Delta t\bigg[{\rm adj}\Big(T^1_{\Delta t}\bm{I}_q-\bm{\tilde{A}}_j\Big)\bm{W}\bm{\overline{F}}^n\bigg]_j,
\end{align}
while it is unclear whether the finite-difference scheme (\ref{GPMFD-overD-D}) has the Corollary 2 or not. We will pay attention to  this problem in the following part. It is clear that the last term on the right-hand side of Eq. (\ref{GPMFD-overD-D}) does not depend on the relaxation matrix $\bm{S}$, thus  we only need to consider the following scheme:
\begin{align}
	\det{\Big(T^1_{\Delta t}\bm{I}_q-\bm{\tilde{A}}_j\Big)}\bm{\overline{m}}^n_j=&\bigg[{\rm adj}\Big(T^1_{\Delta t}\bm{I}_q-\bm{\tilde{A}}_j\Big)\bm{\overline{A}}_j\bm{\overline{m}}^n\bigg]_j \notag \\
	&+\bigg[{\rm adj}\Big(T^1_{\Delta t}\bm{I}_q-\bm{\tilde{A}}_j\Big)\bm{B}\Big(\bm{m}^{eq,n}-\frac{\Delta t}{2}\bm{MF}\Big)\bigg]_j.
\end{align}
If we write the matrix $\bm{B}$ as
\begin{align}
	&\bm{B}=\bm{B}\Big|_{s_{kj}=0,k\geq j}+
	\bm{B}_jI^T_j,
\end{align}
and according to the relation $\bm{\overline{m}}_j=\bm{m}^{eq}-\Delta t/2\bm{MF}$,  one can obtain
\begin{align}
	&\Big[\det{\big(T^1_{\Delta t}\bm{I}_q-\bm{\tilde{A}}_j\big)}- I_j^T{\rm adj}\big(T^1_{\Delta t}\bm{I}_q-\bm{\tilde{A}}_j\big)\bm{B}_j\Big]\bm{\overline{m}}^n_j=\Big[{\rm adj}\big(T^1_{\Delta t}\bm{I}_q-\bm{\tilde{A}}_j\big) \bm{\overline{A}}_j\bm{\overline{m}}^n \Big]_j\notag\\
	&\qquad
	+\bigg[{\rm adj}\big(T^1_{\Delta t}\bm{I}_q-\bm{\tilde{A}}_j\big)\bm{B}\Big|_{s_{kj}=0,k\geq j}\Big( \bm{m}^{eq,n}-\frac{\Delta t}{2}\bm{MF}\Big) \bigg]_j,
\end{align}
which is similar to Eq. (\ref{above2}), and the detailed proof is not presented here. From above discussion,  it can be concluded that the finite-difference scheme (\ref{GPMFD-overD-D}) also has the Corollary 2, like the finite-difference scheme (\ref{det-GPMFD-N-geq1}).\\
\textbf{Remark 4} Based on the Corollary 2,  the finite-difference scheme (\ref{det-GPMFD-N-geq1}) on the $j_{th}$ $\big(j\in\{1\sim N\}\big)$ conservative moment  is consistent with the GPMRT-LB model, while Eq. (\ref{det-GPMFD}) is not. It should also be noted that  the lower triangular elements $s_{il}$ $\big(l\in\{1\sim N\};i\in\{l\sim q\}\big)$ in the relaxation matrix $\bm{S}$ (\ref{S}) do not affect the forms of difference schemes (\ref{GPMFD-N-geq1}) and (\ref{det-GPMFD-N-geq1}). To simplify following   analysis on the  truncation errors and MEs in Section \ref{Truncation-Error-Modified-Equation}, we assume that the diagonal relaxation parameters corresponding to the conservative moments in the relaxation matrix $\bm{S}$  (\ref{S}) are equal to one, and the non-diagonal relaxation parameters associated with the conservative moments are equal to zero, in this case, the finite-difference schemes (\ref{GPMFD-N-geq1}) and (\ref{det-GPMFD-N-geq1}) can be further simplified, see the following Proposition 3 for details.\\
\textbf{Proposition 3  } For a given GPMRT-LB model ($\bm{M},\bm{S}$) with $N\geq 1$ conservative moments, let $\bm{\hat{S}}=$\textbf{diag}$\Big(\bm{I}_N,\bm{S}_r\Big)$ with $\bm{S}_r\in R^{(q-N)\times (q-N)}$ representing the matrix consisting of the $(N+1)_{th}$ to $q_{th}$ rows and columns of matrix $\bm{S}$, then the finite-difference schemes (\ref{GPMFD-N-geq1}) and (\ref{det-GPMFD-N-geq1}) of the GPMRT-LB model (\ref{GPMRT-LB-model}) on the $j_{th}$ $\big(j\in\{1\sim N\}\big)$ conservative moment can be simplified as 
\begin{align}\label{GPMFD-N-geq1-brevity}
	m^{n+1}_{j }
	=&-\sum_{k=N}^{q}\gamma_{k}m^{n+k-q}_{j } \notag\\
	&+\sum_{k=N}^{q}\bigg[\sum_{l=N}^{k}
	\gamma_{q+1+l-k}\bm{ A} ^{l-N}\Big(\bm{B}\bm{m}^{eq|n-k+N}
	+\Delta t\bm{W}\bm{\tilde{F}}^{n-k+N}\Big)\bigg]_{j },
\end{align} 
and 
\begin{align}\label{det-GPMFD-N-geq1-brevity}
	&\det{\big(T^1_{\Delta t}\bm{I}_q-\bm{A}\big)}m_j^n=\bigg[{\rm adj}\big(T^1_{\Delta t}\bm{I}_q-\bm{A}\big)\big(\bm{B}\bm{m}^{eq,n}+\Delta  t\bm{W}\bm{\tilde{F}}^n\big)\bigg]_j, 
\end{align} 
where $\bm{A}=\bm{W}\big(\bm{I}_q-\bm{\hat{S}}\big)$, $\bm{B}=\bm{W}\bm{\hat{S}}$. In addition,   Eq. (\ref{det-GPMFD-N-geq1-brevity}) has the same  form as the $j_{th}$ row of Eq. (\ref{det-GPMFD}).  
\begin{proof} Based on the form of relaxation $\bm{\hat{S}}$ and the definition of matrix $\bm{P}_j$ (\ref{A*P}), one can obtain 
	\begin{align}
		&	\bm{\tilde{A}}_j=\bm{A}\bm{P}_j=\bm{W}(\bm{I}_q-\bm{\hat{S}})=\bm{A},
	\end{align}  which means $\bm{\overline{A}}_j=\bm{A}-\bm{\tilde{A}}_j =\bm{0}$, thus the finite-difference schemes (\ref{GPMFD-N-geq1}) and  (\ref{det-GPMFD-N-geq1}) can be simplified by 
	\begin{align}\label{gammjk}
		m^{n+1}_{j }
		=&-\sum_{k=1}^{q+1-N}\gamma_{j,k}m^{n+N+k-1-q}_{j }\notag\\
		&+\sum_{k=1}^{q+1-N}\bigg[\sum_{l=1}^{k}
		\gamma_{j,q+2-N+l-k}\bm{{A}}^{l-1}\Big(\bm{B}\bm{m}^{eq|n-k+1}
		+\Delta t\bm{W}\bm{\tilde{F}}^{n-k+1}\Big)\bigg]_{j }.
	\end{align}
	and
	\begin{align}\label{det-GPMFD-N-geq1-brevity-proof}
		&\det{\big(T^1_{\Delta t}\bm{I}_q-\bm{A}\big)}m_j^n=\bigg[{\rm adj}\big(T^1_{\Delta t}\bm{I}_q-\bm{A}\big)\big(\bm{B}\bm{m}^{eq,n}+\Delta  t\bm{W}\bm{\tilde{F}}^n\big)\bigg]_j.
	\end{align} 
	Then one can obtain the finite-difference schemes (\ref{GPMFD-N-geq1-brevity}) and (\ref{det-GPMFD-N-geq1-brevity})  through rearranging the Eqs. (\ref{gammjk}) and (\ref{det-GPMFD-N-geq1-brevity-proof}) and with the aid of following results,
	\begin{align}
		&\gamma_k=0,k\in\{1\sim (N-1)\}; \gamma_k=\gamma_{j,k+1-N},k\in\{N\sim (q+1)\},
	\end{align}
	where the relation $p_{\bm{A}}(x)=x^{-1}\sum_{k=1}^{q+1}\gamma_kx^k=p_{\bm{\tilde{A}}_j}(x)=x^{N-2}\sum_{k=1}^{q+2-N}\gamma_{j,k}x^{k}$ has been used.\end{proof} 
Now, we further discuss the issue of the equivalence between two GPMRT-LB models. On one hand, the two GPMRT-LB models $(\bm{M},\bm{S})$ and $(\bm{M},\bm{S}_1)$ can be considered equivalently if the relaxation parameters associated with non-conservative moments in the relaxation matrices $\bm{S}$ and $\bm{S}_1$ are identical regardless of whether the relaxation parameters associated with conservative moments are same or not, and the finite-difference scheme (\ref{det-GPMFD-N-geq1}) also has this feature  (see the Corollary 2). 
On the other hand, the two GPMRT-LB models $(\bm{M},\bm{S})$ and $(\bm{M}_1,\bm{S}_1)$ are also equivalent if the following relations hold:
\begin{subequations}\label{equivalem-mrt-model}
	\begin{align}
		&\bm{M}_1=\bm{N}\bm{M},\bm{S}_1=\bm{N}\bm{S}\bm{N}^{-1},\label{MM1-equivalent}\\
		&\big[\bm{M}_1\big]_{il}=\big[\bm{M}\big]_{il},i\in\{1\sim N\};l\in\{1\sim q\},\label{moment-condition-equivalence}
	\end{align}
\end{subequations}where $\bm{N}$ is an invertible block-lower-triangular matrix. This means that the finite-difference scheme of the GPMRT-LB model ($\bm{M}_1,\bm{S}_1$) should be the same as that of the GPMRT-LB model ($\bm{M},\bm{S}$). In the following, we first show  that the first $N$ rows of Eq. (\ref{det-GPMFD})  corresponding   to the two equivalent GPMRT-LB models ($\bm{M},\bm{S}$) and ($\bm{M}_1,\bm{S}_1$) are identical, and then   present another form of the finite-difference scheme from the GPMRT-LB model $(\bm{M}_1,\bm{S}_1)$, which is identical to the schemes (\ref{GPMFD-N-geq1}) and (\ref{det-GPMFD-N-geq1}).\\
\textbf{Theorem 1} The first $N$ rows of Eq. (\ref{det-GPMFD}) corresponding   to the two equivalent GPMRT-LB models ($\bm{M},\bm{S}$) and ($\bm{M}_1,\bm{S}_1$) satisfying Eq. (\ref{equivalem-mrt-model}) are totally identical.
\begin{proof} According to the relation (\ref{MM1-equivalent}), the matrices $\bm{A} $ and $\bm{A}^1 $ satisfy 
	\begin{align}
		\bm{A}^1=\bm{N}\bm{A}\bm{N}^{-1},
	\end{align}
	where $\bm{A}^1=\bm{W}^1(\bm{I}-\bm{S}_1)$ with $\bm{W}^1=\bm{M}_1\bm{\overline{T}}\bm{M}^{-1}_1$. Thanks to the invertible block-lower-triangular matrix $\bm{N}$, one can obtain 
	\begin{align}\label{detm1m2}
		\det{\big[T^1_{\Delta t} \bm{I}_q-\bm{A}\big]} = \det{\big[T^1_{\Delta t} \bm{I}_q- \bm{A}^1\big]}.
	\end{align}
	Based on  the relation (\ref{moment-condition-equivalence}), we have
	\begin{align}
		&\big[\bm{A}^1\bm{B}^1\bm{M}_1\big]_{il}=\big[\bm{A}\bm{B}\bm{M}\big]_{il};\big[\bm{A}^1\bm{W}^1\bm{M}_1\big]_{il}=\big[\bm{A}\bm{W}\bm{M}\big]_{il} ,
	\end{align} 
	where $\bm{B}^1=\bm{W}^1\bm{S}_1$. Thus, from the algebraic expression of adjugate matrix ${\rm adj}\big(T^1_{\Delta t}\bm{I}_q-\bm{A}\big)$ (\ref{adj-expression}) and for any $j\in\{1\sim N\}$, it is easy to prove  
	\begin{subequations}\label{adjm1m2}
		\begin{align} 
			&\Big[{\rm adj}\big(T^1_{\Delta t}\bm{I}_q-\bm{A}^1\big)\bm{B}^1\bm{M}_1\bm{f}^{eq,n}\Big]_j=\Big[{\rm adj}\big(T^1_{\Delta t}\bm{I}_q-\bm{A}\big)\bm{B}\bm{M}\bm{f}^{eq,n}\Big]_j,\\
			&\Big[{\rm adj}\big(T^1_{\Delta t}\bm{I}_q-\bm{A}^1\big)\bm{W}^1\bm{M}_1\big(\bm{F}^n+\bm{G}^n+\frac{\Delta t}{2}\bm{D}\bm{F}^n\big)\Big]_j\notag\\
			&\qquad=\Big[{\rm adj}\big(T^1_{\Delta t}\bm{I}_q-\bm{A}\big)\bm{W}\bm{M}\big(\bm{F}^n+\bm{G}^n+\frac{\Delta t}{2}\bm{D}\bm{F}^n\big)\Big]_j.
		\end{align}
	\end{subequations}
	According to the relations (\ref{detm1m2}) and (\ref{adjm1m2}), one can prove the Theorem 1.
\end{proof} 
Based on the Theorem 1, we can also conclude that the analysis on the truncation errors and MEs shown below will remain identical whether the transform matrix $\bm{M}$ is independent of the parameter $c$ [see Eq. (\ref{c-a-relation})] or not. The reason is provided in the following Remark 5. 

\noindent\textbf{Remark 5} Considering the following relation between the two GPMRT-LB models $(\bm{M}^c,\bm{S}^c)$  and $(\bm{M}^{o},\bm{S}^o)$:
\begin{subequations}\label{equivalem-cd-mrt-model}
	\begin{align}
		&\bm{M}^c=\bm{C}_d\bm{M}^o,\bm{S}^c=\bm{C}_d\bm{S}^o\bm{C}_c^{-1},\\
		&\big[\bm{M}^c\big]_{il}=\big[\bm{C}_d\big]_{ii}\big[\bm{M}^o\big]_{il},i\in\{1\sim N\};l\in\{1\sim q\},
	\end{align}
\end{subequations}
where $\bm{C}_d$ is an invertible diagonal matrix  associated with the parameter $c$. Similar to the proof in Theorem 1, one can obtain
\begin{align}\label{cd-detm1m2}
	\det{\big[T^1_{\Delta t} \bm{I}_q-\bm{A}^c\big]}\bm{M}^c\bm{f} = \bm{C}_d\det{\big[T^1_{\Delta t} \bm{I}_q- \bm{A}\big]}\bm{M}^o\bm{f},
\end{align}
and
\begin{subequations}\label{cd-adjm1m2}
\begin{align}
	&\Big[{\rm adj}\big(T^1_{\Delta t}\bm{I}_q-\bm{A}^c\big)\bm{B}^c\bm{M}^c\bm{f}^{eq,n}\Big]_j=\bm{C}_d\Big[{\rm adj}\big(T^1_{\Delta t}\bm{I}_q-\bm{A}\big)\bm{B}\bm{M}^o\bm{f}^{eq,n}\Big]_j, \\
	&\Big[{\rm adj}\big(T^1_{\Delta t}\bm{I}_q-\bm{A}^c\big)\bm{W}^c\bm{M}^c\big(\bm{F}^n+\bm{G}^n+\frac{\Delta t}{2}\bm{D}\bm{F}^n\big)\Big]_j\notag\\
	&\qquad=\bm{C}_d\Big[{\rm adj}\big(T^1_{\Delta t}\bm{I}_q-\bm{A}\big)\bm{W}\bm{M}^o\big(\bm{F}^n+\bm{G}^n+\frac{\Delta t}{2}\bm{D}\bm{F}^n\big)\Big]_j.
\end{align}
\end{subequations}

It is obvious that the first $N$ rows of Eq. (\ref{det-GPMFD}) corresponding respectively to the two GPMRT-LB models ($\bm{M}^c,\bm{S}^c$) and ($\bm{M}^o,\bm{S}^o$) satisfying Eq. (\ref{equivalem-cd-mrt-model}) only differ in the constant matrix $\bm{C}_d$, which has no influence on the truncation errors and MEs analysis in Section \ref{Truncation-Error-Modified-Equation}.

\noindent \textbf{Theorem 2} The finite-difference scheme (\ref{det-GPMFD-N-geq1}) corresponding to the GPMRT-LB model $(\bm{M},\bm{S})$ can be rewritten as   
\begin{align}\label{GPMFD-M1-S1}
	&\det{\Big(T^1_{\Delta t}\bm{I}_q-\bm{\tilde{A}}_j^1\Big)}\bm{m}_j^{n,(1)}=\bigg[{\rm adj}\Big(T^1_{\Delta t}\bm{I}_q-\bm{\tilde{A}}_j^1\Big)\bm{\overline{A}}_j^1\bm{m}^{n,(1)}\bigg]_j\notag\\
	&\qquad+\bigg[{\rm adj}\Big(T^1_{\Delta t}\bm{I}_q-\bm{\tilde{A}}_j^1\Big)\bm{B}^1\bm{m}^{eq,n,(1)}\bigg]_j+\Delta t\bigg[{\rm adj}\Big(T^1_{\Delta t}\bm{I}_q-\bm{\tilde{A}}_j^1 \Big)\bm{W} \bm{\tilde{F}} ^{n,(1)}\bigg]_j,
\end{align} 
where 
\begin{subequations}\label{theorem2-parameter}
\begin{align} 
	&\bm{A}^1=\bm{W}^1\big(\bm{I}_q-\bm{S}_1\big),\bm{B}^1=\bm{W}^1\bm{S}_1, \bm{\tilde{A}}_j^1=\bm{A}^1\bm{P}_j^1,\bm{\overline{A}}_j^1=\bm{A}^1\bm{\overline{P}}_j ,\\
	&\bm{m}^{n,(1)}=\bm{M}_1\bm{m}^n,\bm{m}^{eq,n,(1)}=\bm{M}_1\bm{m}^{eq,n},\bm{\tilde{F}}^{n,(1)}=\bm{M}_1\bm{\tilde{F}}^n,
\end{align} 
\end{subequations}
with the matrices $\bm{M}_1$ and $\bm{S}_1$ satisfying Eq. (\ref{equivalem-mrt-model}) and
\begin{align} \label{p1j} 
	&\bm{W}^1:=\bm{M}_1\bm{\overline{T}}\bm{M}^{-1}_1 ,\bm{P}_j^1:=\bm{N}\bm{P}_j\bm{N}^{-1}, \bm{\overline{P}}_j:=\bm{I}_q-\bm{\tilde{P}}_j^1.
\end{align}
The proof is similar to Theorem 1 and the details are not shown here. 

We now give some remarks on  the conclusions in Theorems 1 and 2.\\
\textbf{Remark 6} For two equivalent GPMRT-LB models ($\bm{M},\bm{S}$) and $(\bm{M}_1,\bm{S}_1)$ in the Theorem 1, if we further consider the GPMRT-LB model $(\bm{M}_1,\bm{\hat{S}_1})$ with
\begin{align}\label{hats1s2}
	\bm{\hat{S}}_1=&\bm{N}\bm{\hat{S}}\bm{N}^{-1}=\left(\begin{matrix}
		\bm{N}_1&\bm{0}\\
		\bm{N}_2&\bm{N}_3\\
	\end{matrix}\right)\left(\begin{matrix}
		\bm{I}_N&\bm{0}\\
		\bm{0}&\bm{S}_r\\
	\end{matrix}\right)\left(\begin{matrix}
		\bm{N}^{-1}_1&\bm{0}\\
		-\bm{N}_3^{-1}\bm{N}_2\bm{N}_1^{-1}&\bm{N}_3^{-1}\\
	\end{matrix}\right)\notag\\
	=&\left(\begin{matrix}
		\bm{I}_N&\bm{0}\\
		\bm{N}_2\bm{N}_1^{-1}-\bm{N}_3\bm{S}_r\bm{N}_3^{-1}\bm{N}_2\bm{N}_1^{-1}&\bm{N}_3\bm{S}_r\bm{N}_3^{-1}\\
	\end{matrix}\right),
\end{align}
where matrix $\bm{N}_1\in R^{N\times N}$, $\bm{N}_2\in R^{ (q-N)\times N}$, and $\bm{N}_3\in R^{(q-N)\times (q-N)}$ are the submatrices of the block-lower-triangular matrix $\bm{N}$, and the matrix  $\bm{\hat{S}}$ is defined as that in Proposion 3, one can obtain 
\begin{align}\label{hatS}
	\bm{\hat{S}}=\bm{S}\bm{\tilde{P}}_j+\bm{\overline{P}}_j+\big(\bm{I}_q-\bm{S}\big){I}_j{I}_j^T.
\end{align}
Subsituting above Eq. (\ref{hatS}) into Eq. (\ref{hats1s2}) yields
\begin{align}
	\bm{\hat{S}}_1=\bm{S}_1\bm{\tilde{P}}_j^1+\bm{\overline{P}}^1_j+\big(\bm{I}_q-\bm{S}_1\big)\bm{N}{I}_jI_j^T\bm{N}^{-1},
\end{align}
where matrices $\bm{\tilde{P}}_j^1$ and $\bm{\overline{P}}^1_j$ are those in Eqs. (\ref{theorem2-parameter}) and (\ref{p1j}). It is evident that based on the Corollary 2, Proposition 3, and Eq. (\ref{equivalem-mrt-model}), the four GPMRT-LB models ($\bm{M}_1,\bm{S}_1$), ($\bm{M},\bm{S}$), ($\bm{M},\bm{\hat{S}}$) and ($\bm{M}_1,\bm{\hat{S}}_1$) are equivalent. 
	\begin{figure}[htbp]    
			\centering            
			\includegraphics[width=1.0\textwidth]{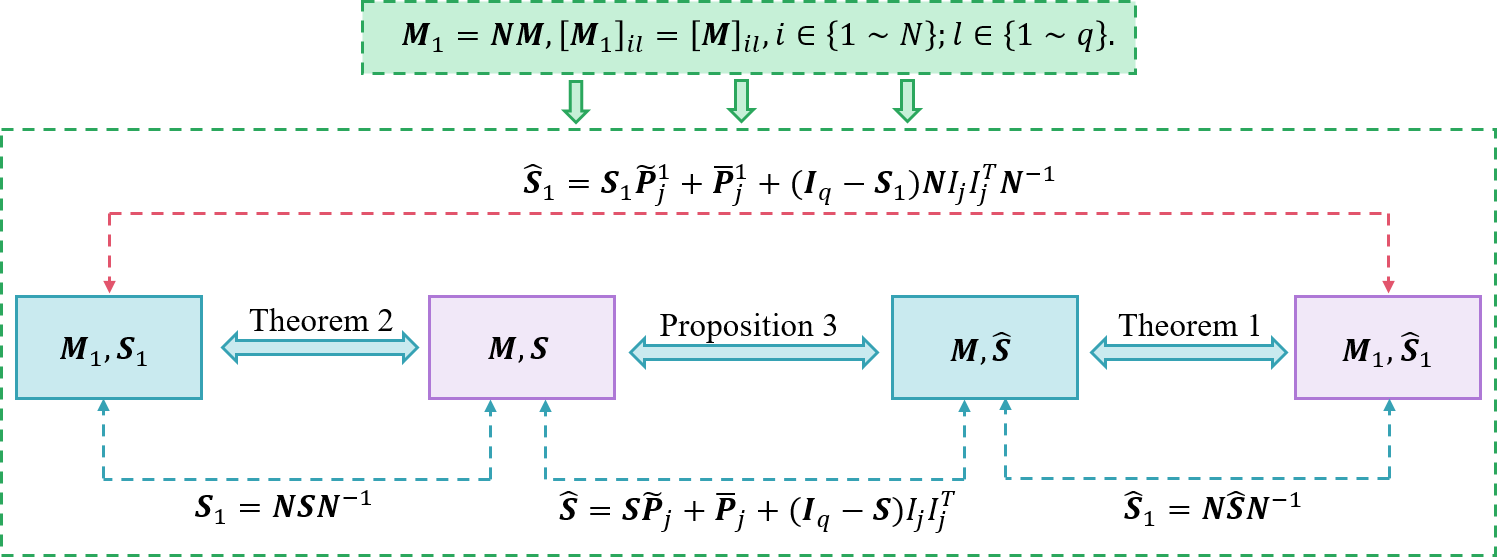}  
			\caption{The relations between the finite-difference schemes of the four equivalent GPMRT-LB models ($\bm{M}_1,\bm{S}_1$), ($\bm{M},\bm{S}$), ($\bm{M},\bm{\hat{S}}$), and ($\bm{M}_1,\bm{\hat{S}}_1$). }    
			\label{relation-GPMRT-model}            
		\end{figure}
Moreover, according to the Proposition 3, Theorems 1 and 2, the four equivalent GPMRT-LB models ($\bm{M}_1,\bm{S}_1$), ($\bm{M},\bm{S}$), ($\bm{M},\bm{\hat{S}}$)  and ($\bm{M}_1,\bm{\hat{S}}_1$)  satisfy the relations presented in Fig. \ref{relation-GPMRT-model} where the double-line arrow connecting the two boxes indicates that the two GPMRT-LB models have the same finite-difference scheme. \\
\textbf{Remark 7} Regarding the finite-difference scheme (\ref{det-GPMFD-N-geq1}) with $N\geq 1$ conservative moments in Proposition 2, the characteristics shown in the Corollary 2 and Theorem 2 are consistent with the GPMRT-LB model. Therefore, we refer to it as the macroscopic finite-difference (GPMFD) scheme of the GPMRT-LB model. Furthermore, since the finite-difference schemes (\ref{det-GPMFD-N-geq1}) and (\ref{det-GPMFD-N-geq1-brevity}) have the same form (see the Proposition 3), we only need to consider the latter in the following discussion, and this will simplify the analysis on truncation errors and MEs.

	\section{Truncation errors and MEs of the GPMRT-LB model and GPMFD scheme }\label{Truncation-Error-Modified-Equation}
In this section, we will conduct some theoretical analysis on the truncation errors and MEs of the GPMRT-LB model (\ref{GPMRT-LB-model}) and GPMFD scheme (\ref{det-GPMFD-N-geq1-brevity}). It is known that the GPMRT-LB model $(\bm{M}^c,\bm{S}^c)$ can be equivalent to a GPMRT-LB model $(\bm{M}^c_N,\bm{S}_N^c)$ [see Eq. (\ref{equivalem-mrt-model})] where the transform matrix $\bm{M}_N^c$ is based on the natural moment and $\bm{S}^c_N$ is a block-lower-triangular relaxation matrix. Then from the two purple boxes shown in Fig. \ref{relation-GPMRT-model}, one can find that the  GPMFD schemes corresponding  to the equivalent GPMRT-LB models ($\bm{M}^c,\bm{S}^c$) and ($\bm{M}^c_N,\bm{\hat{S}}^c_N$)  where the relaxation matrices $\bm{S}^c$ and $\bm{\hat{S}}^c_N$ [see Eq. (\ref{hats1s2})] are assumed to be independent on the parameter $c$, are identical. According to the Remark 5 and for convenience of the  analysis on the truncation errors and MEs, we only need to consider the GPMRT-LB model $(\bm{M}_N,\hat{\bm{S}}_N)$ satisfying the relations $\bm{M}_N=\bm{C}_d^{-1}\bm{M}^c_N$ and $\hat{\bm{S}}_N=\bm{C}^{-1}_d\hat{\bm{S}}^c_N\bm{C}_d$. Without loss of generality, here we consider that the degree of parameter $c$ does not decrease with the number of rows (columns) in the diagonal matrix $\bm{C}_d$. Therefore, for a given GPMRT-LB model ($\bm{M}^c,\bm{S}^c$), we focus on the GPMRT-LB model ($\bm{M}_N,\bm{\hat{S}}_N$) with the transform matrix $\bm{M}_N$ that is  independent on the parameter $c$ and its GPMFD scheme (\ref{det-GPMFD-N-geq1-brevity}). In addition, for the inverse of the collision matrix $\bm{\overline{\Lambda}}:=\bm{M}_N\bm{\hat{S}}_N\bm{M}_N^{-1}$, the following requirements are needed:
\begin{align}\label{natural-condition}
	&\sum_{j=1}^{q}\bm{e}_{j\alpha} \bm{\overline{\Lambda}}_{jk}= S^{10}_{\alpha }\bm{i}_k+S^{ 1}_{\alpha \xi_1}\bm{e}_{k\xi_1 }, 	\sum_{j=1}^{q}\bm{e}_{j\alpha}\bm{e}_{j\beta}\bm{\overline{\Lambda}}_{jk}= S^{20}_{\alpha\beta}\bm{i}_k+S^{21}_{\alpha\beta\xi_1}\bm{e}_{k\xi_1}+S^{2}_{\alpha\beta\xi_1\xi_2}\bm{e}_{k\xi_1}\bm{e}_{k\xi_2},\notag\\
	&\sum_{j=1}^{q }\bm{e}_{j\alpha}\bm{e}_{j\beta}\bm{e}_{j\gamma}\bm{\overline{\Lambda}}_{jk}=S^{30}_{\alpha\beta\gamma}\bm{i}_k+S^{31}_{\alpha\beta\gamma\xi_1}\bm{e}_{k\xi_1}+S^{32}_{\alpha\beta\gamma\xi_1\xi_2}\bm{e}_{k\xi_1}\bm{e}_{k\xi_2}+S^{33}_{\alpha\beta\gamma\xi_1\xi_2\xi_3}\bm{e}_{k\xi_1}\bm{e}_{k\xi_2}\bm{e}_{k\xi_3},\notag\\
	&\sum_{j=1}^{q }\bm{e}_{j\alpha}\bm{e}_{j\beta}\bm{e}_{j\gamma}\bm{e}_{j\eta}\bm{\overline{\Lambda}}_{jk}=S^{40}_{\alpha\beta\gamma\eta}\bm{i}_k+S^{41}_{\alpha\beta\gamma\eta\xi_1}\bm{e}_{k\xi_1}+S^{42}_{\alpha\beta\gamma\eta\xi_1\xi_2}\bm{e}_{k\xi_1}\bm{e}_{k\xi_2}\notag\\
	&\qquad\qquad\qquad\qquad\qquad+S^{43}_{\alpha\beta\gamma\eta\xi_1\xi_2\xi_3}\bm{e}_{k\xi_1}\bm{e}_{k\xi_2}\bm{e}_{k\xi_3}+S^{44}_{\alpha\beta\gamma\eta\xi_1\xi_2\xi_3\xi_4}\bm{e}_{k\xi_1}\bm{e}_{k\xi_2}\bm{e}_{k\xi_3}\bm{e}_{k\xi_4}, 
\end{align}	 
where  $\bm{i}_k$ indicates the $k_{th}$ element of vector $\bm{i}=\big(1,1,\ldots,1)\in R^q$, $\bm{S}^{li}$ is a $d^l\times d^i$ matrix $\big(l\in\{1\sim 2\};i\in\{0\sim (l-1)\}$ and $l\in\{3\sim 4\};i\in\{0\sim l\}\big)$,  $\bm{S}^{1}$
is an invertible $d\times d$ relaxation matrix associated with the
diffusion coefficient matrix of the NACDE, and $\bm{S}^{2}$ is a $d^2\times d^2$ relaxation matrix associated with the viscosity coefficient of the NSEs. Specifically, we take $\bm{S}^{21}=\bm{0}$ for the NACDE  and  $\bm{S}^{32}=\bm{0}$ for the NSEs.

Due to the equivalence between the GPMRT-LB model (\ref{GPMRT-LB-model}) and GPMFD scheme (\ref{det-GPMFD-N-geq1}) discussed in Section \ref{Derivation-GPMFD} and Remark 7, we will focus on the GPMFD scheme (\ref{det-GPMFD-N-geq1-brevity}) and present some details to derive its truncation error and ME. First, we decompose the matrix $\bm{W}$ in Eq.  (\ref{FD-scheme-quan}) as
\begin{equation}
	\bm{W}=p_0\bm{I}_q+p_{-1}\bm{W}_{-1}+p_1\bm{W}_1,
\end{equation}
where $\bm{W}_{-1}=\bm{M}_N\bm{T}_{-1}\bm{M}_N^{-1}$ and $\bm{W}_1=\bm{M}_N\bm{T}_1\bm{M}_N^{-1}$. Thanks to the space shift operator with a series form (\ref{series-space-oper}), the matrices $\bm{W}_{-1}$ and $\bm{W}_1$ can be rewritten as  
\begin{subequations}\label{W-1W1}
	\begin{align}
		&\bm{W}_{-1}= \bigg(\sum_{k=0}^{+\infty}\frac{(-\Delta x)^k}{k!} \bm{\mathcal{W}}_{0}^k\bigg)=\exp{\Big(- \Delta x \bm{\mathcal{W}_{0}}\Big)}, \\
		&\bm{W}_1= \bigg(\sum_{k=0}^{+\infty}\frac{(\Delta x)^k}{k!}\bm{\mathcal{W}}_{0}^k\bigg) =\exp{\Big( \Delta x \bm{\mathcal{W}_{0}}\Big)}, 
	\end{align}
	\end{subequations}
where $\bm{\mathcal{W}}_{0}=\bm{M}_N\Big[$\textbf{diag}$\big(\bm{e}_1\cdot \nabla,\bm{e}_2\cdot \nabla,\ldots,\bm{e}_{q}\cdot \nabla\big) \Big]\bm{M}_N^{-1} $,  and the inverse of matrix $\bm{W}$ can be further expressed as
\begin{align}\label{inverse-W}
	\bm{W}^{-1} =\sum_{k=0}^{+\infty}\bigg[\bm{I}_q-p_0\bm{I}_q-p_{-1}\exp{\Big(- \Delta x \bm{\mathcal{W}_{0}}\Big)}-p_1\exp{\Big(\Delta x \bm{\mathcal{W}_{0}}\Big)} \bigg]^k ,
\end{align} 
which will be used below. 

Based on the Maxwell iteration method  \cite{Ikenberry1956, Yong2016} and the relation between adjugate matrix and determinant in Eq.  (\ref{relation-det-adj}), we set $\bm{\psi}=\bm{\psi}^n$ with $\bm{\psi}$ representing $\{\bm{m},\bm{m}^{eq},\bm{\tilde{F}} \}$,  and substitute the  matrices $\bm{M}_N$ and $\bm{\hat{S}}_N$ into the GPMFD scheme (\ref{det-GPMFD-N-geq1-brevity}), one can derive
\begin{align}\label{detadjmaxwell-2}
	\bm{0}&=	\det{\big(T^1_{\Delta t}\bm{I}-\bm{A}\big)}\bm{m} -\Big({\rm adj}\big(T^1_{\Delta t}\bm{I}-\bm{A}\big)\big(\bm{B}\bm{m}^{eq }+\Delta t\bm{W}\bm{\tilde{F}}\big)\Big)\notag \\
	&=\det{\big(T^1_{\Delta t}\bm{I}-\bm{A}\big)}\bigg[\bm{m} -\Big(T^1_{\Delta t}\bm{I}-\bm{W}(\bm{I}-\bm{\hat{S}}_N)\Big)^{-1}\Big(\bm{W}\bm{\hat{S}}_N\bm{m}^{eq,n} +\bm{W}\Delta t\bm{\tilde{F}} \Big) \bigg]\notag \\
	&=\det{\big(T^1_{\Delta t}\bm{I}-\bm{A}\big)}\Bigg[\bm{m} -\bigg(\bm{\hat{S}}_N^{-1}\bm{W}^{-1}\Big(T^1_{\Delta t}\bm{I}-\bm{W}(\bm{I}-\bm{\hat{S}}_N)\Big)\bigg)^{-1}\Big( \bm{m}^{eq } +\bm{\hat{S}}_N^{-1}\Delta t\bm{\tilde{F}}  \Big) \Bigg]\notag \\
	&=\det{\big(T^1_{\Delta t}\bm{I}-\bm{A}\big)}\Bigg[\bm{m} -\bigg(
	\bm{I}+\bm{\hat{S}}_N^{-1}\Big(T^1_{\Delta t}\bm{W}^{-1}-\bm{I}\Big)\bigg)^{-1}\Big( \bm{m}^{eq } +\bm{\hat{S}}_N^{-1}\Delta t\bm{\tilde{F}} \Big) \Bigg] \notag\\
	&:=\det{\big(T^1_{\Delta t}\bm{I}-\bm{A}\big)}\bm{\Xi}. 
\end{align}
Here $\bm{m}=\bm{M}_N\bm{f}, \bm{m}^{eq}=\bm{M}_N\bm{f}^{eq}, \bm{\tilde{F}}=\bm{M}_N\big(\bm{F}+\bm{G}+\Delta t/2\bm{\overline{D}}\bm{F}\big)$, and $\bm{\Xi}$ is defined as   
\begin{align}\label{Xi0}
	\bm{\Xi}  =  \bm{m} - \Bigg(\sum_{k=0}^{+\infty}{\bm{\Gamma}}^{k}\Bigg)\Big( \bm{m}^{eq } +\bm{\hat{S}}_N^{-1}\Delta t\bm{\tilde{F}}  \Big) ,
\end{align} 
where $\bm{\Gamma}=-\bm{\hat{S}}_N^{-1}\big(T^1_{\Delta t}\bm{W}^{-1}-\bm{I}_q\big)$. With the help of  $\bm{\Gamma}^0=\bm{I}_q$, Eq. (\ref{inverse-W}), and the series form  (\ref{series-time-oper}) of time operator $T^1_{\Delta t}$, the expression of  $\bm{\Xi}$ in Eq. (\ref{Xi0}) can be given as
\begin{align}\label{Xi} 
	\bm{\Xi}  & = \bm{m}-\big(\bm{m}^{eq}+ \Delta t\bm{\hat{S}}_N^{-1}\bm{\tilde{F}}\big)-\sum_{k=1}^{+\infty}\Delta x^k \bm{\Xi}^{(k)}, 
\end{align} 
where  $\bm{\Xi}^{(k)} =\bm{\Gamma}^{(k)}\Big( \bm{m}^{eq } +\bm{\hat{S}}_N^{-1}\Delta t\bm{\tilde{F}}  \Big)$ $(k\geq 1)$ is the coefficient before the $k_{th}$-order term of the series expansion of $\bm{\Xi}$ with $\bm{\Gamma}^{(k)}$ denoting that of the series expansion of $\bm{\Gamma}$.

Based on  the fact that $\det{\big(T^1_{\Delta t}\bm{I}_q-\bm{A}\big)}=\det{(\bm{\hat{S}}_N)}+O(\Delta x)$  \cite{Bellotti-1-2022}, we will show that the analysis on the truncation error and ME of scheme (\ref{detadjmaxwell-2}) is equivalent to the discussion  on $\bm{\Xi}$ (\ref{Xi}). The reason is as follows.

According to Eq. (\ref{detadjmaxwell-2}), we have
\begin{align}\label{maxwell-eq1}
	\bm{0}=&\det{\big(T^1_{\Delta t}\bm{I}-\bm{A}\big)}\bm{\Xi} =  \det{(\bm{\hat{S}}_N)} \bigg[\bm{m}-\big(\bm{m}^{eq}+\Delta t\bm{\hat{S}}_N^{-1}\bm{\tilde{F}}\big)-\sum_{k=1}^{+\infty}\Delta x^k\bm{\Xi}^{(k)} \bigg]\notag\\
	&+ \bigg[\Big(\bm{m}-\big(\bm{m}^{eq}+ \Delta t\bm{\hat{S}}_N^{-1}\bm{\tilde{F}}\big)\Big)\times O(\Delta x)-\sum_{k=1}^{+\infty}\Delta x^{k+1}\bm{\Xi}^{(k)} \bigg].
\end{align} 
As mentioned previously, to perform the analysis on the truncation errors and the MEs, it is necessary to specify the scaling relationship between $\Delta t$ and $\Delta x$. In particular, it should be noted that at the diffusive scaling ($\Delta t\sim\Delta x^2$), $\bm{\Xi}^{(k)}$ is at least of order $O(1)$, which arises from the fact that  $\bm{m}^{eq}$, $\bm{\tilde{F}}$  and $\bm{\hat{S}}_N^{-1}$ (\ref{hats1s2})  are at least of order  $O(1)$ $\big[$
the main diagonal element of the block-lower-triangular matrix $\bm{\hat{S}}_N^{-1}$ is of order $O(1)$, while the lower triangular element of $\bm{\hat{S}}_N^{-1}$ is at least of order $O(1)$$\big]$. Due to the conservation law, i.e., $m_{i\in\{1\sim N\}}=m^{eq}_{i\in\{1\sim N\}}$, for any $i_{th}$ $\big(i\in\{1\sim N\}\big)$ row of Eq. (\ref{maxwell-eq1}), we have
\begin{subequations}
	\begin{align}
		&{\rm diffusive\:\: scaling}\:(\Delta t\sim \Delta x^2):\: \det{(\bm{\hat{S}}_N)}\big[\Delta x\bm{\Xi}^{(1)} \big]_i=O(\Delta x^2),\\
		&{\rm acoustic\:\: scaling}\:(\Delta t\sim \Delta x):\:\det{(\bm{\hat{S}}_N)}\big[\Delta t\bm{\hat{S}}_N^{-1}\bm{\tilde{F}}+\Delta x\bm{\Xi}^{(1)} \big]_i=O(\Delta x^2), \label{maxwell-eq2}
	\end{align} 
\end{subequations}  
where $\det{(\bm{\hat{S}}_N)}$ is of order $O(1)$, this is because the main diagonal element of the block-lower-triangular matrix $\bm{\hat{S}}_N$ is of order $O(1)$. Based on Eq. (\ref{maxwell-eq2}), one can obtain the third-order ME at the acoustic scaling,
\begin{align}\label{acou-third-order}
	&  \Big[\Delta t\bm{\hat{S}}_N^{-1}\bm{\tilde{F}}+\Delta x\bm{\Xi}^{(1)}+\Delta x^2\bm{\Xi}^{(2)}\Big]_i=O(\Delta x^3).
\end{align} Then according to Eq. (\ref{acou-third-order}), one can further derive fourth- (and higher-) order truncation errors and MEs of the GPMRT-LB model (\ref{GPMRT-LB-model}) and GPMFD scheme (\ref{det-GPMFD-N-geq1-brevity}) at the acoustic scaling. In addition, it can also be observed that the analysis on the truncation errors and MEs   of GPMFD scheme (\ref{detadjmaxwell-2}) are actually equivalent to the discussion on $\bm{\Xi}$ (\ref{Xi}), which is also true at the diffusive scaling. 
\subsection{MEs at the diffusive scaling}
Let us begin the analysis at the diffusive scaling, i.e., $\Delta t=\eta\Delta x^2,\eta\in\mathcal{R}$.  It should be noted that the lattice velocity $\lambda$ is of order $O(1/\Delta x)$, $\bm{m}^{eq}$, $\bm{\tilde{F}}$, and $\bm{\hat{S}}^{-1}_N$ are at least of order $O(1)$. We would also like to point out that $\Delta x^{-k_i}m_i$ $\big(i\in\{1\sim N\}\big)$ is of order $O(1)$ owing to the fact that  the transform matrix $\bm{M}$ is of order $O(1)$ and the form of the equilibrium distribution function for some specific problems [see Eq. (\ref{NACDE-feq}) for NACDE and Eq. (\ref{NSEs-feq}) for NSEs, $k_i\geq 0$ denotes the order of the $i_{th}$ conservative moment $m_i$ in space]. Thus, a higher-order expansion of $\bm{\Xi}$ (\ref{Xi}) beyond $\bm{\Xi}^{(k)}$ is necessary for the $k_{th}$-order ME of the GPMRT-LB model (\ref{GPMRT-LB-model}) and GPMFD scheme (\ref{det-GPMFD-N-geq1-brevity}). Here we expand $\bm{\Xi}$ (\ref{Xi}) up to $\bm{\Xi}^{(4)}$ and consider the first- to third-order MEs on conservative moment $m_i$ $\big(i\in\{1\sim N\}\big)$, 
\begin{align}\label{five-order-diffusive}
	&    \Big[ \Delta t\bm{\tilde{F}}+\Delta x\bm{\hat{S}}_N\bm{\Xi}^{(1)}+\Delta x^2\bm{\hat{S}}_N\bm{\Xi}^{(2)}+\Delta x^3\bm{\hat{S}}_N\bm{\Xi}^{(3)}   
	+\Delta x^4\bm{\hat{S}}_N\bm{\Xi}^{(4)} \Big]_i=O(\Delta x^5),
\end{align}
where the $i_{th}$ row of matrix $\bm{\hat{S}}_N$ (\ref{hats1s2}) is identical to $I_i^T$, which is of order $O(1)$ and has been used.

Moreover, it is also necessary to expand the inverse of matrix $\bm{W}$ (\ref{inverse-W}) as
\begin{align}
	\bm{W}^{-1} =&\bm{I}_q+\Delta xa \bm{\mathcal{W}}_0+\Delta x^2a ^2\Big(1-\frac{b}{2a^2}\Big)\bm{\mathcal{W}}_0^2+\Delta x^3a ^3\Big(\frac{1}{6a^3} -\frac{b}{a^2}+1\Big)\bm{\mathcal{W}}_0^3\notag\\
	&+\Delta x^4a^4\Big(\frac{b^2}{4a^4}+\frac{1}{3a^3}-\frac{3b}{2a^2}+1-\frac{b}{24a^4}\Big)\bm{\mathcal{W}}_0+O(\Delta x^5). 
\end{align}
By adopting the series form (\ref{series-time-oper}) of the time shift operator $T^1_{\Delta t}$ and Eq. (\ref{W-1W1}), one can obtain
\begin{align}\label{SGamma0}
	-\bm{\hat{S}}_N\bm{\Gamma}=&\Delta xA_{1}+\Delta x^2\Big(A_{12}+  A_{2 }\Big)\notag\\
	&+\Delta x^3\Big(A_{31} +A_{3}\Big)+\Delta x^4\Big(A_{41}+A_{42}+A_4\Big)+O(\Delta x^5),
\end{align}
where 
\begin{subequations} \label{SGamma}
	\begin{align}
		&A_{1}=a\bm{\mathcal{W}}_0,\label{A1}\\
		&A_{21}=\eta\partial_t\bm{I}_q,A_{2}=\bm{\mathcal{W}}_0^2\Big(a^2-\frac{b}{2}\Big),\label{A2}\\
		&A_{31}=a\eta\bm{\mathcal{W}}_0\partial_t,A_{3}=\Big(\frac{1}{6}-ab+a^3\Big)\bm{\mathcal{W}}_0^3,\label{A3}\\
		&A_{41}=\eta\Big(a^2-\frac{b}{2}\Big)\bm{\mathcal{W}}_0^2\partial_t,A_{42}=\eta^2\frac{\partial_{tt}\bm{I}_q}{2},A_4=\Big(\frac{b^2}{4}+\frac{a}{3}-\frac{3a^2b}{2}+a^4-\frac{b}{24 }\Big)\bm{\mathcal{W}}_0^4,\label{A4}
	\end{align}
\end{subequations} 
In the following, the analysis are based on Eqs.  (\ref{five-order-diffusive}) and (\ref{SGamma}).

\subsubsection{First-order ME of the GPMRT-LB model and  GPMFD scheme on conservative moment $m_i$}\label{first-diffusive}
For any $i\in\{1\sim N\}$, to obtain the first-order ME  of the GPMRT-LB model and  GPMFD scheme on the moment $m_i$ at the diffusive scaling, one needs to consider the following truncation equation of Eq. (\ref{five-order-diffusive}):
\begin{equation}\label{Xi2}
	\Big[\Delta t\bm{\tilde{F}}+\Delta x\bm{\hat{S}}_N\bm{\Xi}^{(1)}+\Delta x^2\bm{\hat{S}}_N\bm{\Xi}^{(2)}\Big]_i=O(\Delta x^3).
\end{equation}

With the aid of Eqs. (\ref{Xi}) and (\ref{SGamma0}), Eq. (\ref{Xi2}) can be written as
\begin{align}\label{modified-equation-first-order-diffusive-0}
	\Bigg[\Delta t\bm{\tilde{F}}-\bigg[\Delta xA_1+  \Delta x^2A_{12}+ \Delta x^2  \Big(A_2-A_1\bm{\hat{S}}_N^{-1}A_1\Big)  \bigg]  \bm{m}^{eq}\Bigg]_i+O(\Delta x^3), 
\end{align} 
multiplying above equation by $1/\Delta t$ yields 
\begin{align}
	\Bigg[\bm{\tilde{F}}- \bigg[\lambda A_1+ \frac{A_{12}}{\eta}+  \lambda^2\Delta t\Big(A_2-A_1\bm{\hat{S}}_N^{-1}  A_1\Big) \bigg]  \bm{m}^{eq}\Bigg]_{i}=O(\Delta x), 
\end{align} 
then substituting Eqs. (\ref{A1}) and (\ref{A2}) into  Eq. (\ref{modified-equation-first-order-diffusive-0}), we have 
\begin{align}\label{modified-equation-first-order-diffusive}
	& \Big(\partial_t \bm{m}^{eq}+ c \bm{\mathcal{W}}_0 \bm{m}^{eq}\Big)_{i}=\Bigg[ \bm{\tilde{F}}+    \frac{a^2}{\eta}  \bm{\mathcal{W}}_0\bigg[\bm{\hat{S}}_N^{-1}+\Big(\frac{b}{2a^2}-1\Big)\bm{I}\bigg]    \bm{\mathcal{W}}_0  \bm{m}^{eq}\Bigg]_{i}+O(\Delta x).
\end{align}

\subsubsection{Second-order ME of the GPMRT-LB model and  GPMFD scheme on conservative moment $m_i$}\label{second-diffusive}
Similar to the discussion in the previous part, we consider the following truncation equation of Eq. (\ref{five-order-diffusive}):
\begin{equation}\label{Xi3}
	\Big[\Delta t\bm{\tilde{F}}+\Delta x\bm{\hat{S}}_N\bm{\Xi}^{(1)}+\Delta x^2\bm{\hat{S}}_N\bm{\Xi}^{(2)}+\Delta x^3\bm{\hat{S}}_N\bm{\Xi}^{(3)}\Big]_i=O(\Delta x^4),
\end{equation}
for any $i\in\{1\sim N\}$, one can obtain
\begin{align}\label{modified-equation-second-order-diffusive-0}
	&\Bigg[\Delta t\bm{\tilde{F}}-\Delta x\Delta tA_1\bm{\hat{S}}_N^{-1}\bm{\tilde{F}}\notag\\
	&-\Delta x\bigg[A_1+\Delta x\Big(A_{12}+A_{2}\Big)+\Delta x^2\Big(A_{31}+A_3\Big)\bigg]\bm{m}^{eq}\notag\\
	& +\Delta x^3 \bigg[A_1\bm{\hat{S}}_N^{-1}\Big(A_{21}+A_2\Big) +\Big(A_{21}+A_2\Big)\bm{\hat{S}}_N^{-1}A_1 - A_1\bm{\hat{S}}_N^{-1}A_1\bm{\hat{S}}_N^{-1}A_1\bigg]\bm{m}^{eq} \Bigg]_i =O(\Delta x^4) .
\end{align}
Substituting Eqs. (\ref{A1}$\sim$\ref{A3}) into above Eq. (\ref{modified-equation-second-order-diffusive-0}) yields 
\begin{align}\label{modified-equation-second-order-diffusive}
	&\Big[\big(\partial_t\bm{I}+c \bm{\mathcal{W}}_0\big)\bm{m}^{eq} \Big]_i=\notag\\
	&\qquad\biggg[\bm{\tilde{F}}+\frac{a^2}{\eta} \bm{\mathcal{W}}_0 \bigg(\bm{\hat{S}}_N^{-1}+ \Big(\frac{b}{2a^2}-1\Big)\bm{I}\bigg)    \bm{\mathcal{W}}_0   \bm{m}^{eq}-a\Delta x  \partial_t\big( \bm{I}-\bm{\hat{S}}_N^{-1}\big) \bm{\mathcal{W}}_0 \bm{m}^{eq}	  \notag\\ 
	& \qquad+a\Delta x\bm{\mathcal{W}}_0  \bm{\hat{S}}_N^{-1} \Bigg[  \partial_t \bm{m}^{eq}- \bm{\tilde{F}}-\frac{a^2}{\eta} \bm{\mathcal{W}}_0\bigg(\bm{\hat{S}}_N^{-1} +\Big(\frac{b}{2a^2}-1\Big)\bm{I}\bigg) \bm{\mathcal{W}}_0 \bm{m}^{eq}\Bigg]\notag\\
	&\qquad-\frac{a^3\Delta x}{\eta}\Big(\frac{b}{2a^2}-1\Big) \bm{\mathcal{W}}_0^2 \bm{\hat{S}}_N^{-1}\bm{\mathcal{W}}_0  \bm{m}^{eq} -\frac{a^3\Delta x}{\eta}\Big(\frac{1}{6a^3}-\frac{b}{a^2}+1\Big)\bm{\mathcal{W}}_0^3 \bm{m}^{eq}\biggg]_i+O(\Delta x^2) .
\end{align} 
\subsubsection{Third-order ME of the GPMRT-LB model and  GPMFD scheme on conservative moment  $m_i$}\label{third-diffusive}
For any $i\in\{1\sim N\}$, after some manipulations by using Eqs. (\ref{Xi}) and (\ref{SGamma0}), we can derive the third-order ME, 
\begin{align}\label{modified-equation-third-order-diffusive-0}
	&\Bigg[\Delta t\bm{\tilde{F}} -\Big[ \Delta xA_{1}+\Delta x^2\big(A_{12}+  A_{2}\big)+\Delta x^2\big( A_{31} + A_{3}\big) +\Delta x^4\Big(A_{41}+A_{42}+A_{4}\big) \Big] \bm{m}^{eq}\notag\\
	&- \Big[ \Delta xA_{1}+\Delta x^2\big(A_{21}+ A_{2}\big)  \Big]\Delta t \bm{\hat{S}}_N^{-1}\bm{\tilde{F}} +\Delta x^2\Big(\bm{\hat{S}}_N^{-1}A_{1}\Big)^2 \Delta t\bm{\hat{S}}_N^{-1}\bm{\tilde{F}} \notag\\
	&+\bm{\hat{S}}_N\bigg(\bm{\hat{S}}_N^{-1}\Big[\Delta xA_{1}+\Delta x^2\big(A_{21}+A_{2}\big)+\Delta x^3\big(
	A_{31}+A_{3}\big) \Big]\bigg)^2 \bm{m}^{eq} \notag \\
	&-\bm{\hat{S}}_N\bigg(\bm{\hat{S}}_N^{-1}\Big[\Delta xA_{1}+\Delta x^2\big(A_{21}+ A_{2}\big) \Big]\bigg)^3\bm{m}^{eq}  +\Delta x^4\Big(\bm{\hat{S}}_N^{-1}A_1\Big)^4 \bm{m}^{eq}\Bigg]_i =O(\Delta x^5).
\end{align} 
Substituting Eqs. (\ref{A1}$\sim$\ref{A4}) into  Eq. (\ref{modified-equation-third-order-diffusive-0}) gives
\begin{align}\label{modified-equation-third-order-diffusive}
	&\Bigg[\bm{\tilde{F}}-\bigg[c\bm{\mathcal{W}}_0+\partial_t\bm{I}+\frac{a^2}{\eta}\Big(1-\frac{b}{2a^2}\Big)\bm{\mathcal{W}}_0^2 +a\Delta x\bm{\mathcal{W}}_0\partial_t+\frac{a^3\Delta x}{\eta}\Big(1-\frac{b}{2a^2}\Big)\bm{\mathcal{W}}_0^2\partial_t\bm{I}\notag\\
	&\qquad\qquad+\frac{\Delta x}{2\lambda}\partial_{tt}\bm{I}+\frac{a^3\Delta x}{\eta}\Big(\frac{1}{6a^2}-\frac{b}{a^2}-1\Big)\bm{\mathcal{W}}_0^3\notag\\
	&\qquad\qquad+\frac{a^4\Delta x^2}{\eta}\Big(\frac{b^2}{4a^2}+\frac{1}{3a^3}-\frac{3b}{2a^2}+1-\frac{b}{24a^4}\Big)\bm{\mathcal{W}}_0^4\bigg]\bm{m}^{eq}\notag\\
	&-\Big(c\bm{\mathcal{W}}_0+\partial_t\bm{I}+\frac{a^2}{\eta}\Big(1-\frac{b}{2a^2}\Big)\bm{\mathcal{W}}_0^2\Big)\frac{\Delta x}{\lambda}\bm{\hat{S}}_N^{-1}\bm{\tilde{F}} +\Big(\frac{a^2}{\eta}\bm{\mathcal{W}}_0\bm{\hat{S}}_N^{-1}\bm{\mathcal{W}}_0\Big)\frac{\Delta x}{\lambda}\bm{\hat{S}}_N^{-1}\bm{\tilde{F}}\notag\\
	&+\bigg[\frac{a^2}{\eta}\bm{\mathcal{W}}_0\bm{\hat{S}}_N^{-1}\bm{\mathcal{W}}_0+a\Delta x\bm{\mathcal{W}}_0\bm{\hat{S}}_N^{-1}\partial_t\bm{I}+\frac{a^3\Delta x}{\eta}\Big(1-\frac{b}{2a^2}\Big)\bm{\mathcal{W}}_0\bm{\hat{S}}_N^{-1}\bm{\mathcal{W}}_0^2\notag\\
	&\qquad+a^2\Delta x^2\bm{\mathcal{W}}_0\bm{\hat{S}}_N^{-1}\bm{\mathcal{W}}_0\partial_t+a\partial_t\bm{I}\bm{\hat{S}}_N^{-1}\bm{\mathcal{W}}_0+\frac{\Delta x}{\lambda}\partial_t\bm{\hat{S}}_N^{-1}\bm{I}\partial_t\bm{I}\notag\\
	&\qquad+a^2\Delta x^2\Big(1-\frac{b}{2a^2}\Big)\partial_t\bm{I}\bm{\hat{S}}_N^{-1}\bm{\mathcal{W}}_0^2\notag+\frac{a^3\Delta x}{\eta}\Big(1-\frac{b}{2a^2}\Big)\bm{\mathcal{W}}_0^2\bm{\hat{S}}_N^{-1}\bm{\mathcal{W}}_0\notag\\
	&\qquad+a^2\Delta x^2\Big(1-\frac{b}{2a^2}\Big)\bm{\mathcal{W}}_0^2\bm{\hat{S}}_N^{-1}\partial_t\bm{I}+\frac{a^4\Delta x^2}{\eta}\Big(1-\frac{b}{2a^2}\Big)^2\bm{\mathcal{W}}_0^2\bm{\hat{S}}_N^{-1}\bm{\mathcal{W}}_0^2\notag\\
	&\qquad +a^2\Delta x^2\bm{\mathcal{W}}_0\partial_t\bm{I}\bm{\hat{S}}_N^{-1}\bm{\mathcal{W}}_0\bigg]\bm{m}^{eq}-\bigg[\frac{a^3\Delta x}{\eta}\bm{\mathcal{W}}_0\bm{\hat{S}}_N^{-1}\bm{\mathcal{W}}_0\bm{\hat{S}}_N^{-1}\bm{\mathcal{W}}_0\notag\\
	&+a^2\Delta x^2\bm{\mathcal{W}}_0\bm{\hat{S}}_N^{-1}\bm{\mathcal{W}}_0\bm{\hat{S}}_N^{-1}\partial_t\bm{I}+a^2\Delta x^2\bm{\mathcal{W}}_0\partial_t\bm{I}\bm{\hat{S}}_N^{-1}\bm{\hat{S}}_N^{-1}\bm{\mathcal{W}}_0\notag\\
	&\qquad+\frac{a^4\Delta x^2}{\eta}\Big(1-\frac{b}{2a^2}\Big)\bm{\mathcal{W}}_0\bm{\hat{S}}_N^{-1}\bm{\mathcal{W}}_0\bm{\hat{S}}_N^{-1}\bm{\mathcal{W}}_0^2\notag\\
	&\qquad+\frac{a^4\Delta x^2}{\eta}\Big(1-\frac{b}{2a^2}\Big)\bm{\mathcal{W}}_0\bm{\hat{S}}_N^{-1}\bm{\mathcal{W}}_0^2\bm{\hat{S}}_N^{-1}\bm{\mathcal{W}}_0\notag\\
	&\qquad+a^2\Delta x^2\partial_t\bm{I}\bm{\hat{S}}_N^{-1}\bm{\mathcal{W}}_0\bm{\hat{S}}_N^{-1}\bm{\mathcal{W}}_0+\frac{a^4\Delta x^2}{\eta}\Big(1-\frac{b}{2a^2}\Big)\bm{\mathcal{W}}_0^2\bm{\hat{S}}_N^{-1}\bm{\mathcal{W}}_0\bm{\hat{S}}_N^{-1}\bm{\mathcal{W}}_0\bigg]\bm{m}^{eq} \notag\\
	&\qquad+\frac{a^4\Delta x^2}{\eta}\bm{\mathcal{W}}_0\bm{\hat{S}}_N^{-1}\bm{\mathcal{W}}_0\bm{\hat{S}}_N^{-1}\bm{\mathcal{W}}_0\bm{\hat{S}}_N^{-1}\bm{\mathcal{W}}_0\bm{m}^{eq}\Bigg]_i=O(\Delta x^3).
\end{align} 
Now we give a remark on these   MEs at the diffusive scaling.\\
\textbf{Remark 8} According to the results shown in Parts \ref{first-diffusive}$\sim$\ref{third-diffusive}, we now consider a specific problem, i.e., the $d$-dimensional NSEs with $d+1$ conservative moments. In this case, the moment $m_1$ is $O(1)$, while $m_i$ $\big(i \in \{2 \sim (d+1)\}\big)$ is of order $O(\Delta x^1)$. Thus, Eqs. (\ref{modified-equation-first-order-diffusive}), (\ref{modified-equation-second-order-diffusive}) and (\ref{modified-equation-third-order-diffusive}) correspond to the first- to third-order MEs of the GPMRT-LB model and GPMFD scheme on conservative moment $\rho$, but the zeroth- to second-order MEs of the GPMRT-LB model and GPMFD scheme on conservative moment $\rho\bm{u}$.
\subsection{MEs  at the acoustic scaling}
At the acoustic scaling, all the lattice velocity $\lambda$, $\bm{m}^{eq}$,  $\bm{\tilde{F}}$,  and $\hat{\bm{S}}^{-1}_N$ are of order $O(1)$, and only a $k_{th}$-order expansion of $\bm{\Xi}$ is needed for the derivation of the $k_{th}$ order MEs of the GPMRT-LB model (\ref{GPMRT-LB-model}) and GPMFD scheme (\ref{det-GPMFD-N-geq1-brevity}). We now expand $\bm{\Xi}$ (\ref{Xi}) up to $\bm{\Xi}^{(2)}$ and consider the first- and second-order MEs on conservative moments $m_i$ for any $i\in\{1\sim N\}$,
\begin{equation}\label{second-order-acoustic-xi}
	\Big[\Delta t\bm{\tilde{F}}+\Delta x\bm{\hat{S}}_N\bm{\Xi}^{(1)}+\Delta x^2\bm{\hat{S}}_N\bm{\Xi}^{(2)}\Big]_i=O(\Delta x^3),
\end{equation}
subsequently, the second-order expansion  of $\bm{W}^{-1}$ (\ref{inverse-W}) can be given by
\begin{equation}\label{W-expand-2}
	\bm{W}^{-1}=\bm{I}_q+\Delta xa\bm{\mathcal{W}}_0+\Delta x^2\Big(1-\frac{b}{2a^2}\Big)a^2\bm{\mathcal{W}}_0^2+O(\Delta x^3).
\end{equation}
Similar to the analysis at  diffusive scaling,  one can obtain the expansion  of $-\bm{\hat{S}}_N\bm{\Gamma}$ based on above Eq. (\ref{W-expand-2}),
\begin{equation}
	-\bm{\hat{S}}_N\bm{\Gamma}=\Delta x\Big(B_{11}+B_1\Big)+\Delta x^2\Big(B_{21}+B_{22}+B_2\Big)+O(\Delta x^3),
\end{equation}
where
\begin{subequations}
	\begin{align}
		&B_{11}=\frac{a\partial_t}{c}\bm{I}_q,B_1=a\bm{\mathcal{W}}_0,\label{B1}\\
		&B_{21}=\frac{a^2\partial_{tt}}{2c^2}\bm{I}_q,B_{22}=\frac{a^2\bm{\mathcal{W}}_0\partial_t}{c},B_2=a^2\bm{\mathcal{W}}_0^2\Big(1-\frac{b}{2a^2}\Big).\label{B2}
	\end{align}
\end{subequations}
\subsubsection{First-order ME of the GPMRT-LB model and GPMFD scheme on conservative moment $m_i$}\label{first-acoustic}
Considering the first-order truncation equation of Eq. (\ref{second-order-acoustic-xi}), for any $i\in\{1\sim N\}$, one can obtain 
\begin{align}\label{first-modified-equation-acoustic-0}
	& \Big[\Delta t \bm{\tilde{F}}- \Delta x\big(B_{11}+B_1\big)\bm{m}^{eq}\Big]_i=O(\Delta x^2),
\end{align} 
substituting Eq. (\ref{B1}) into above Eq. (\ref{first-modified-equation-acoustic-0}) gives rise to the  first-order ME: 
\begin{align}\label{first-modified-equation-acoustic}
	&   \Big[\partial_t\bm{m}^{eq}+c \bm{\mathcal{W}}_0\bm{m}^{eq}\Big]_i=\big[ \bm{\tilde{F}}\big]_i+O(\Delta x ). 
\end{align}
\subsubsection{Second-order ME of the GPMRT-LB model and GPMFD scheme on conservative moment $m_i$}\label{second-acoustic}
Considering the second-order equation of Eq. (\ref{second-order-acoustic-xi}), for any $i\in\{1\sim N\}$, we have
\begin{align}\label{second-modified-equation-acoustic-0}
	&\Bigg[\Delta t\bm{\tilde{F}}-\Delta x\big(B_{11}+B_1\big)\big(\bm{m}^{eq}+\bm{\hat{S}}_N^{-1}\Delta t\bm{\tilde{F}}\big)\notag\\
	&+\Delta x^2\bigg(\bm{\hat{S}}_N\Big[\bm{\hat{S}}_N^{-1}\big(B_{11}+B_1\big)\Big]^2- \big(B_{21}+B_{22}+B_2\big)\bigg)\bm{m}^{eq}
	\Bigg]_i=O(\Delta x^3),
\end{align}
substituting Eqs. (\ref{B1}) and (\ref{B2}) into   Eq. (\ref{second-modified-equation-acoustic-0}) yields the following second-order ME: 
\begin{align}\label{second-modified-equation-acoustic}
	&\Bigg[\Delta t\bm{\tilde{F}}-  \Delta t\big(\partial_t\bm{I}+c\bm{\mathcal{W}}_0\big)\big(\bm{m}^{eq}\notag \\
	&-\frac{\Delta t^2}{2 }\Big(1-\frac{b}{a^2}\Big) c^2\bm{\mathcal{W}}_0^2\bm{m}^{eq} +\frac{\Delta t^2}{2 } \Big[\big(2\bm{\hat{S}}_N^{-1}-\bm{I}\big)\partial_{tt}\notag\\
	&\qquad+2c\big(\bm{\hat{S}}_N^{-1}\bm{\mathcal{W}}_0+\bm{\mathcal{W}}_0\bm{\hat{S}}_N^{-1}-\bm{\mathcal{W}}_0\big)\partial_t+c^2\bm{\mathcal{W}}_0\big(2\bm{\hat{S}}_N^{-1}-\bm{I}\big)\bm{\mathcal{W}}_0\Big]\bm{m}^{eq}\Bigg]_i =O(\Delta x^3).
\end{align}  
From above results, one can see that Eqs. (\ref{first-modified-equation-acoustic}) and (\ref{second-modified-equation-acoustic})  would reduce to the results presented in Ref. \cite{Bellotti-1-2022} when   $a=b=1$ and the term $\bm{\tilde{F}}$ is neglected.
\subsection{NACDE: MEs of the GPMRT-LB model and GPMFD scheme}
The $d$-dimensional  NACDE  with a source term can be expressed as 
\begin{equation}\label{NACDE}
	\partial_t\phi+\nabla\cdot\bm{B}=\nabla\cdot\big[\bm{\kappa}\cdot(\nabla\cdot \bm{D})\big]+R,
\end{equation}
where $\phi$ is a scalar variable related to both time $t$ and space $\bm{x}$, $R$ denotes the source term. $\bm{B} = (B_{\alpha})$ is a vector function, $\bm{\kappa}=(\kappa_{\alpha\beta})$ and $\bm{D}=(D_{\alpha\beta})$ are symmetric tensors (matrices), and they can be functions of $\phi$, $\bm{x}$, and $t$.

In order to recover the NACDE (\ref{NACDE}) correctly, some requirements  or moment conditions on the equilibrium, auxiliary, and source distribution functions, denoted by $f^{eq}_i$, $G_i$, and $F_i$, should be satisfied. For a general D$d$Q$q$ lattice structure, the moment conditions are given by
\begin{subequations}\label{NACDE-moment-conditions}
\begin{align}
	&\sum_{i=1}^qf_i=\sum_{i=1}^qf_i^{eq}=\phi,\sum_{i=1}^qF_i=R,\sum_{i=1}^qG_i=0,\\
	&\sum_{i=1}^q\bm{c}_if_i^{eq}=\bm{B},\sum_{i=1}^q\bm{c}_iF_i=\bm{0}, \sum_{i=1}^q\bm{c}_iG_i=\bm{M}_{1G},\\
	&\sum_{i=1}^q\bm{c}_i\bm{c}_if_i^{eq}=\chi c_s^2\bm{D}+ \bm{C},
\end{align}
\end{subequations}
where $c_s$ is a model parameter related to the lattice velocity $\lambda$. The parameter  $\chi$ is  used to adjust the relaxation matrix [see Eq. (\ref{kappa})], $\bm{C}$ is an auxiliary moment \cite{Chai2020}, and $\bm{M}_{1G}=\Big(\bm{I}-(\bm{S}^1)^{-1}/2\Big)\partial_t\bm{B}+\Big(\bm{I}-b(\bm{S}^1)^{-1}/(2a^2)\Big)\nabla\cdot\bm{C}$ is the first-order moment of $G_i$. From Eq. (\ref{NACDE-moment-conditions}) one can determine the expressions of $f^{eq}_i$, $G_i$, and $F_i$, while for simplicity, we only consider the following commonly used forms \cite{Chai2020}: 
\begin{subequations}\label{NACDE-feq}
\begin{align}
	&f_i^{eq}=w_i\biggg[\phi+\frac{\bm{c}_i\cdot\bm{B}}{c_s^2}+\frac{\Big(\chi c_s^2\bm{D}+ \bm{C}-c_s^2\phi\bm{I}\big):\big(\bm{c}_i\bm{c}_i-c_s^2\bm{I}\big)}{2c_s^4}\biggg],\\
	& G_i=w_i\bigg[ \frac{\bm{c}_i\cdot \bm{M}_{1G}}{c_s^2}\bigg], F_i=w_i R.
\end{align} 
\end{subequations}
According to the results shown in the Parts \ref{first-diffusive} and \ref{second-diffusive}, 
it is easy to obtain the first- and second-order MEs of the GPMRT-LB model (\ref{GPMRT-LB-model}) and GPMFD scheme (\ref{det-GPMFD-N-geq1-brevity}) on the conservative moment $\phi=O(1)$ (see the Appendix \ref{appendix-NACDEs-diffusive} for details) at the diffusive scaling,
\begin{subequations}
	\begin{align}
		&   \partial_t \phi +\partial_{\alpha}B_{\alpha} -\Delta t \frac{\partial}{\partial x_{\beta}}\Bigg[\chi \bigg(S^{1}_{\beta\gamma}+\Big(\frac{b}{2a^2}-1\Big)\delta_{\beta\gamma}\bigg)c_s^2\frac{\partial D_{\gamma \theta}}{ \partial x_{\theta}}\Bigg] -R=O(\Delta x) , \\
		&\partial_t \phi +\partial_{\alpha}B_{\alpha} - 
		\Delta t\frac{\partial}{\partial x_{\beta}}\Bigg[\chi \bigg(S^{1}_{\beta\gamma}+\Big(\frac{b}{2a^2}-1\Big)\delta_{\beta\gamma}\bigg)c_s^2\frac{\partial D_{\gamma \theta}}{ \partial x_{\theta}}\Bigg]-R =O(\Delta x^2) ,\label{second-diffusive-NACDEs}
	\end{align}
\end{subequations}
and additionally, from the Parts \ref{first-acoustic} and \ref{second-acoustic}, one can derive  the first- and second-order MEs at the acoustic scaling (see the Appendix \ref{appendix-NCDEs-acoustic} for details),
\begin{subequations}
	\begin{align}
		&   \partial_t \phi +\partial_{\alpha}B_{\alpha}  -R=O(\Delta x) , \\
		&\partial_t \phi +\partial_{\alpha}B_{\alpha} - 
		\Delta t\frac{\partial}{\partial x_{\beta}}\Bigg[\chi \bigg(S^{1}_{\beta\gamma}+\Big(\frac{b}{2a^2}-1\Big)\delta_{\beta\gamma}\bigg)c_s^2\frac{\partial D_{\gamma \theta}}{ \partial x_{\theta}}\Bigg]-R =O(\Delta x^2) .\label{second-acoustic-NACDEs}
	\end{align}
\end{subequations}
It is clear that  Eqs. (\ref{second-diffusive-NACDEs}) and (\ref{second-acoustic-NACDEs}) are consistent with the NACDE (\ref{NACDE}) with
\begin{equation}\label{kappa}
	\bm{\kappa}=\chi c_s^2\Delta t\bigg[\bm{S}^{1}+\Big(\frac{b}{2a^2}-1\Big)\bm{I}\bigg].
\end{equation}
When   the MRT-LB model with the orthogonal moments and   D2Q9 lattice structure is considered, it is easy to show that the second-order modified equation (\ref{second-diffusive-NACDEs}) is consistent with the result in Ref. \cite{Zhao2019}.

\subsection{NSEs: MEs of the GPMRT-LB model and GPMFD scheme}
We now consider the following $d$-dimensional NSEs with a force term,
\begin{subequations}\label{NSEs}
	\begin{align}
		&\partial_t\rho+\nabla\cdot(\rho \bm{u})=0,\label{continuous-equation}\\
		&\partial_t(\rho \bm{u})+\nabla\cdot(\rho \bm{uu})=-\nabla p+\nabla\cdot\bm{\sigma}+\bm{\hat{F}},\label{momentum-equation}
	\end{align}
\end{subequations}
where $p=c_s^2\rho$ is the pressure, $\bm{\hat{F}}=\Big(\hat{F}_{x_{\alpha-1}}\Big)_{\alpha=2}^{d+1}$ is the force term, and $\bm{\sigma}$ is  the shear stress defined by 
\begin{align}
	\bm{\sigma}&=\mu\big[\nabla\bm{u}+(\nabla\bm{u})^T\big]+\Big(\mu_b-\frac{2\mu}{d}\Big)(\nabla\cdot\bm{u})\bm{I}\notag\\
	&=\mu\bigg[\nabla\bm{u}+(\nabla\bm{u})^T-\frac{2}{d}(\nabla\cdot\bm{u})\bm{I}\bigg]+\mu_b\big(\nabla\cdot\bm{u}\big)\bm{I},
\end{align} 
here $\mu$ and $\mu_b$ are the dynamic and bulk viscosity, respectively.

To recover the macroscopic NSEs (\ref{NSEs}) from the GPMRT-LB model (\ref{GPMRT-LB-model}), the equilibrium, auxiliary, and source distribution functions, i.e., $f_i^{eq}$, $F_i$, and $G_i$, should satisfy the following moment conditions:
\begin{subequations}\label{NSEs-moment-conditions}
\begin{align}
	& \sum_{i=1}^qf_i=\sum_{i=1}^qf_i^{eq}=\rho,\sum_{i=1}^q\bm{c}_if_i=\sum_{i=1}^q\bm{c}_if_i^{eq}=\rho\bm{u},\\
	&\sum_{i=1}^q\bm{c}_i\bm{c}_if_i^{eq}=\rho\bm{uu}+ c_s^2\rho\bm{I},\sum_{i=1}^q\bm{c}_i\bm{c}_i\bm{c}_if_i^{eq}=c_s^2\rho\Delta \cdot\bm{u},\\
	&\sum_{i=1}^qG_i=0,\sum_{i=1}^q\bm{c}_iG_i=\bm{0},\sum_{i=1}^q\bm{c}_i\bm{c}_iG_i=\bm{0}, \\
	&\sum_{i=1}^qF_i=0,\sum_{i=1}^q\bm{c}_iF_i=\rho\bm{\hat{F}} ,
	\sum_{i=1}^q\bm{c}_i\bm{c}_iF_i=\rho\Big(\bm{\hat{F}u}+\big(\bm{\hat{F}u}\big)^T\Big), 
\end{align} 
\end{subequations}
for the D$d$Q$q$ lattice structure, the explicit expressions of $f^{eq}_i$, $G_i$, and $F_i$ can be given by \cite{Chai2020}
\begin{subequations}\label{NSEs-feq}
\begin{align}
	&f_i^{eq}=w_i\rho\Bigg[1+ \frac{\bm{c}_i\cdot\bm{u}}{c_s^2} + \frac{\bm{uu}:\big(\bm{c}_i\bm{c}_i-c_s^2\bm{I}\big)}{2c_s^4} \Bigg], \\
	&G_i=0, F_i=w_i\rho\biggg[ \frac{\bm{c}_i\cdot\bm{\hat{F}}}{c_s^2}+\frac{\Big(\bm{\hat{F}u}+\big(\bm{\hat{F}u}\big)^T\Big):\big(\bm{c}_i\bm{c}_i-c_s^2\bm{I}\big)}{2c_s^4}\biggg]. 
\end{align}
\end{subequations}
It should be noted that at the diffusive scaling,  the terms $w_i\rho \bm{c}_i\cdot\bm{u}/c_s^2$ and $w_i\rho \big[\bm{uu}:(\bm{c}_i\bm{c}_i-c_s^2\bm{I})\big]/(2c_s^4)$ in the expression of $f_i^{eq}$ are of order $O(\Delta x)$ and $O(\Delta x^2)$, respectively.

For the continuity equation (\ref{continuous-equation}), the derivation process is similar to that of NACDE (we refer the reader to Appendices \ref{appendix-NSEs-diffusive} and \ref{appendix-NSEs-acoustic} for details), and the first- to second-order MEs at the diffusive and acoustic scalings are given by 
\begin{align}\label{diff-con-second}
	&\Delta t\sim \Delta x^2:\:
	\left\{\begin{aligned}
		&   \partial_{\alpha}(\rho u_{\alpha})=O(\Delta x) , \\
		&\partial_t \rho +\partial_{\alpha}(\rho u_{\alpha}) =O(\Delta x^2) , 
	\end{aligned}\right.
\end{align} 
 
\begin{align}\label{acou-con-second}
	&\Delta t\sim \Delta x:\:
	\left\{\begin{aligned}
		&   \partial_t \rho +\partial_{\alpha}(\rho u_{\alpha})=O(\Delta x) , \\
		&\partial_t \rho +\partial_{\alpha}(\rho u_{\alpha})+\Delta t\Big(\frac{1}{2}-\frac{b}{2a^2}\Big)\partial_{\beta}\partial_{\theta}\big(\rho u_{\beta}u_{\theta}+\rho c_s^2\delta_{\beta\theta}\big) =O(\Delta x^2) . 
	\end{aligned}\right.
\end{align} 
We would like to point out that the term $\Delta t\big[1/2-b/(2a^2)\big]\partial_{\beta}\partial_{\theta}\big(\rho u_{\beta}u_{\theta}+\rho c_s^2\delta_{\beta\theta}\big) $ in Eq. (\ref{acou-con-second}) is of order $O(\Delta t Ma^2)$ with $Ma:=u/c_s$ being the Mach number, which can be eliminated  when $a^2=b$. However, this term  does not appear in Eq. (\ref{diff-con-second}) at the diffusive scaling, this is because $Ma$ is of order $O(\Delta x)$, in this case, it can be rearranged into the truncation error $O(\Delta x^2)$. This also indicates that the LB method is suitable  for nearly incompressible flows at both the diffusive and acoustic scalings.

With respect to the momentum equation (\ref{momentum-equation}), the first- and second-order MEs at the diffusive and acoustic scalings are given by (see Appendixes \ref{appendix-NSEs-diffusive} and \ref{appendix-NSEs-acoustic} for details)
\begin{align}\label{mon-second-diff}  
	&\Delta t\sim \Delta x^2:\:
	\left\{\begin{aligned}
		&   \partial_t(\rho u_{\alpha})+ \partial_{\beta}\big(\rho u_{\alpha}u_{\beta}+\rho c_s^2\delta_{\alpha\beta}) +\Delta t\Big(1-\frac{b}{2a^2}\Big)\partial_{\beta}\big(\rho c_s^2\partial_{\theta}u_{\theta}\delta_{\alpha\beta}\big) - \rho  {\hat{F}}_{x_{\alpha}} \\
		&\qquad -\Delta t \partial_{\beta}\bigg[S^{2}_{\alpha\beta\xi_1\xi_2}-\Big(1-\frac{b}{2a^2}\Big)\delta_{\xi_1\alpha}\delta_{\xi_2\beta}\bigg]\big(\rho c_s^2\partial_{\xi_1}u_{\xi_2}+\rho c_s^2\partial_{\xi_2}u_{\xi_1}\big) \\
		&\qquad=O(\Delta x ) , \\
		& \partial_t(\rho u_{\alpha})+\partial_{\beta}\big(\rho u_{\alpha}u_{\beta}+\rho c_s^2\delta_{\alpha\beta}\big)+ \Delta t \Big(1-\frac{b}{2a^2}\Big)\partial_{\beta}\big(\rho c_s^2\partial_{\theta}u_{\theta}\delta_{\alpha\beta}\big) - \rho  \hat{F}_{x_{\alpha}}  \\
		&\qquad-\Delta t \partial_{\beta}\bigg[S^{2}_{\alpha\beta\xi_1\xi_2}-\Big(1-\frac{b}{2a^2}\Big)\delta_{\xi_1\alpha}\delta_{\xi_2\beta}\bigg]\big(\rho c_s^2\partial_{\xi_1}u_{\xi_2}+\rho c_s^2\partial_{\xi_2}u_{\xi_1}\big)   \\
		&\qquad =O(\Delta x^2) , 
	\end{aligned}\right.
\end{align} 
and 
\begin{align} \label{mon-second-acou}
	&\Delta t\sim \Delta x:\:
	\left\{\begin{aligned}
		&   \partial_t (\rho u_{\alpha})+\partial_{\beta}\big(\rho u_{\alpha}u_{\beta}+\rho c_s^2\delta_{\alpha\beta}\big)-\rho \hat{F}_{x_{\alpha}}=O(\Delta x) , \\
		& \partial_t(\rho u_{\alpha})+\partial_{\beta}\big(\rho u_{\alpha}u_{\beta}+\rho c_s^2\delta_{\alpha\beta}\big)+\frac{\Delta t}{2 }\Big(1-\frac{b}{a^2}\Big)\partial_{\beta}\big(\rho c_s^2\partial_{\theta}u_{\theta}\delta_{\alpha\beta}\big) - \rho  \hat{F}_{x_{\alpha}}  \\
		&\qquad\qquad-\Delta t \partial_{\beta}\bigg[S^{2}_{\alpha\beta\xi_1\xi_2}-\Big(1-\frac{b}{2a^2}\Big)\delta_{\xi_1\alpha}\delta_{\xi_2\beta}\bigg]\big(\rho c_s^2\partial_{\xi_1}u_{\xi_2}+\rho c_s^2\partial_{\xi_2}u_{\xi_1}\big) \\
		&\qquad\qquad =O(\Delta x^2+\Delta xMa^3 ),
	\end{aligned}\right.
\end{align} 
where  
\begin{subequations}
	\begin{align}
		&\partial_t(\rho u_{\xi_1}u_{\xi_2})=\hat{F}_{x_{\xi_1}}u_{\xi_2}+\hat{F}_{x_{\xi_2}}u_{\xi_1}-c_s^2\big(u_{\xi_1}\partial_{\xi_2}\rho+u_{\xi_2}\partial_{\xi_1}\rho\big)+O(\Delta x+Ma^3),\\
		&\partial_{\theta}(\rho c_s^2\Delta_{\alpha\beta\theta\zeta}u_{\zeta})= \partial_{\theta}(\rho c_s^2u_{\theta}\delta_{\alpha\beta})+ \partial_{\alpha}(\rho c_s^2u_{\beta})+\partial_{\beta}(\rho c_s^2u_{\alpha})   ,
	\end{align} 
\end{subequations}
$\partial_t (\rho \bm{u}\bm{u}\bm{u})=O(Ma^3)$, $\bm{u}=O(Ma)$ and $\nabla\rho=O(Ma^2)$ are used at the acoustic scaling. From Eqs. (\ref{mon-second-diff}) and (\ref{mon-second-acou}), one can see that the viscous terms are different at the diffusive and acoustic scalings. This is because  the terms $\Delta t\partial_{tt}\hat{m}_{\alpha+1}^{eq}$ and $\Delta t\big[\bm{\mathcal{W}}_0\partial_t\bm{\hat{m}}^{eq}\big]_{\alpha+1}$ $\big(\alpha\in\{1\sim d\}\big)$, where $\bm{\hat{m}}^{eq}:=\bm{m}^{eq}/\Delta x^k=O(1)$ with $k=1$ at the diffusive scaling while $k=0$ at the acoustic scaling, are of order $O(\Delta x^2)$ at the diffusive scaling while are of order $O(\Delta x)$ term at the acoustic scaling. 

For the term $S^2_{\alpha\beta\xi_1\xi_2}$ shown in Eqs. (\ref{mon-second-diff}) and (\ref{mon-second-acou}), we consider the following commonly used form \cite{Chai2020}:
	\begin{equation}
S^2_{\alpha\beta\xi_1\xi_2}=\begin{cases}
		\frac{1}{s_{2s}}+\frac{1}{d}\Big(\frac{1}{s_{2b}}-\frac{1}{s_{2s}}\Big),&\xi_1=\xi_2=\alpha=\beta, \\ 
		\frac{1}{d}\Big(\frac{1}{s_{2b}}-\frac{1}{s_{2s}}\Big),&\alpha=\beta,\xi_1=\xi_2,\xi_1\neq \alpha,\\
		\frac{1}{s_{2s}} ,&\xi_1=\alpha,\xi_2=\beta,\alpha\neq \beta, \\
		0,&{\rm others}. \\	
	\end{cases}
\end{equation}
Then the second-order ME (\ref{mon-second-diff}) at the diffusive scaling is the same as the momentum euqation (\ref{momentum-equation}) with \begin{equation}
	\upsilon=\Big(\frac{1}{s_{2s}}+\frac{b}{2a^2}-1\Big) c_s^2\Delta t,\upsilon_b=\bigg[\frac{2}{d}\Big(\frac{1}{s_{2b}}+\frac{b}{2a^2}-1\Big)+\Big(\frac{b}{2a^2}-1\Big)\bigg]  c_s^2\Delta t,\mu=\rho\upsilon,\mu_b=\rho\upsilon_b,
\end{equation}
while the second-order ME (\ref{mon-second-acou}) at the acoustic scaling would reduce to the momentum euqation (\ref{momentum-equation}) under the following condition, \begin{equation}
	\upsilon=\Big(\frac{1}{s_{2s}}+\frac{b}{2a^2}-1\Big) c_s^2\Delta t,\upsilon_b=\bigg[\frac{2}{d}\Big(\frac{1}{s_{2b}}+\frac{b}{2a^2}-1\Big)+\Big(\frac{b}{2a^2}-\frac{1}{2}\Big)\bigg] c_s^2\Delta t,\mu=\rho\upsilon,\mu_b=\rho\upsilon_b.
\end{equation}

\section{Fourth-order GPMRT-LB model and GPMFD scheme for one-dimensional CDE}\label{Fourth-Order-for-CDE}
Based on the above results, the GPMFD scheme  of a given GPMRT-LB model can be directly derived from Eq. (\ref{det-GPMFD-N-geq1-brevity}), then one can further conduct the accuracy and stability analysis with the help of the traditional tools adopted in the finite-difference method.    In this section, we will present the fourth-order GPMRT-LB model and GPMFD scheme at the diffusive scaling for the one-dimensional CDE,	\begin{equation}\label{CDE}
	\frac{\partial \phi}{\partial t}+u\frac{\partial \phi}{\partial x}=\kappa\frac{\partial ^2\phi}{\partial x^2},
\end{equation}
where $\phi$ is a scalar function of the position $x$ and time $t$, $u$ and $\kappa$ are two constants. For CDE (\ref{CDE}), the evolution equation of the GPMRT-LB model can be written as \cite{Chen2023, Guo2018} 
\begin{align}\label{GPMRT-LB-CDE}
	&f_i^{\star}(x,t)=f_i(x,t)-\Big(\bm{M}^{-1}\bm{SM}\Big)_{ik}\big[f_k-f_k^{eq}\big](x,t),\\
	&f_i(x,t+\Delta t)=p_0f_i^{\star}(x,t)+p_{-1}f_i^{\star}(x-\lambda_i\Delta t,t)+p_1f_i^{\star}(x+\lambda_i\Delta t,t),i=-1,0,1,
\end{align} 
where the D1Q3 lattice structure is considered, the orthogonal transform matrix $\bm{M}$ and diagonal relaxation matrix $\bm{S}$ can be expressed as
\begin{align}
	&\lambda_0=0,\lambda_1=\lambda,\lambda_{-1}=-\lambda, \bm{M}=\left(\begin{matrix}
		1&1&1\\
		0&c&-c\\
		-2c^2&c^2&c^2\\
	\end{matrix}\right),\bm{S}=\left(\begin{matrix}
		s_0&0&0\\
		0&s_1&0\\
		0&0&s_2\\
	\end{matrix}\right).
\end{align}
Here  the diagonal element $s_i$ of the relaxation matrix $\bm{S}$ is the relaxation parameter corresponding to the $i_{th}$ moment of the distribution function $f_i$.  

To derive the correct CDE (\ref{CDE}), the equilibrium distribution function should satisfy the following moment conditions,
\begin{align}
	&\sum_if_i=\sum_if_i^{eq}=\phi, \sum_ic_if_i^{eq}=\phi u, \sum_ic_ic_if_i^{eq}=(c_s^2+u^2)\phi,
\end{align}
from which one can determine  the equilibrium distribution as
\begin{equation}
	f_i^{eq}=w_i\phi\Big[1+\frac{c_iu}{c_s^2}+\vartheta\frac{u^2(c_i^2-c_s^2)}{2c_s^4}\Big],
\end{equation}
where 
\begin{align}
	&c_s^2=(1-w_0)c^2,w_1=w_{-1}=\frac{1-w_0}{2},
\end{align}
$w_i$ is the weight coefficient,  $\vartheta=\zeta\xi$ with $\zeta=2(1-w_0)/w_0$ and $\xi=(1/s_1-1/2)/\big[1/s_1+b/(2a^2)-1\big]$ \cite{Chen2023}. 

With the help of the  scheme (\ref{GPMFD-N-eq1}), one can easily obtain the  GPMFD scheme on the variable $\phi$, 
\begin{align}\label{fourth-scheme}
	\phi^{n+1}_j=&\alpha_1\phi^n_j+\alpha_2\phi^n_{j-1}+\alpha_3\phi^{n}_{j+1} +\beta_1\phi^{n-1}_j+\beta_2\phi^{n-1}_{j-1}+\beta_3\phi^{n-1}_{j+1}+\beta_4\phi^{n-1}_{j-1}+\beta_5\phi^{n-1}_{j+2}\notag\\
	&+\gamma_1\phi^{n-2}_j+\gamma^{n-2}_{j-1}+\gamma_3\phi^{n-2}_{j+1}+\gamma_4\phi^{n-2}_{j-2}+\gamma_5\phi^{n-2}_{j+2},
\end{align} 
where $\phi_j^n$ represnts $\phi(j\Delta x,t_n)$, $j\in\mathcal{Z}$ and $n\in\mathcal{N}$, the parameters $\alpha_i(i\in\{1\sim 3\}),\beta_k$ and $\gamma_k(k\in\{1\sim 5\})$ can be found in Appendix \ref{coe-fourth-order-fd-scheme}.

Due to the fact that the accuracy analysis on the above scheme (\ref{fourth-scheme}) is similar to our previous work \cite{Chen2023}, here we only present some results, and do not show more details.  Actually, the second-order ME of the GPMFD scheme (\ref{fourth-scheme}) can be given by
\begin{align}
	\Big[\frac{\partial\phi}{\partial t}+u\frac{\partial \phi}{\partial x}\Big]^n_j=&\varepsilon\frac{\Delta x^2}{\Delta t}\Big[\frac{\partial ^2\phi}{\partial x^2}\Big]^n_j-\frac{uTR_3}{6s_1}\Delta x^2\Big[\frac{\partial ^3\phi}{\partial x^3}\Big]^n_j\notag\\
	&+\frac{TR_4(1-w_0)}{24s_1s_2}\frac{\Delta x^4}{\Delta t}\Big[\frac{\partial ^4\phi}{\partial x^4}\Big]^n_j+O(\Delta t^2),
\end{align}
where $\varepsilon=\kappa\Delta t/\Delta x^2$. To derive a fourth-order GPMRT-LB model (\ref{GPMRT-LB-CDE}) and (\ref{fourth-scheme}), the following conditions need to be satisfied,
\begin{subequations}\label{fourth-condition}
	\begin{align}
		&\varepsilon=a^2\Big(\frac{1}{s_1}+\frac{b}{2a^2}-1\Big)(1-w_0) , \\		
		&TR_3=s_1^2s_2(1-3b)+12a^2s_2(s_1-1)\notag\\
		&\qquad+3w_0\Big[2a^2(s_1+2s_2+s_1^2(s_2-1)-3s_1s_2)+
		bs_1\big(s_1(1-s_2)+s_2\big)\Big]=0, \\
		&TR_4= 6a^4s_2(6s_1-4-2s_1^3) + bs_1^3s_2(1-3b)+ 8a^2s_1^2s_2(1-s_1)(1-3b) \notag
		\\
		&\qquad  + 6a^4w_0\Big[4\big(s_1+s_2+s_1^3+2s_1^2(s_2-1)\big)-10s_1s_2-2s_1^3s_2\Big] 
		+ 3b^2s_1^3w_0(2-s_2) \notag \\
		&\qquad + 12a^2bs_1w_0\Big[2s_1(1-s_1)+s_2\big((s_1-2)s_1+1\big)\Big]=0,
	\end{align}
\end{subequations}
where $u/c=O(\Delta x)$ has been used, while it has not been take into account in the previous work  \cite{Chen2023}. In addition, for the special case with $a=b=1$, the solution of the fourth-order conditions (\ref{fourth-condition}) can be derived,
\begin{align}
	&
	\left\{\begin{aligned}
		&  s_1=\frac{12\varepsilon}{6\varepsilon+1}, \\
		&s_2=\frac{2}{6\varepsilon +1},\\
		&w_0=1-12\varepsilon^2,
	\end{aligned}\right. 
\end{align}
which is different from that in Ref. \cite{Chen2023}.

\begin{figure}[htbp]    
		\centering            
		\subfloat[$a=1,b=1$]    
		{
			\label{fig:subfig1}\includegraphics[width=0.4\textwidth]{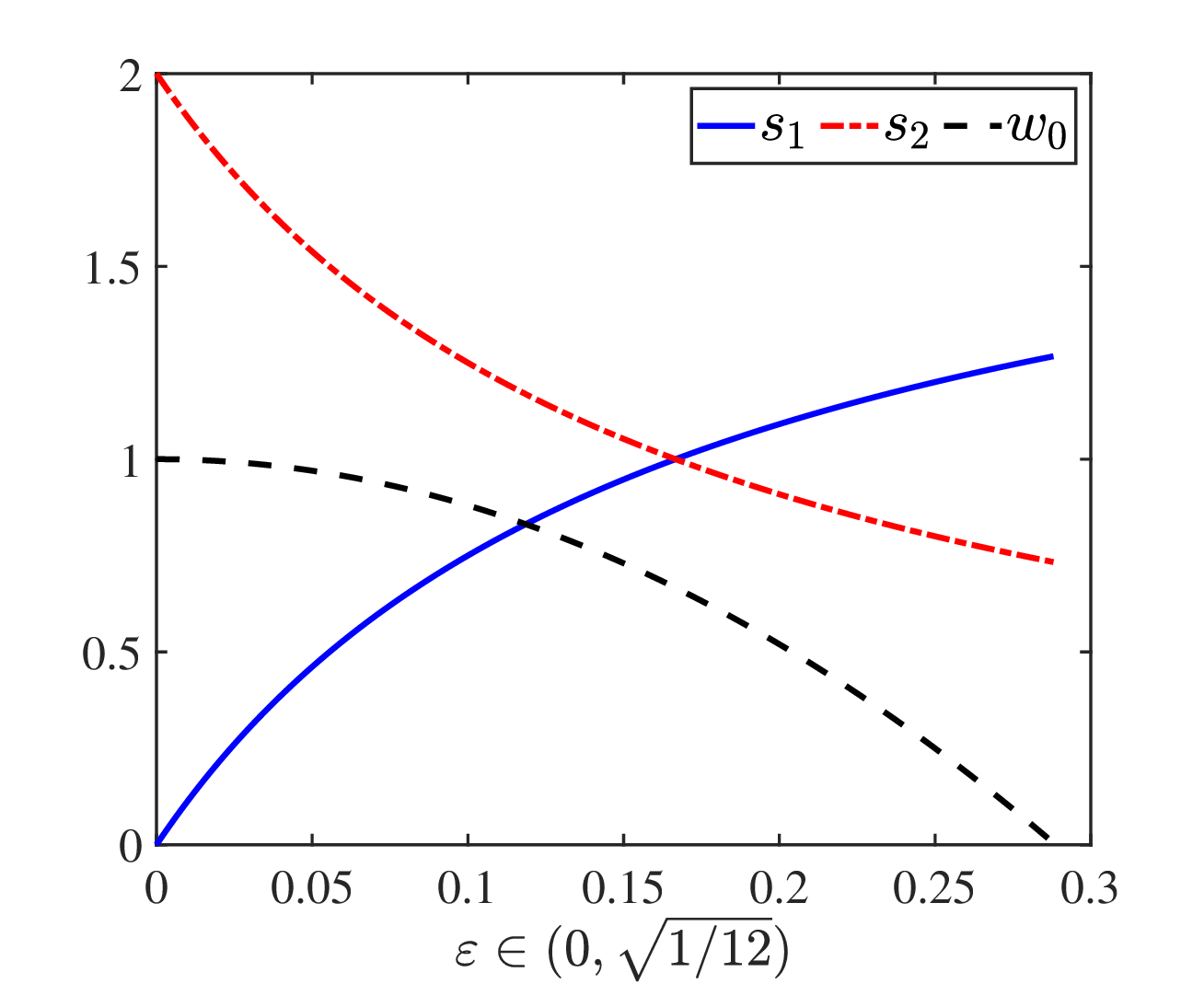} 
		}
		\subfloat[$a=0.6,b=0.9$]
		{
			\label{fig:subfig2}\includegraphics[width=0.4\textwidth]{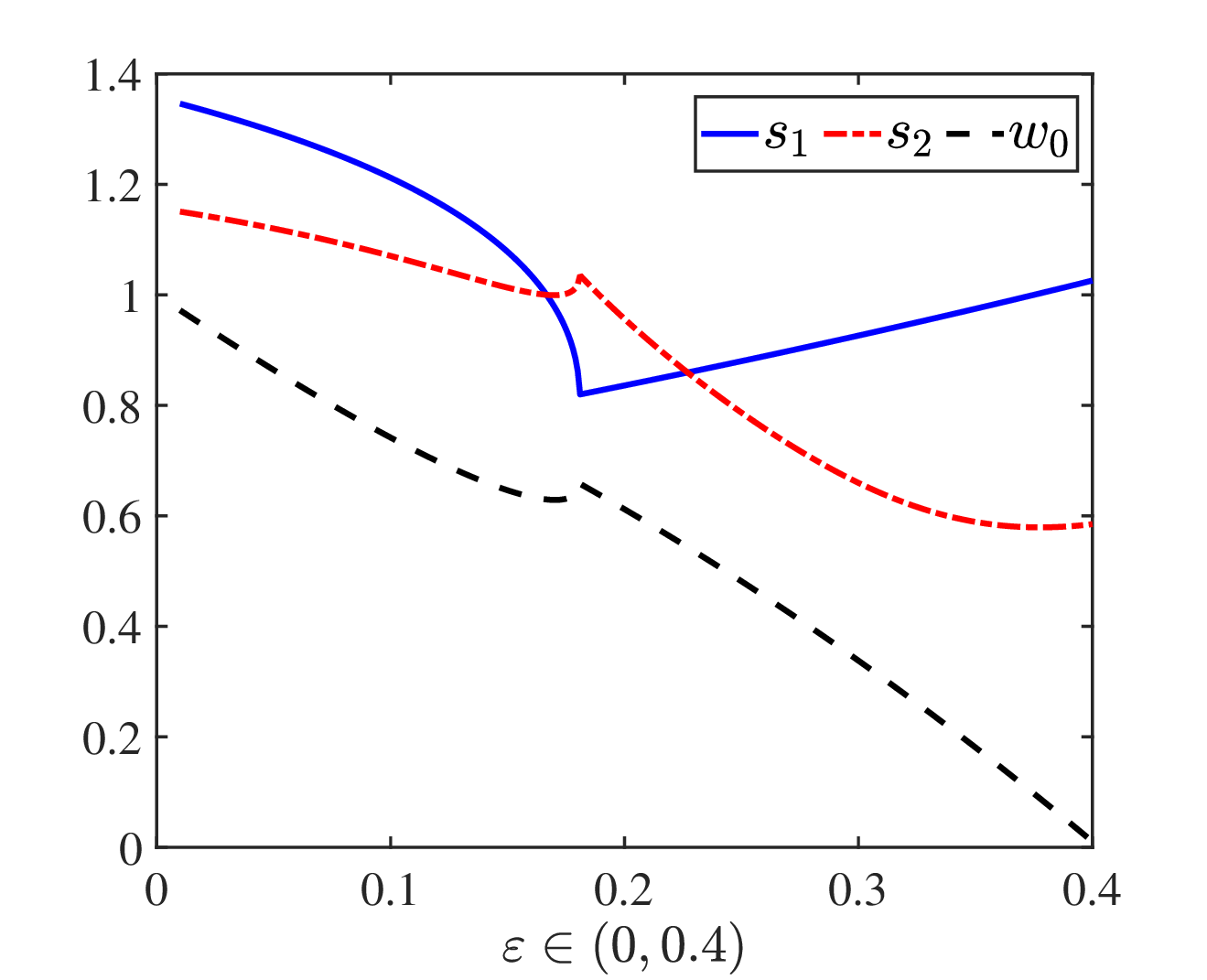}
		}
		\quad
		\subfloat[$a=0.9,b=0.8$]
		{
			\label{fig:subfig2}\includegraphics[width=0.4\textwidth]{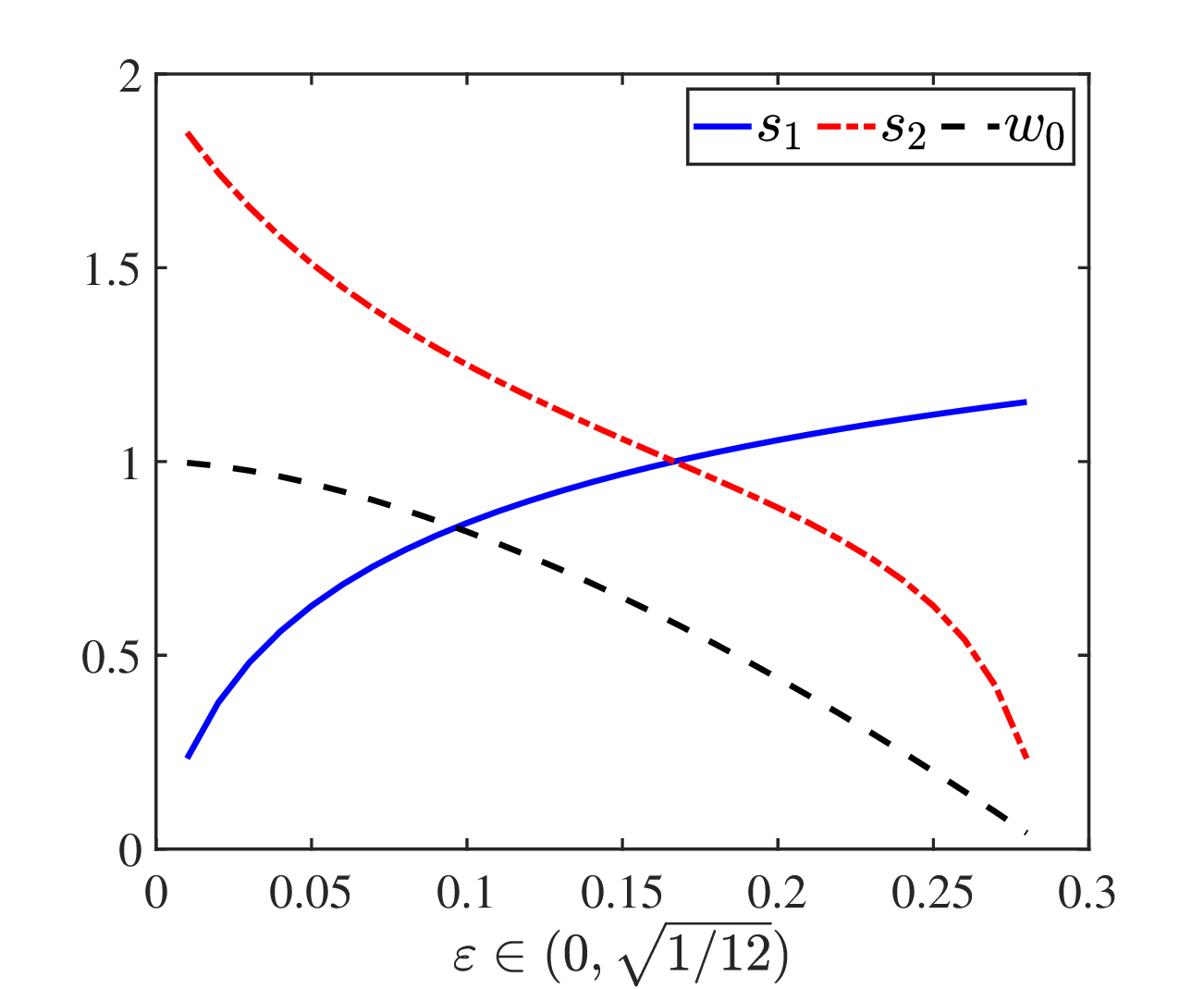}
		}
		\subfloat[$a=0.6,b=0.6$]
		{
			\label{fig:subfig2}\includegraphics[width=0.4\textwidth]{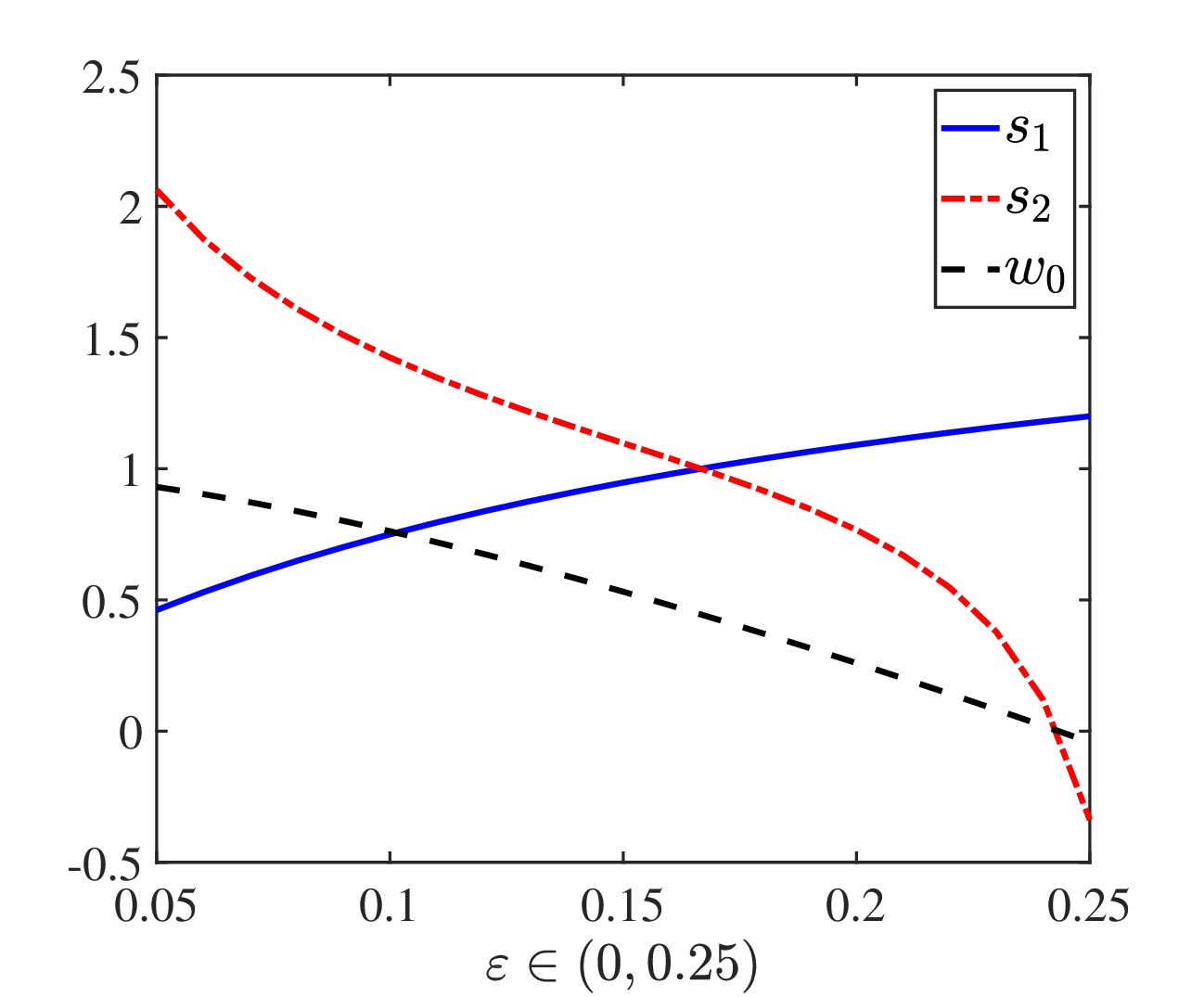}
		}
		\caption{The weight coefficient $w_0$, the relaxation parameters $s_1$ and $s_2$ as a function of the parameter $\varepsilon$ under different values of $a$ and $b$. }    
		\label{relation-w0s1s2}            
	\end{figure}

\begin{figure}[htbp]    
			\centering            
			\subfloat[$a=1,b=1$]    
			{
				\label{fig:subfig1}\includegraphics[width=0.4\textwidth]{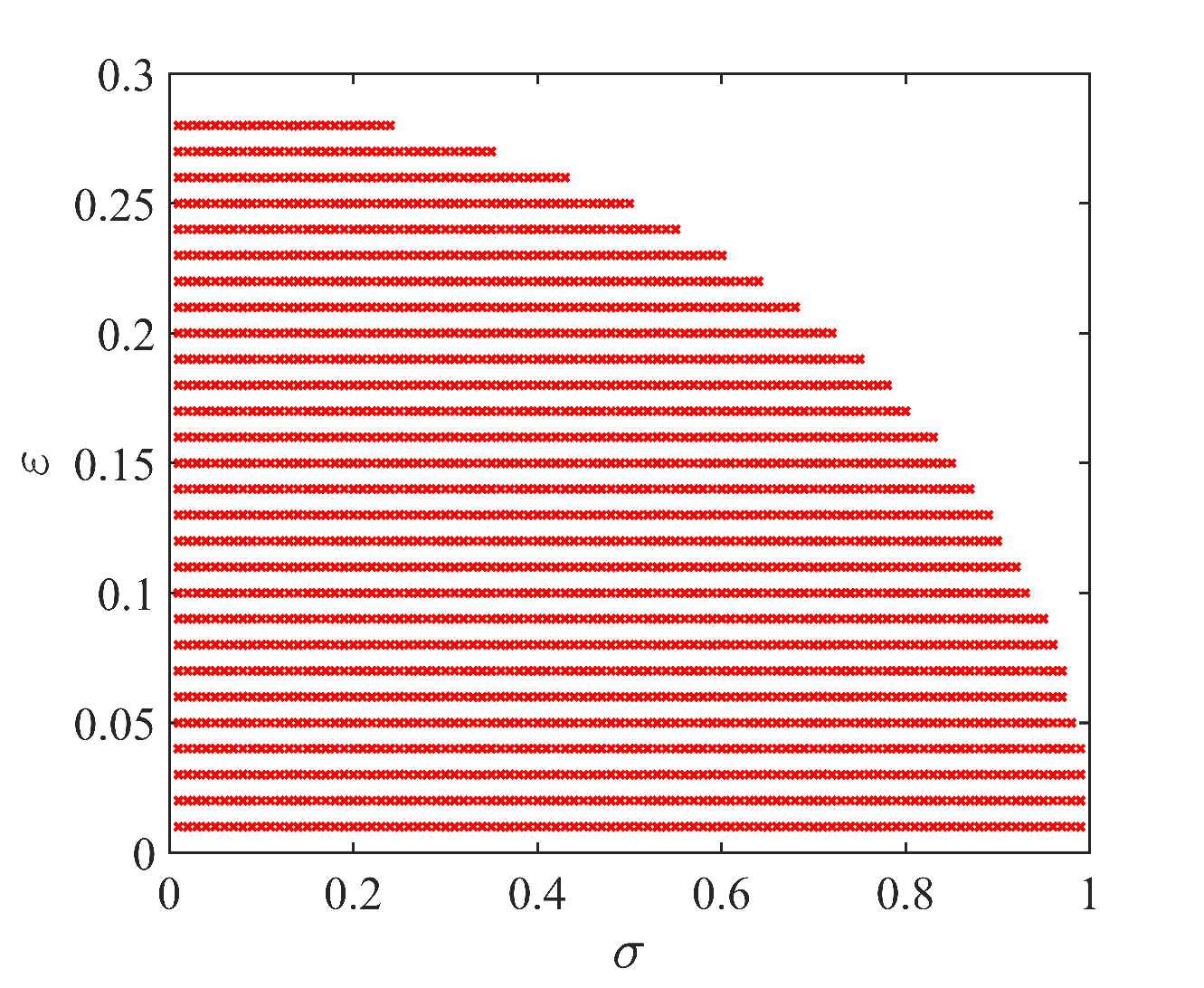} 
			}
			\subfloat[$a=0.6,b=0.9$]
			{
				\label{fig:subfig2}\includegraphics[width=0.4\textwidth]{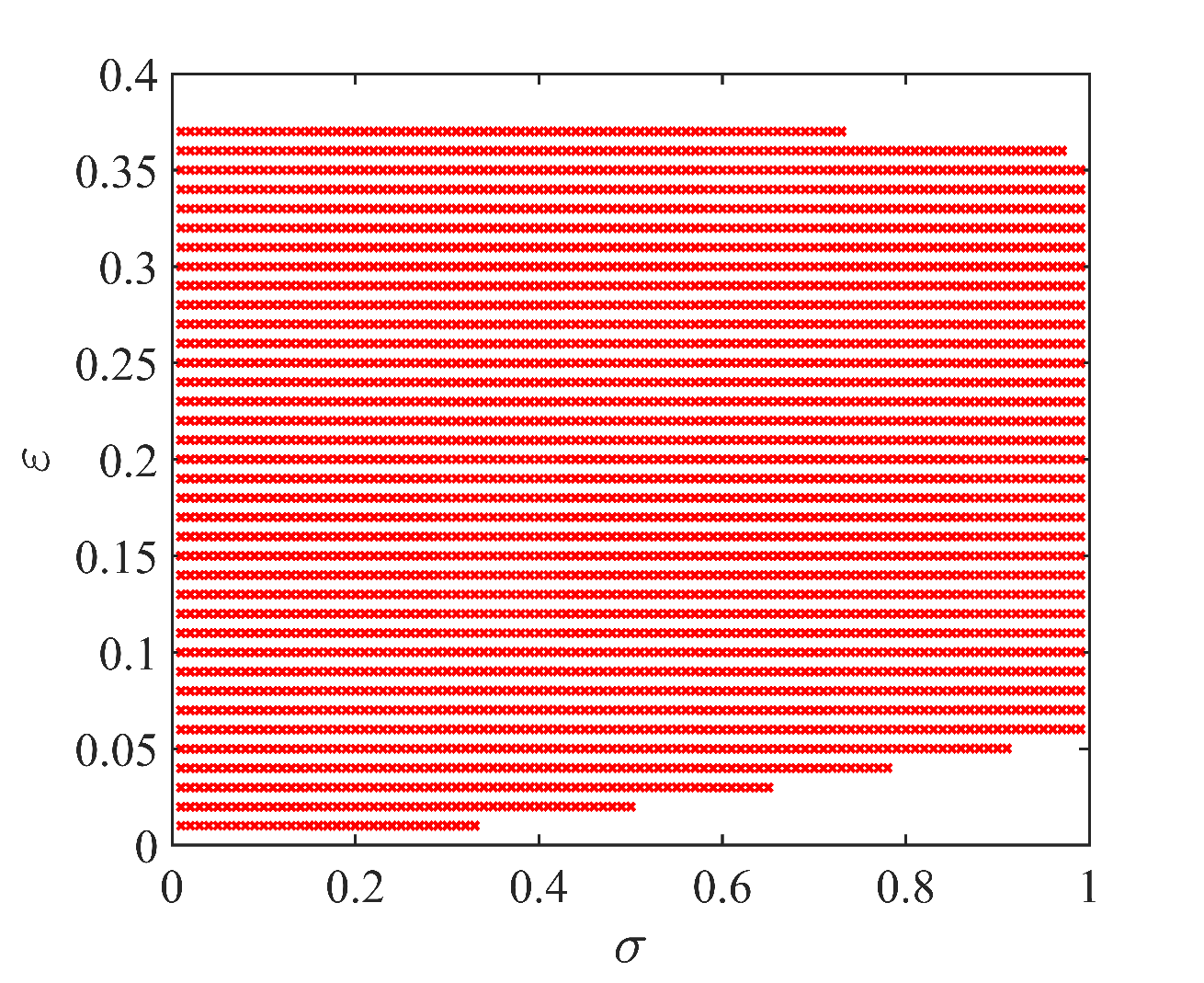}
			}
			\quad
			\subfloat[$a=0.9,b=0.8$]
			{
				\label{fig:subfig2}\includegraphics[width=0.4\textwidth]{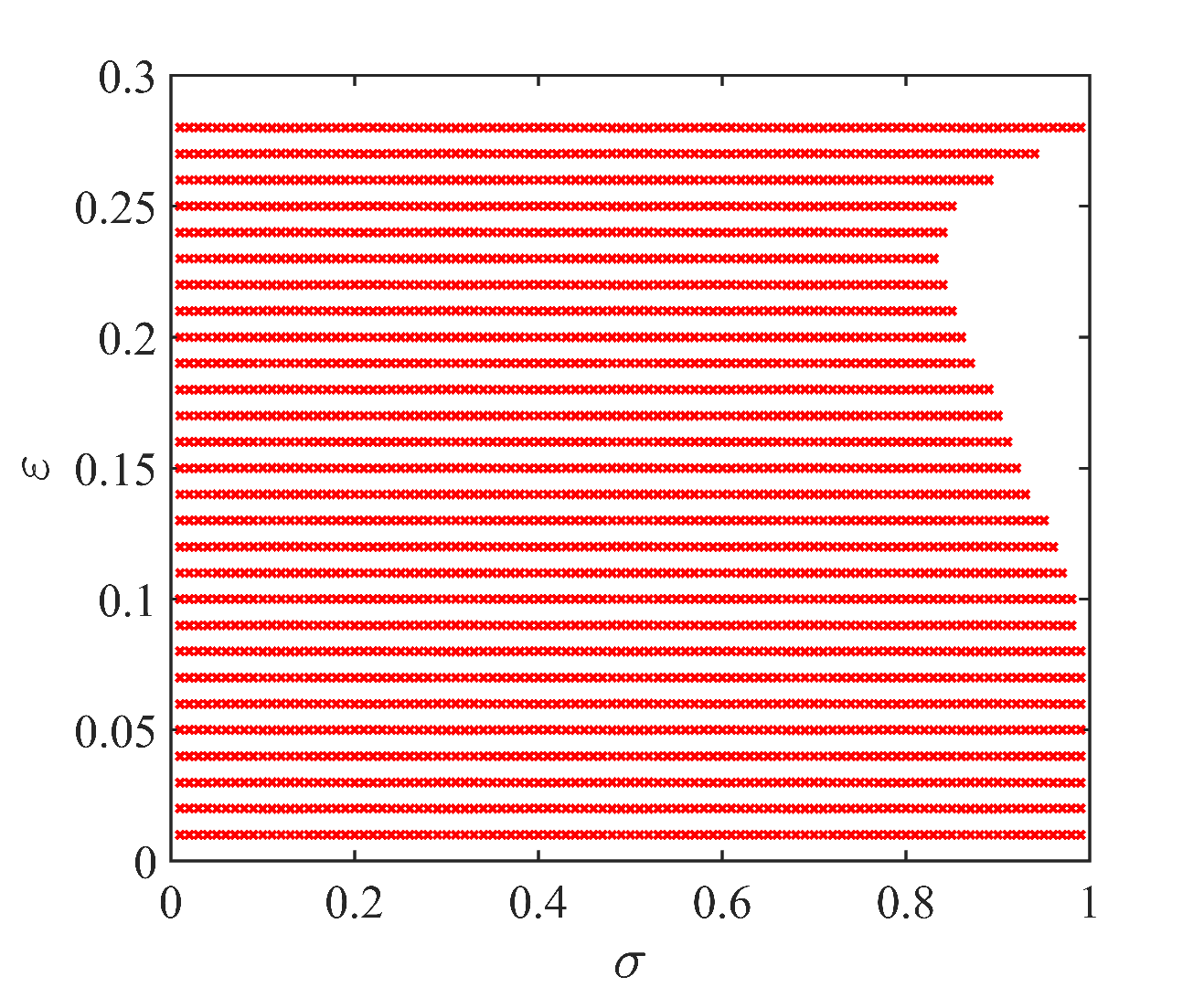}
			}
			\subfloat[$a=0.6,b=0.6$]
			{
				\label{fig:subfig2}\includegraphics[width=0.4\textwidth]{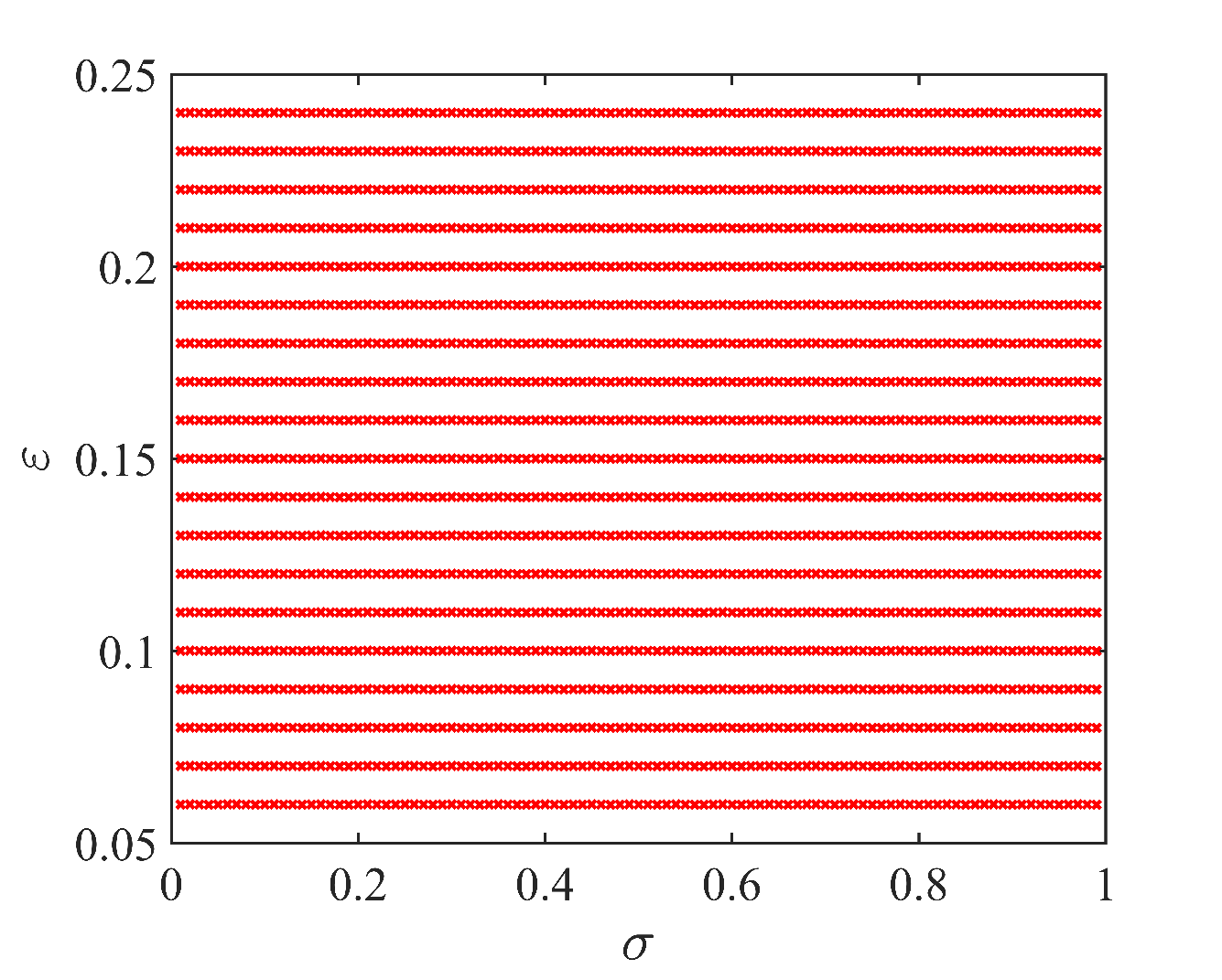}
			}
			\caption{The stability regions of F-GPMRT-LB model and F-GPMFD scheme for CDE (\ref{CDE}) under different values of $a$ and $b$. }    
			\label{stability-region}            
		\end{figure}

It should be noted that the F-GPMRT-LB model and F-GPMFD scheme can be obtained once the weight coefficient $w_0$, relaxation parameters $s_1$ and $s_2$ satisfy the fourth-order conditions (\ref{fourth-condition}) for the given parameters $\varepsilon$, $a$, and $b$. However, owing  to existence of nonlinearity and coupling, it is difficult to derive analytical  solution of the  fourth-order conditions (\ref{fourth-condition}), thus here we only plot the relation between the parameters $s_1$, $s_2$, $w_0$ and  $\varepsilon$ through selecting four cases of parameters $a$ and $b$ in Fig. \ref{relation-w0s1s2}. Furthermore, based on the solution of Eq.  (\ref{fourth-condition}), we can  also determine the corresponding stability regions, as shown in Fig. \ref{stability-region}. From this  figure, one can observe  that stability regions  of the F-GPMRT-LB model and F-GPMFD scheme can be larger than that of the MRT-LB model  through adjusting parameters $a$ and $b$ properly.
	
\section{Numerical results and discussion}\label{Numerical-Simulation}
In this section, we conduct some numerical  experiments of the Gauss hill problem and Poiseuille flow  since they have the analytical solutions, which can also be used to test the convergence rates (CRs) of the GPMRT-LB model (\ref{GPMRT-LB-model}) and GPMFD scheme (\ref{det-GPMFD-N-geq1-brevity}) for NACDE and NSEs. To measure the difference between the numerical result and analytical solution, we adopt  the following root-mean-square error (RMSE) \cite{Kruger2017},
\begin{equation}
	{\rm RMSE}:=\sqrt{\frac{\sum_i[\psi(\bm{x}_i,t_n)-\psi^{\star}(\bm{x}_i,t_n)]^2}
		{\prod_{j=1}^dN_{x_j}}},
\end{equation}
where $N_{x_j}$ is the number of gird points in the $j$ direction, $\bm{x}_i$ denots the grid point, $\psi$ and $\psi^{\star}$ are the numerical and analytical solutions, respectively. Based on the definition of RMSE, one can estimate the CR with the following formula:
\begin{equation}
	{\rm CR}=\frac{\log({\rm RMSE}_{\Delta x}/{\rm RMSE}_{\Delta x/2})}{\log 2}.
\end{equation}
Here we consider the popular D2Q9 lattice structure with $c_s=c/\sqrt{3}$  for two-dimensional problems. For the MRT-LB model, we adopt the orthogonal transform matrix,
\begin{equation}
	\bm{M}^c=\left(	\begin{matrix}
		1&1&1&1&1&1&1&1&1\\
		0&c&0&-c&0&c&-c&-c&c\\
		0&0&c&0&-c&c&c&-c&-c\\
		-4c^2&-c^2&-c^2&-c^2&-c^2&2c^2&2c^2&2c^2&2c^2\\
		0&c^2&-c^2&c^2&-c^2&0&0&0&0\\
		0&0&0&0&0&c^2&-c^2&c^2&-c^2\\	
		4c^4&-2c^4&-2c^4&-2c^4&-2c^4&c^4&c^4&c^4&c^4\\
		0&-2c^3&0&2c^3&0&c^3&-c^3&-c^3&c^3\\
		0&0&-2c^3&0&2c^3&c^3&c^3&-c^3&-c^3\
	\end{matrix}\right),
\end{equation}
the relaxation matrix $\bm{S}^c=$\textbf{diag}$\big(1,(\bm{S}^1)^{-1},1,1,1,1,1,1\big)$ [see $\bm{S}^1$ in Eq. (\ref{natural-condition})] is used  for the Gauss hill problem while $\bm{S}^c=$\textbf{diag}$\big(1,1,1,s_{2s},s_{2s},s_{2s},1,1,1\big)$ is applied for the Poiseuille flow. To satisfy the
stability condition \cite{Zhao2020}, 
the parameter $b$ is located in  the range $[a^2,\min\{1,6\upsilon+a^2\}]$. In addition, we also consider the one-dimensional CDE (\ref{CDE}) with the periodic boundary condition to test the F-GPMRT-LB model (\ref{GPMRT-LB-CDE}) and F-GPMFD scheme (\ref{fourth-scheme}), where the parameter $\varepsilon$, weight coefficient $w_0$, and relaxation parameters $s_1$ and $s_2$ are determined by Eq. (\ref{fourth-condition}) and the stability region  shown in Fig. \ref{stability-region}.\\
\textbf{Example 1} We first consider Gauss hill problem. 
With the following initial condition:
\begin{equation}
	\phi(\bm{x},0)=\frac{\phi_0}{2\phi\Upsilon_0^2}\exp\bigg\{\Big(-\frac{\bm{x}^2}{2\Upsilon_0^2}\Big)\bigg\},
\end{equation}
  one can obtain the analytical solution of this problem under the constant velocity $\bm{u}=(u_x,u_y)^T$ and diffusion coefficient matrix $\bm{\kappa}$,
\begin{align}
	\phi(\bm{x},t)=\frac{\phi_0}{2\pi{|\det{(\bm{\Upsilon})}|}^{1/2}}\exp\Big\{-\frac{\bm{\Upsilon}^{-1}:\big[(\bm{x}-\bm{u}t)(\bm{x}-\bm{u}t)\big]}{2}\Big\},
\end{align}
where $\bm{x}=(x,y)^T$, $\bm{\Upsilon}=\Upsilon_0^2\bm{I}+2\bm{\kappa}t$, and $\bm{\Upsilon}^{-1}$  $\det{(\bm{\Upsilon})}$ represent the inverse matrix and determinant value of $\bm{\Upsilon}$, respectively. In our simulations, the computational domain is $[-1,1]\times[-1,1]$ and the total concentration is set as $\phi_0=2\phi(\Upsilon_0)^2$ with $\Upsilon_0=0.01$, which should be small enough when applying the periodic boundary condition.

\begin{figure}[htbp]    
		\centering            
		\subfloat[GPMRT-LB model]    
		{
			\label{fig:subfig1}\includegraphics[width=0.4\textwidth]{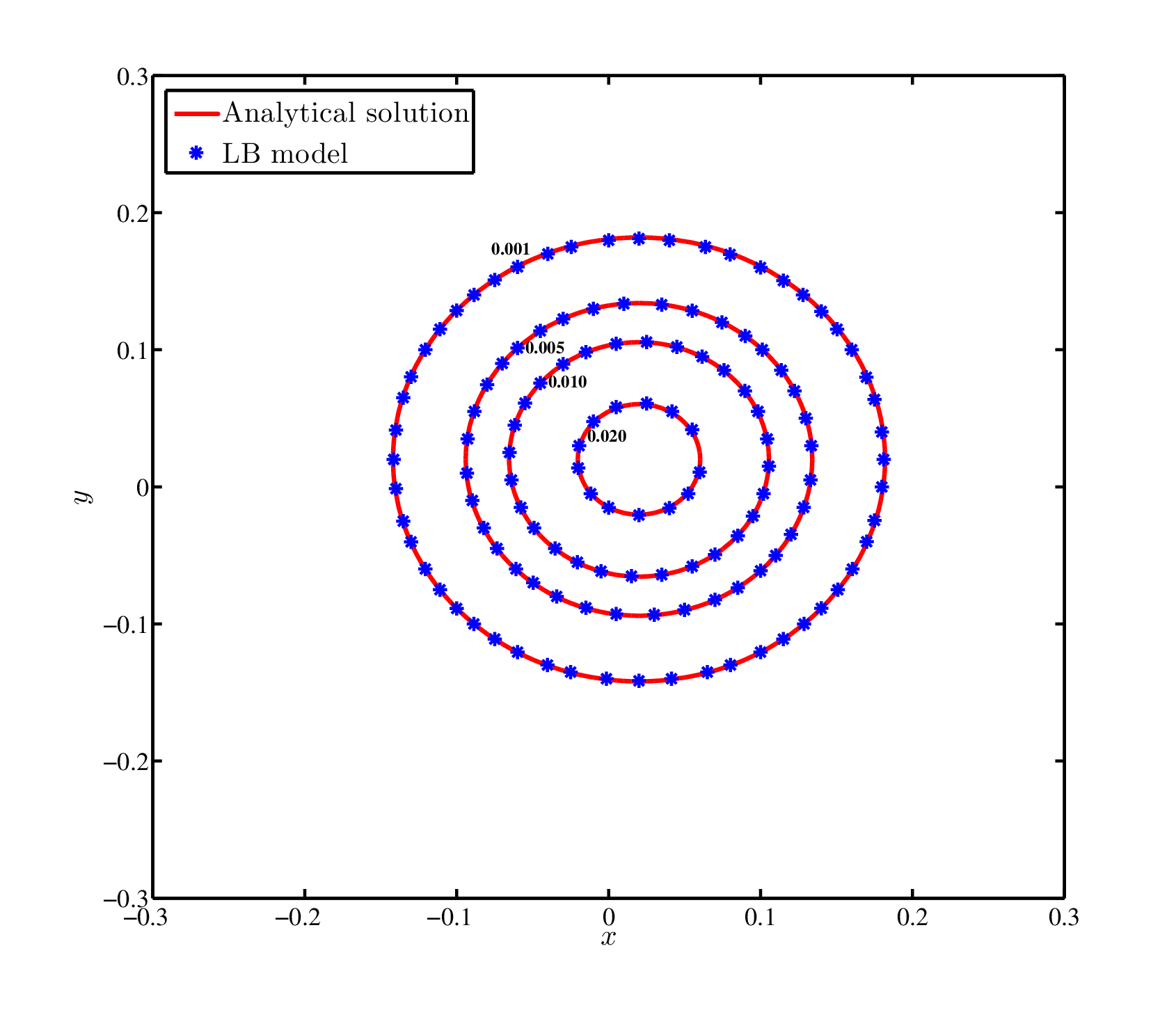} 
		}
		\subfloat[GPMFD scheme]
		{
			\label{fig:subfig2}\includegraphics[width=0.4\textwidth]{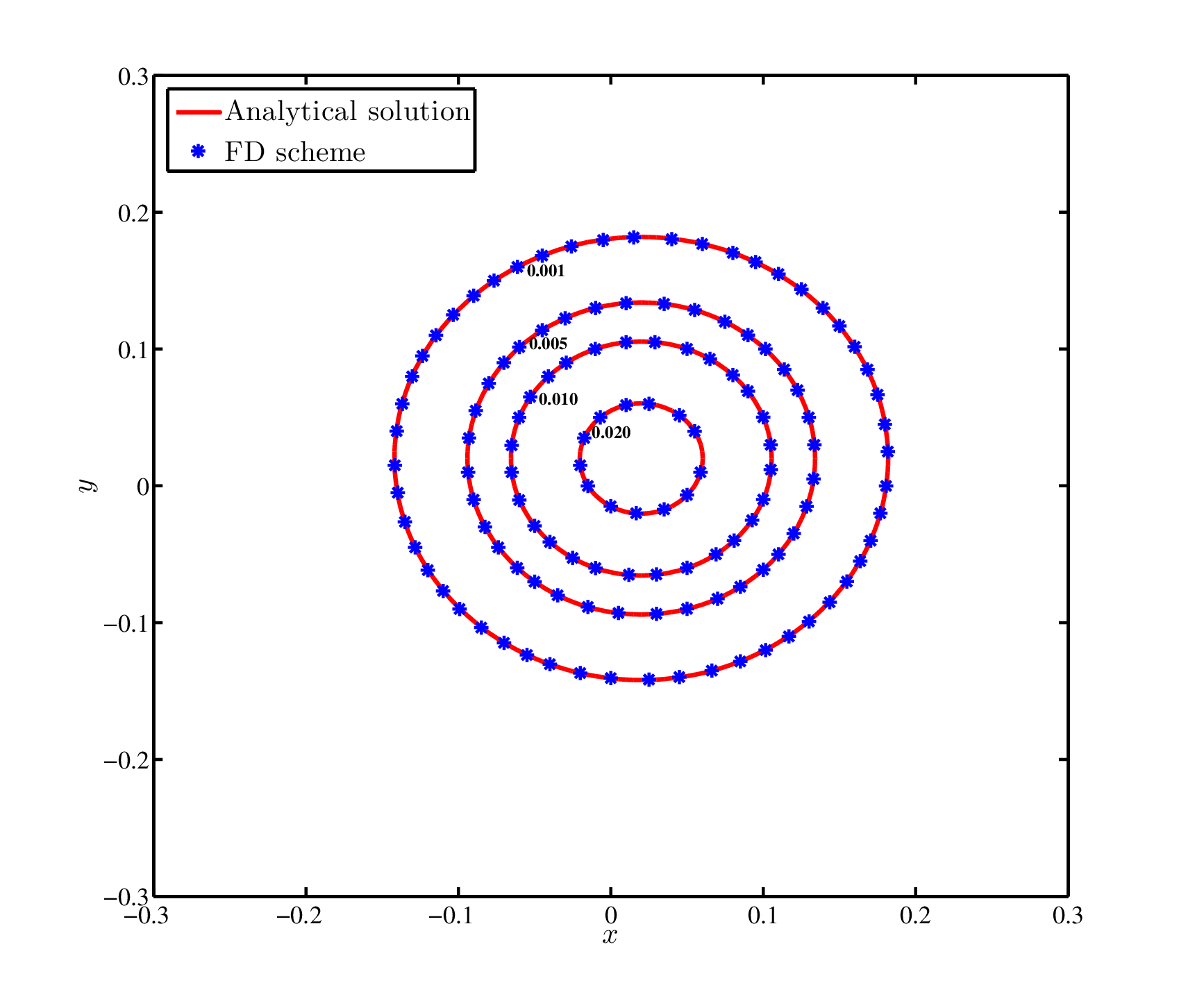}
		}
		\quad
		\subfloat[GPMRT-LB model]
		{
			\label{fig:subfig2}\includegraphics[width=0.4\textwidth]{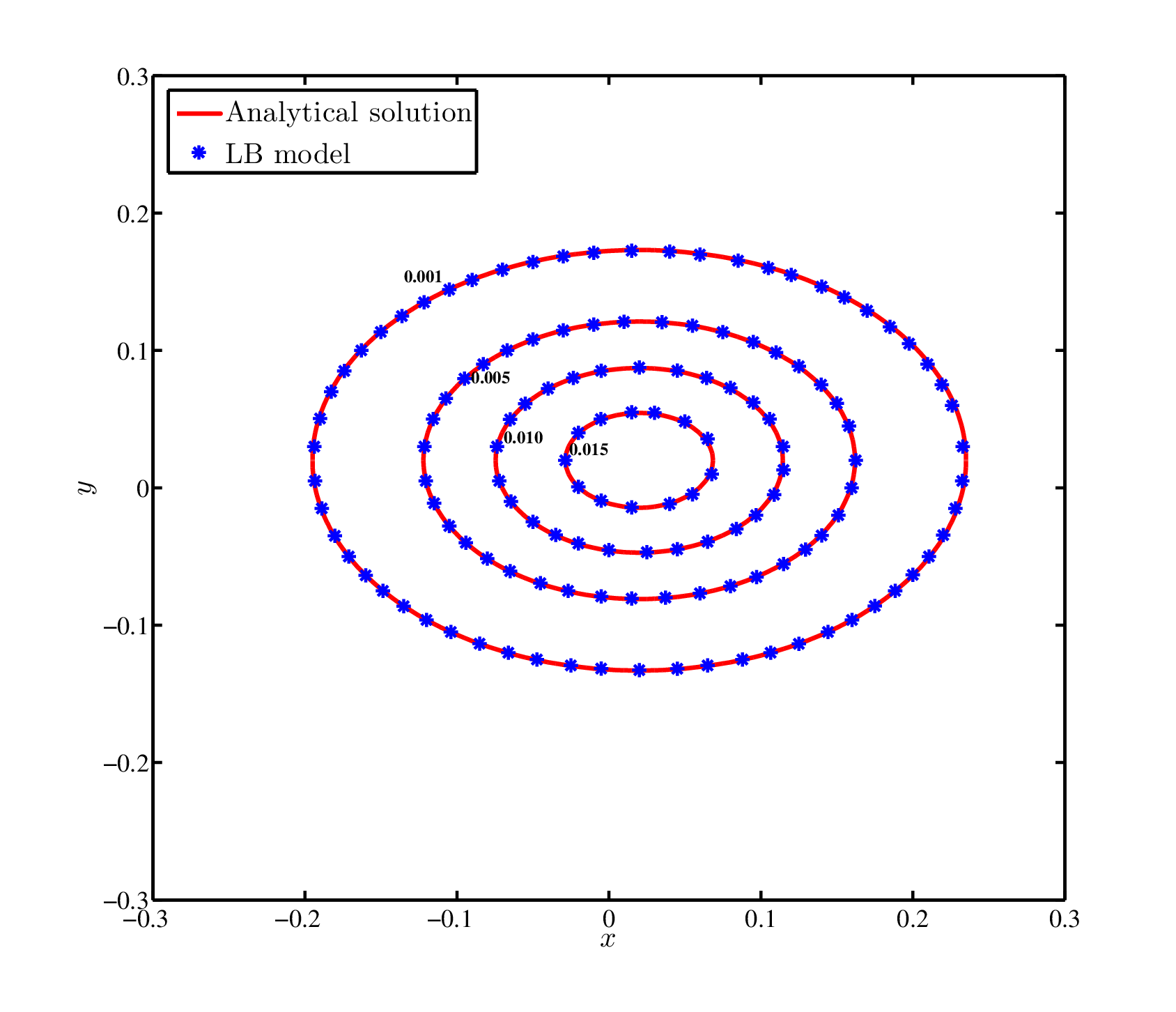}
		}
		\subfloat[GPMFD scheme]
		{
			\label{fig:subfig2}\includegraphics[width=0.4\textwidth]{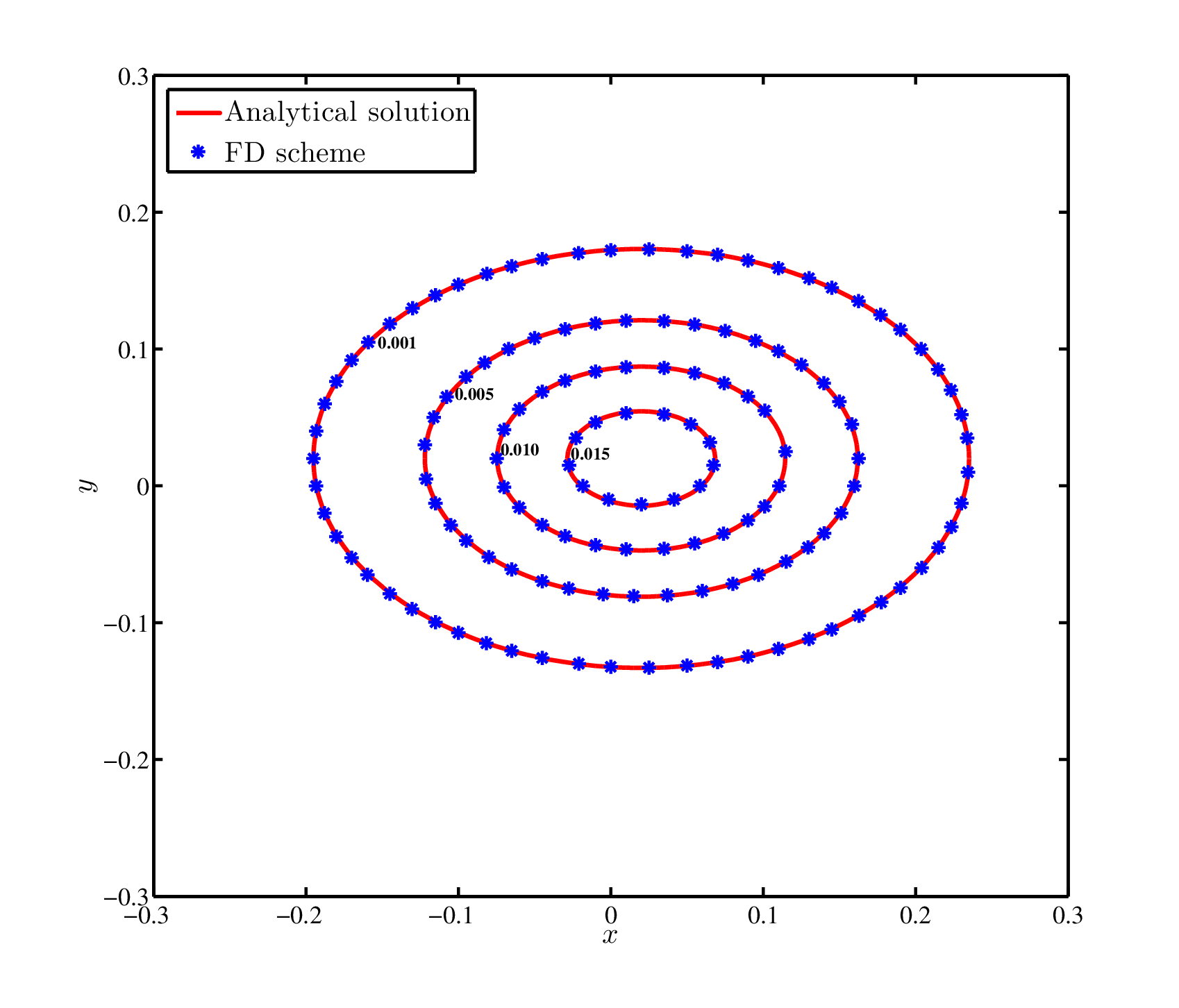}
		}
		\quad
		\subfloat[GPMRT-LB model]
		{
			\label{fig:subfig2}\includegraphics[width=0.4\textwidth]{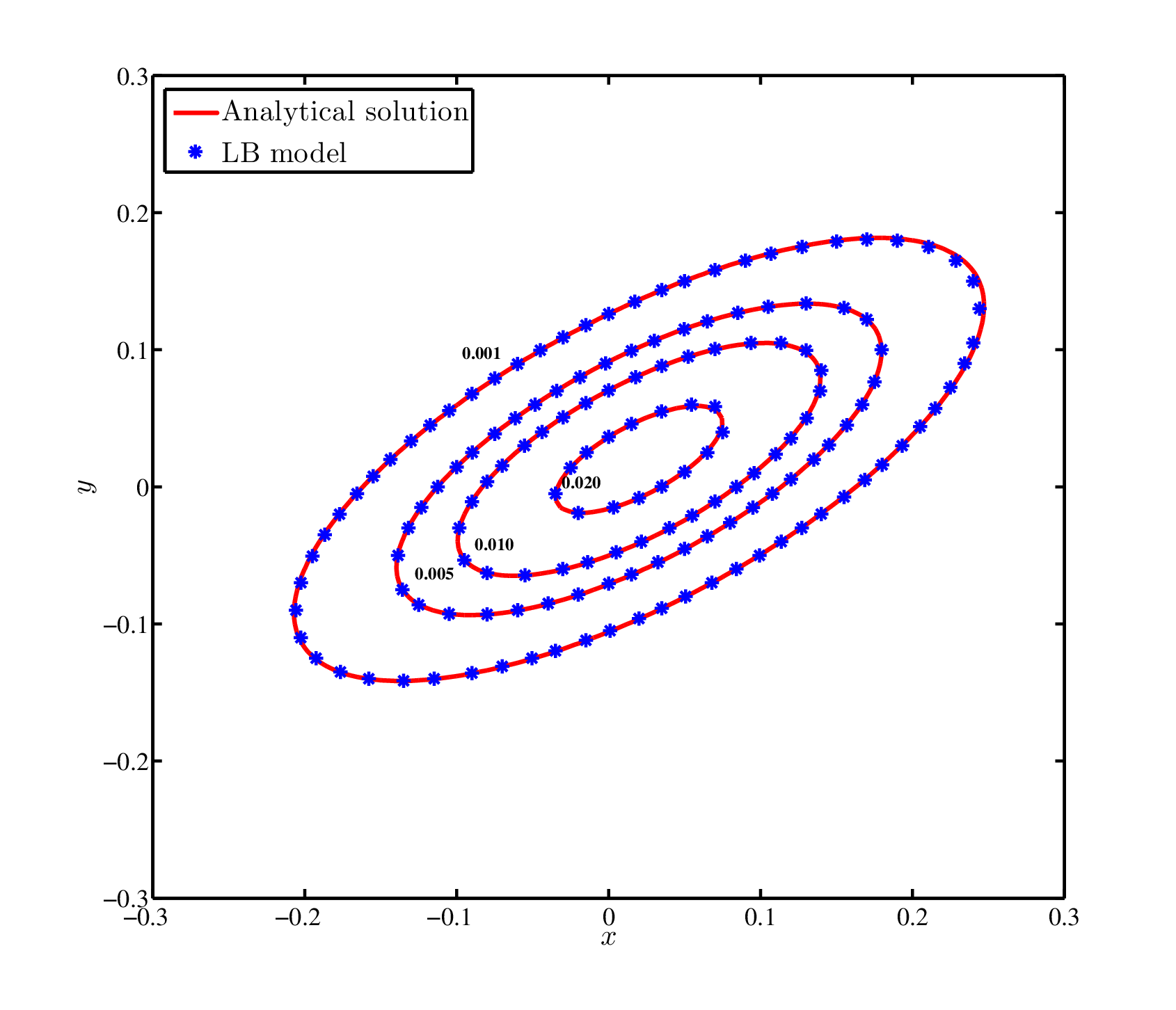}
		}
		\subfloat[GPMFD scheme]
		{
			\label{fig:subfig2}\includegraphics[width=0.4\textwidth]{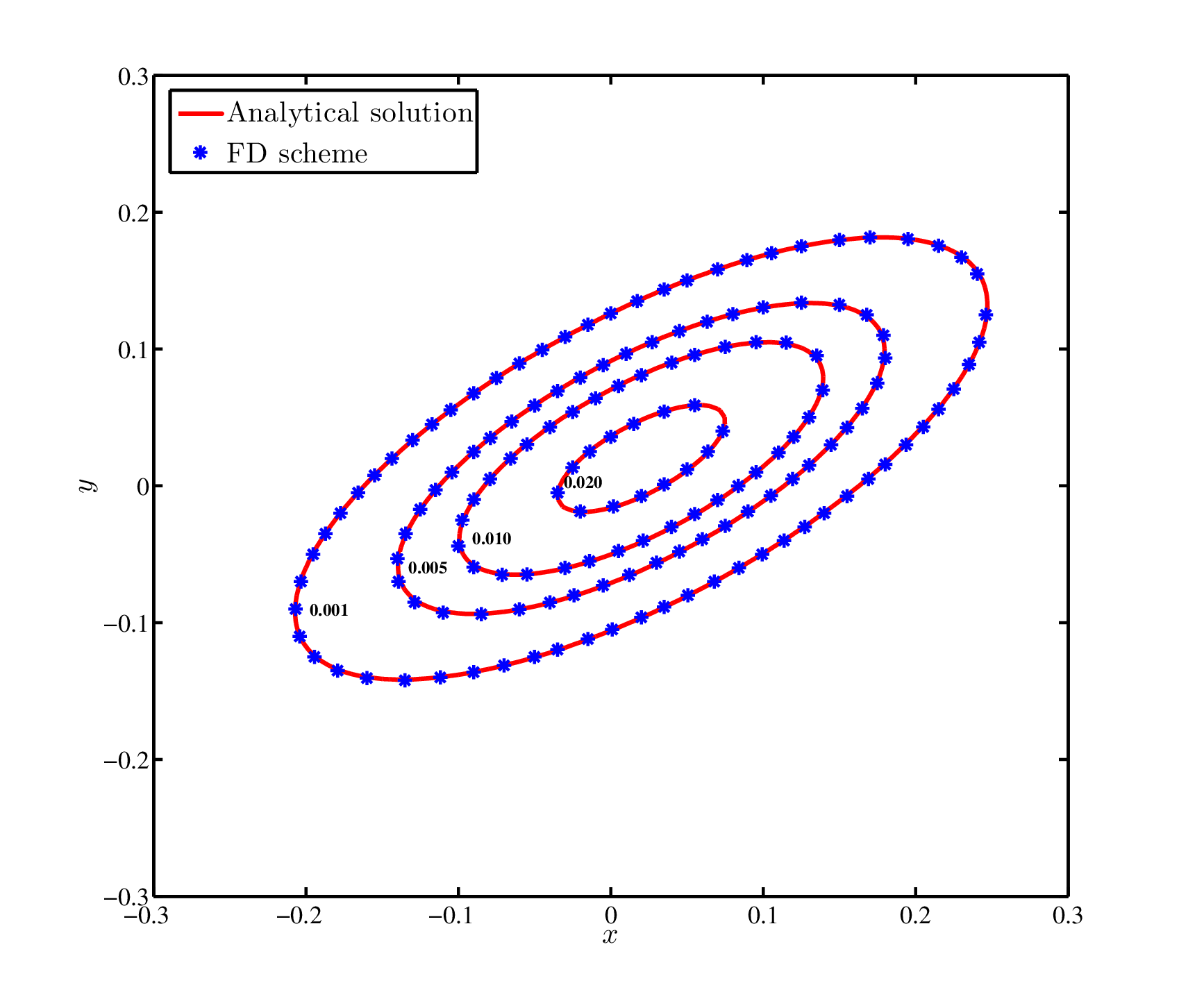}
		}
		\caption{Contour lines of the scalar variable $\phi$  at the time $t=2$ and $\bm{u}=(0.01,0.01)^T$  [(a-b) isotropic diffusion problem, (c-d) diagonal anisotropic diffusion problem, and (e-f) full anisotropic diffusion problem]. }     
		\label{Ex1-contour-line}            
	\end{figure}

We first conduct some tests with $u_x=u_y=0.01$, $\Delta x=\Delta t=1/200$ and following three types of diffusion coefficient matrices,
\begin{align}
	&	\bm{\kappa}= \left[\begin{matrix}
		1&1\\
		1&2\\
	\end{matrix}\right]\times 10^{-3},   	\bm{\kappa}=  \left[\begin{matrix}
		1&0\\
		0&2\\
	\end{matrix}\right]\times 10^{-3},  		\bm{\kappa}=  \left[\begin{matrix}
		1&0\\
		0&1\\
	\end{matrix}\right]\times 10^{-3}, 
\end{align}
which represent the isotropic diffusion, diagonal anisotropic diffusion and full anisotropic diffusion problems,
and present the results of the LW scheme ($b=a^2$)
at the time $t=2$ in Fig. \ref{Ex1-contour-line} where $a=0.5$. As shown in this figure,  the numerical results obtained from both the GPMRT-LB model and GPMFD scheme are in good agreement with the analytical solutions.

	In order to test the CRs of the GPMRT-LB model and GPMFD scheme at the acoustic scaling, we conduct some simulations of the following three different cases: 
	\begin{subequations}\label{three-case-u-kappa}
	\begin{align}
		&{\rm Case}\:1:\: u_x=u_y=0.1,			\bm{\kappa}=\alpha\times \left[\begin{matrix}
			1&1\\
			1&2\\
		\end{matrix}\right]\times 10^{-4},\alpha=0.625, \\
		&{\rm Case}\:2:\: u_x=u_y=0.1,			\bm{\kappa}= \alpha\times\left[\begin{matrix}
			1&1\\
			1&2\\
		\end{matrix}\right]\times 10^{-4},\alpha=1.000, \\
		&{\rm Case}\:3:\: u_x=u_y=1.0,			\bm{\kappa}=\alpha\times \left[\begin{matrix}
			1&1\\
			1&2\\
		\end{matrix}\right]\times 10^{-4},\alpha=0.625. 
	\end{align} 
	\end{subequations}
	For the diffusion coefficient matrices in Eq. (\ref{three-case-u-kappa}), we take the relaxation parameters $S^1_{xx}=0.8$, $S^1_{xy}=S^1_{yx}=0.4$, and $S^1_{yy}=1.2$ [see Eqs. (\ref{natural-condition}) and (\ref{kappa})], and measure the RMSEs and CRs under different values of the lattice spacing $\big(\Delta x=1/200,1/250,1/300,1/350\big)$ and time step $\big(\Delta t=1/200,1/250,1/300,1/350\big)$  in Table \ref{Tab-acou-NCDE-slope} where $t = 1$.
	\begin{sidewaystable}[htbp]
		\centering
			\caption{The RMSEs and CRs of GPMRT-LB model and GPMFD scheme at the acoustic scaling for different values of diffusion coefficient matrix $\bm{\kappa}$ and velocity $\bm{u}$.}	\label{Tab-acou-NCDE-slope}
			\begin{tabular*}{\textheight}{@{\extracolsep\fill}lccccccccccc}\toprule%
				\multirow{2}{*}{$\alpha$}&\multirow{2}{*}{$u_{x}$}&\multirow{2}{*}{$\Delta x$}&	\multirow{2}{*}{$a$}&\multirow{2}{*}{$b$}&\multicolumn{2}{c}{GPMRT-LB model}&&\multicolumn{2}{c}{GPMFD scheme}\\
				\cmidrule{6-7}\cmidrule{9-10}
				&&&&&RMSE&CR&&RMSE&CR \\ \midrule
				\multirow{4}{*}{0.625}&\multirow{4}{*}{0.1}
				&	 $1/200$   &0.306186217847897&0.112500000000000& 1.7185$\times 10^{-4}$ & - && 1.6970$\times 10^{-4}$ &	 - \\ 
				\cmidrule{4-10}
				&	& $1/250$ 	&0.342326598440729&0.140625000000000& 1.1439$\times 10^{-4}$ & 1.8241&&1.1322$\times 10^{-4}$&1.8136\\
				\cmidrule{4-10}
				&	& $1/300$ &0.375000000000000&0.168750000000000& 8.0130$\times 10^{-5}$ & 1.9522&&7.9437$\times 10^{-5}$&1.9436\\
				\cmidrule{4-10}
				&	& $1/350$    &0.405046293650491&0.196875000000000& 5.8337$\times 10^{-5}$ & 2.0592 && 5.7898$\times 10^{-5}$ & 2.0518 \\  
				\cmidrule{1-10}
				\multirow{4}{*}{1.000}&\multirow{4}{*}{0.1}
				& $1/200$  &0.387298334620742&0.180000000000000& 1.4564$\times 10^{-4}$&-&&1.4285$\times 10^{-4}$&-\\
				\cmidrule{4-10}
				&	& $1/250$ 	&0.433012701892219&0.225000000000000&9.4597$\times 10^{-5}$&1.9334&&9.3633$\times 10^{-5}$&1.9244\\
				\cmidrule{4-10}
				&	& $1/300$ &0.474341649025257&0.270000000000000&6.4280$\times 10^{-5}$&2.1192&&6.3717$\times 10^{-5}$&2.1112\\ 
				\cmidrule{4-10}
				&	& $1/350$   &0.512347538297980&0.315000000000000&  4.5156$\times 10^{-5}$ &  2.2908 && 4.4808$\times 10^{-5}$ & 2.2839 \\ 
				\cmidrule{1-10}
				\multirow{4}{*}{0.625}&\multirow{4}{*}{1.0	}
				& $1/200$    & 0.306186217847897&0.112500000000000& 1.7213$\times 10^{-4}$ & - && 1.6998$\times 10^{-4}$ &	- \\ 
				\cmidrule{4-10}
				&	& $1/250$ &0.342326598440729&0.140625000000000& 1.1456$\times 10^{-4}$ & 1.8247  && 1.1393$\times 10^{-4}$ & 1.8142 \\ 
				\cmidrule{4-10}
				&	& $1/300$  &0.375000000000000&0.168750000000000 & 8.0248$\times 10^{-5}$ & 1.9525 && 7.9554$\times 10^{-5}$ & 1.9439 \\ 
				\cmidrule{4-10}
				&	& $1/350$   &0.405046293650491&0.196875000000000& 5.8422$\times 10^{-5}$ & 2.0593 && 5.7982$\times 10^{-5}$ & 2.0519 \\ 
				\botrule
							\end{tabular*}
							\centering
	\end{sidewaystable}

	In addition, we also conduct some simulations at the diffusive scaling under  $u_x=u_y=0.01$ and the full anisotropic diffusion coefficient matrix,
	\begin{align}
		\bm{\kappa}=\left[\begin{matrix}
			1&1\\
			1&2\\
		\end{matrix}\right]\times 10^{-3},
	\end{align}
	and consider  the following five cases with different values of parameters $a$ and $b$: 
	\begin{subequations}\label{ncde-five-cses}
		\begin{align}
		&{\rm Case}\:1:\: a=b=1, {\rm the\:MRT-LB\:model},\\
		&{\rm Case}\:2:\: a^2<b<a, a=0.4,b=0.2 ,\\
		&{\rm Case}\:3:\: a<b<1, a=0.4,b=0.45,\\
		&{\rm Case}\:4:\: b=a^2, a=0.6,{\rm the\:LW\:scheme},\\
		&{\rm Case}\:5:\: b=a=0.5, {\rm the\:FP\:scheme}. 
	\end{align} 
	\end{subequations}

		\begin{table} \caption{The RMSEs and CRs of  GPMRT-LB model at the diffusive scaling for five cases of parameters $a$ and $b$ (\ref{ncde-five-cses}) ($t=2$). }\label{Tab-diff-NCDE-LB-slope}  
		\begin{tabular*}{\textwidth}{@{\extracolsep\fill}lccccccccc}\toprule%
			$\Delta x$&
			$\Delta t$ &$(a,b)$&RMSE$_{\Delta x}$&RMSE$_{\Delta x/2}$&RMSE$_{\Delta x/4}$&RMSE$_{\Delta x/8}$&CR\\
			\midrule
			$\frac{1}{80}$&$\frac{1}{50}$&(1,1)&9.5641$\times 10^{-6}$&2.3468$\times 10^{-6}$&1.0396$\times 10^{-6}$&5.8420$\times 10^{-7}$&$\sim$2.0108\\
				$\frac{1}{80}$&$\frac{1}{50}$&(0.4,0.2)&2.8713$\times 10^{-5}$&6.7951$\times 10^{-6}$&2.9940$\times 10^{-6}$&1.6794$\times 10^{-6}$&$\sim$2.0368\\
				$\frac{1}{80}$&$\frac{1}{50}$&(0.4,0.45)&2.5173$\times 10^{-5}$&6.0406$\times 10^{-6}$&2.6671$\times 10^{-6}$&1.4971$\times 10^{-6}$&$\sim$2.0274\\
				$\frac{1}{80}$&$\frac{1}{50}$&(0.6,0.36)&1.4732$\times 10^{-5}$&3.5532$\times 10^{-6}$&1.5702$\times 10^{-6}$&8.8162$\times 10^{-7}$&$\sim$2.0240\\
				$\frac{1}{80}$&$\frac{1}{50}$&(0.5,0.5)&1.7176$\times 10^{-5}$&4.1824$\times 10^{-6}$&1.8512$\times 10^{-6}$&1.0400$\times 10^{-6}$&$\sim$2.0175\\
			\botrule
		\end{tabular*}
	\end{table}

	\begin{table} 
	\caption{The RMSEs and CRs of GPMFD scheme at the diffusive scaling  for five cases of parameters $a$ and $b$ (\ref{ncde-five-cses}) ($t=2$).}
	\label{Tab-diff-NCDE-FD-slope}  
		\begin{tabular*}{\textwidth}{@{\extracolsep\fill}lccccccccc}\toprule%
			$\Delta x$&
			$\Delta t$ &$(a,b)$&RMSE$_{\Delta x}$&RMSE$_{\Delta x/2}$&RMSE$_{\Delta x/4}$&RMSE$_{\Delta x/8}$&CR\\
			\midrule
			$\frac{1}{80}$&$\frac{1}{50}$&(1,1)&6.1121$\times 10^{-6}$&1.5280$\times 10^{-6}$&6.7934$\times 10^{-7}$& 3.8222$\times 10^{-7}$&$\sim$1.9993\\
				$\frac{1}{80}$&$\frac{1}{50}$&(0.4,0.2)&2.7548$\times10^{-5}$&6.8421$\times 10^{-6}$&3.0392$\times 10^{-6}$&1.7095$\times 10^{-6}$&$\sim$2.0036\\
					$\frac{1}{80}$&$\frac{1}{50}$&(0.4,0.45)&1.0589$\times 10^{-5}$&2.6788$\times 10^{-6}$&1.1932$\times 10^{-6}$&6.7179$\times 10^{-7}$&$\sim$1.9041 \\
					$\frac{1}{80}$&$\frac{1}{50}$&(0.6,0.36)&7.6951$\times 10^{-6}$&1.9097$\times 10^{-6}$&8.4799$\times 10^{-7}$&4.7691$\times 10^{-7}$&$\sim$2.0045 \\
					$\frac{1}{80}$&$\frac{1}{50}$&(0.5,0.5)&9.2178$\times 10^{-6}$&2.3532$\times 10^{-6}$&1.0501$\times 10^{-6}$&5.9156$\times 10^{-7}$&$\sim$1.9849 \\
			\botrule
		\end{tabular*}
	\end{table}

As seen from the Tables \ref{Tab-acou-NCDE-slope}, \ref{Tab-diff-NCDE-LB-slope} and \ref{Tab-diff-NCDE-FD-slope}, both the GPMRT-LB model and GPMFD scheme can achieve a second-order CR at both acoustic and diffusive scalings, which is consistent with the theoretical analysis.\\
\textbf{Example 2}  For the second problem,  the Poiseuille flow between two parallel no-slip walls driven by a constant body force $\bm{\hat{F}}=(G,0)$, and   the analytical solution can be expressed as \begin{align}
	u_x(y)=4U\Big(1-\frac{y}{H}\Big)\frac{y}{H},u_y=0,
\end{align}where $y\in[0, H]$,  $H$ is the channel width and $U=GH^2/8\upsilon$ is the maximal velocity along the centerline of the channel.

In our simulations, $G=0.08\upsilon$ and $H=1$, the periodic boundary condition is used in the horizontal direction. In the vertical direction, however, we adopt the following scheme for the GPMRT-LB model \cite{Zhao2020}: 
\begin{align}
	&f_i(\bm{x},t_{n+1})=p_{-1}f_{\overline{i}}^{\star}(\bm{x},t_n)+p_0f_i^{\star}(\bm{x},t_n)+p_1f_i^{\star}(\bm{x}+\bm{\lambda}_i\Delta t,t_n),\\
	&f_{\overline{i}}(\bm{x},t_{n+1})=p_{-1}f_{\overline{i}}^{\star}(\bm{x}+\bm{\lambda}_i\Delta t,t_n)+p_0f_{\overline{i}}^{\star}(\bm{x},t_n)+p_1f_i^{\star}(\bm{x},t_n),
\end{align} 
where $\bm{x}-\bm{\lambda}_i\Delta t\notin\Omega$ and $\bm{x}+\bm{\lambda}_i\Delta t\in\Omega$  are the boundary points at the top and bottem straight boundaries, $\Omega$ represents the computational domain,  $\overline{i}$ is defined as $\bm{\lambda}_{\overline{i}}=-\bm{\lambda}_i$, and each boundary point is located at the middle of two neighboring grid points, ensuring that the scheme has a second-order accuracy \cite{Kruger2017}. Here the convergence criterion  for this steady state problem is defined as follows:
\begin{align}
	\sqrt{\frac{\sum_i[\bm{u}(\bm{x},t_n)-\bm{u}(\bm{x},t_n-100\Delta t)]^2}{\sum_i[\bm{u}(\bm{x},t_n)]^2}}<10^{-10}.
\end{align}

\begin{figure}[htbp]    
		\centering            
		\subfloat[GPMRT-LB model]    
		{
			\label{fig:subfig1}\includegraphics[width=0.4\textwidth]{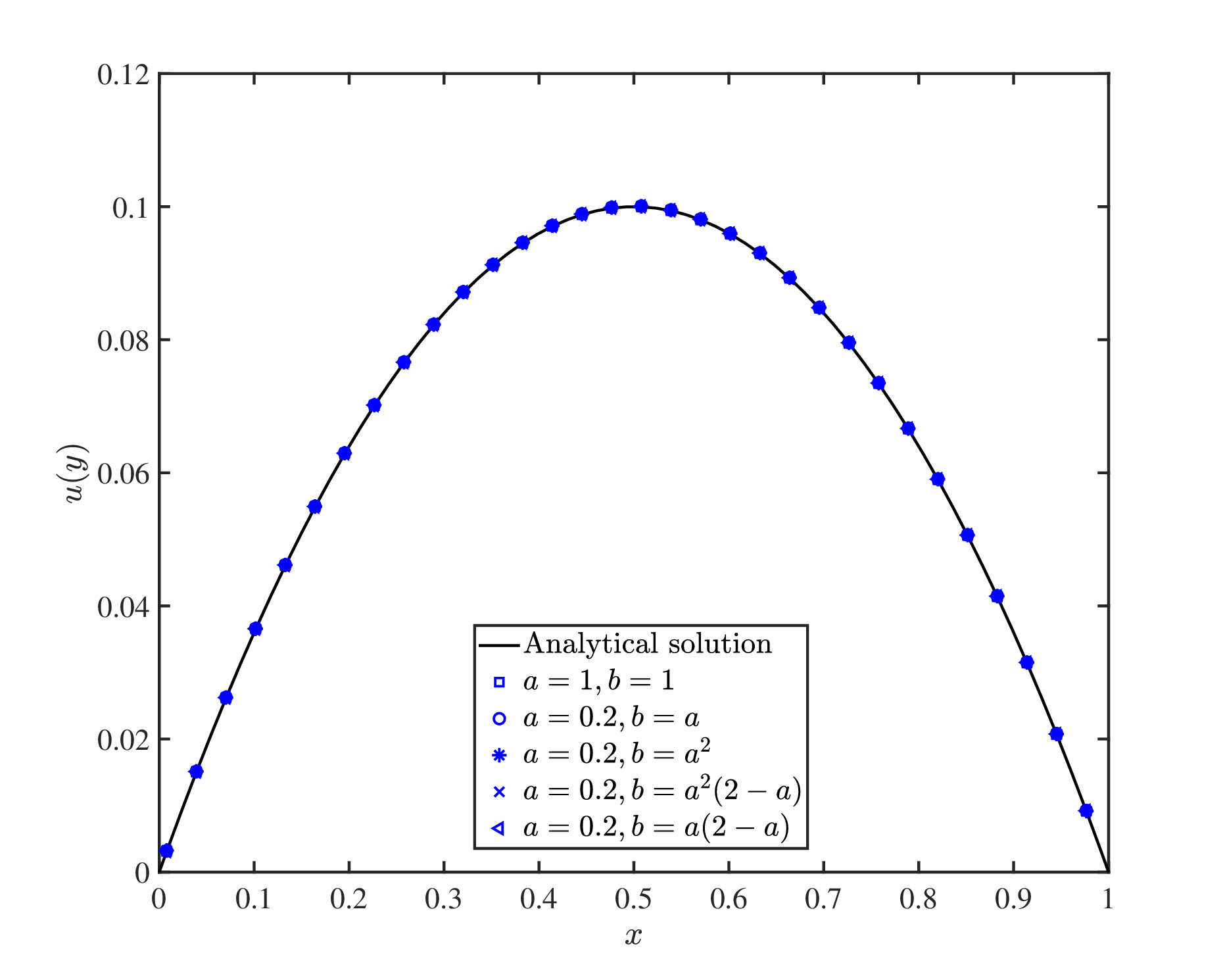} 
		}
		\subfloat[GPMFD scheme]
		{
			\label{fig:subfig2}\includegraphics[width=0.4\textwidth]{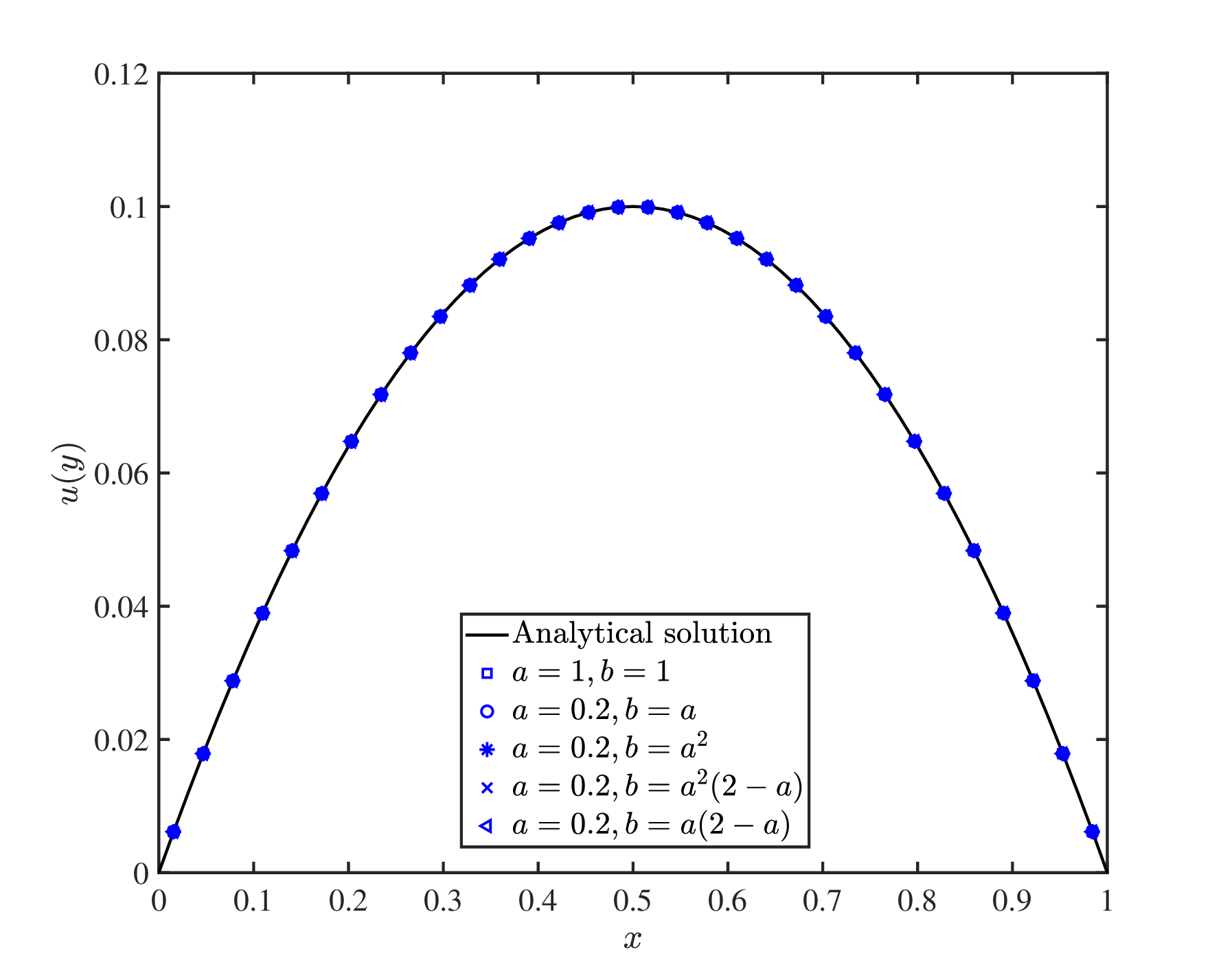}
		} 
		\quad\subfloat[GPMRT-LB model]    
		{
			\label{fig:subfig1}\includegraphics[width=0.4\textwidth]{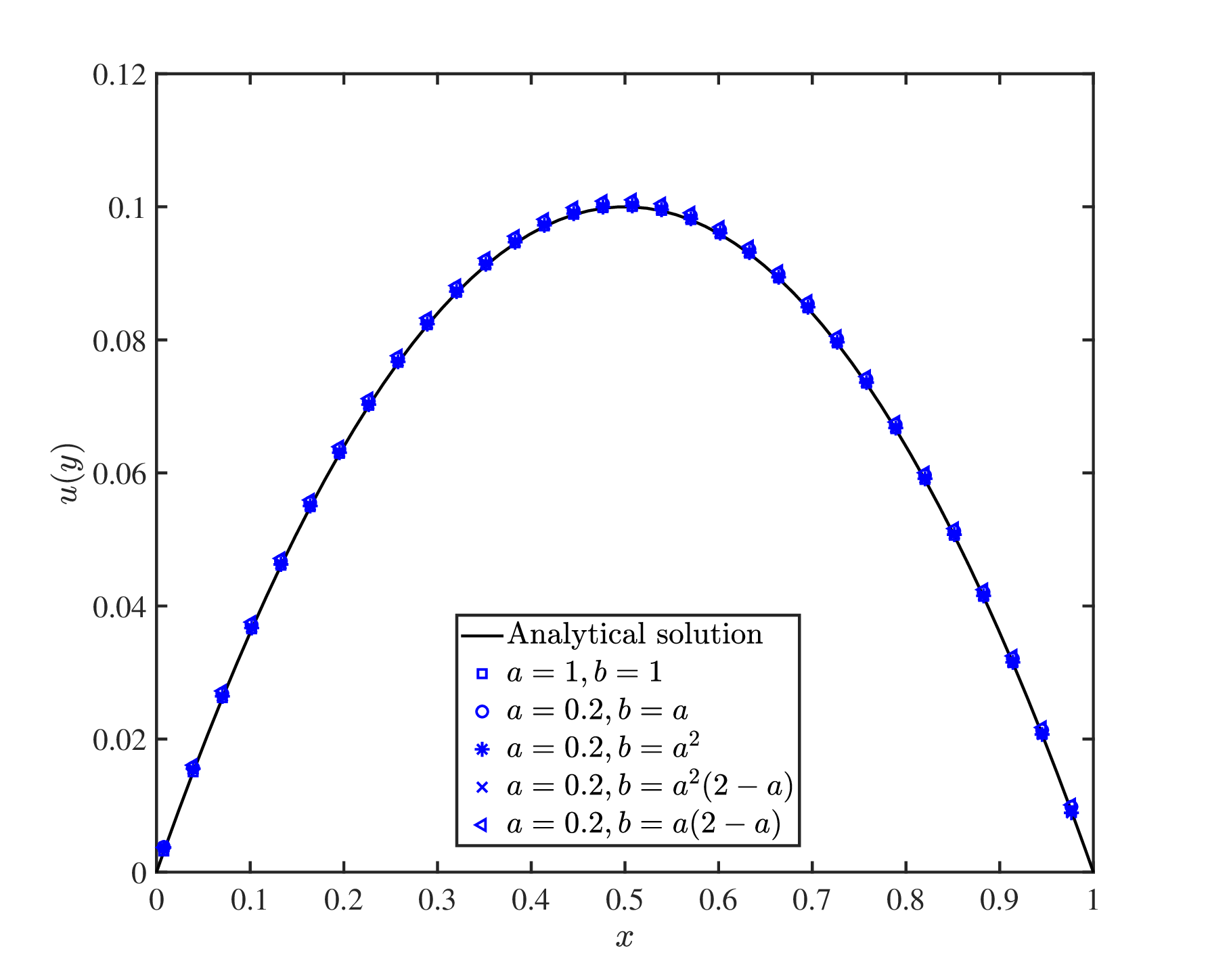} 
		}
		\subfloat[GPMFD scheme]
		{
			\label{fig:subfig2}\includegraphics[width=0.4\textwidth]{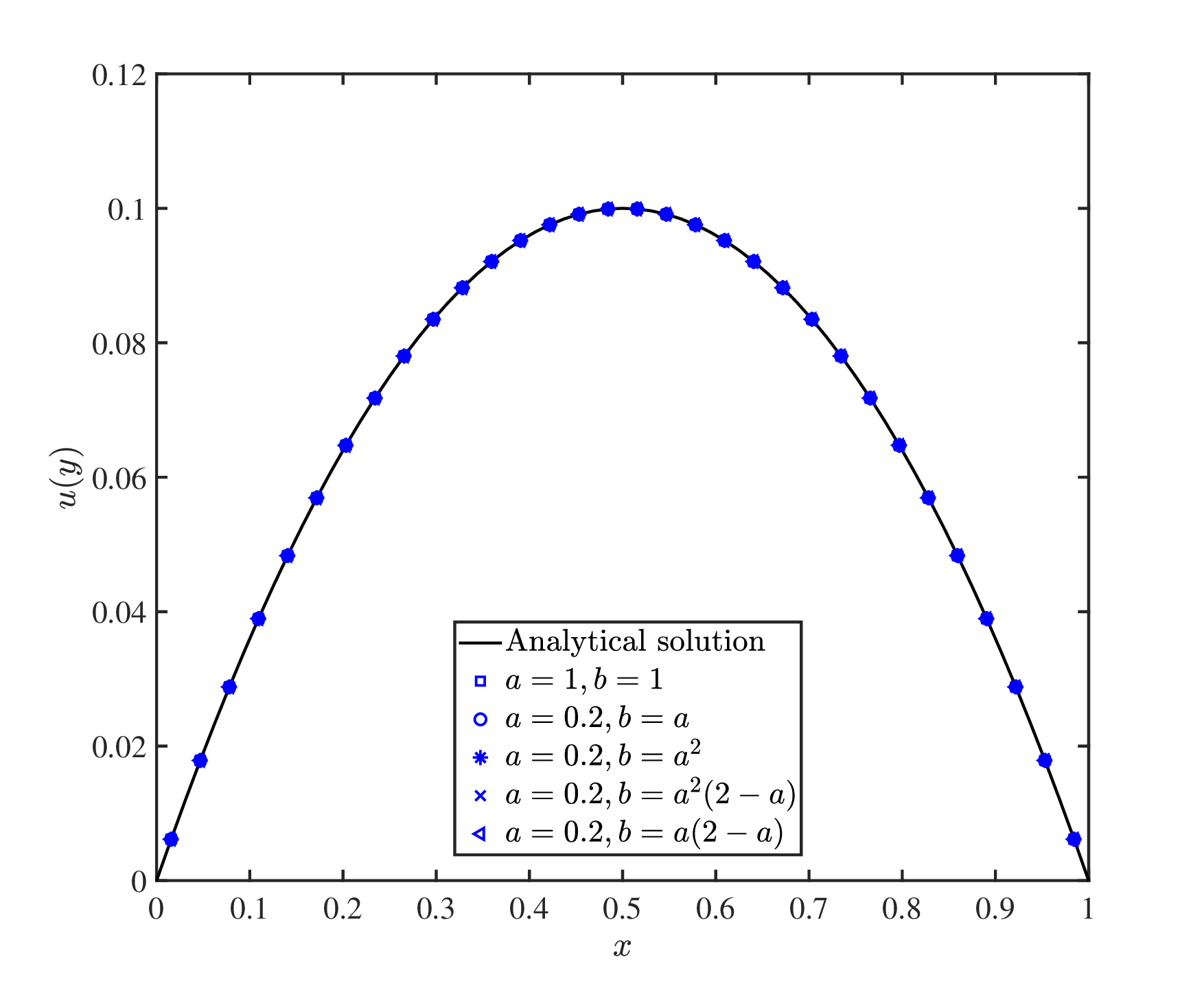}
		} 
		\caption{The numerical and analytical solutions of the velocity $u_x$ along the vertical direction   under different values of  $a$ and $b$ [(a-b) $\upsilon=0.02$ and (c-d) $\upsilon=0.10$]. }    
		\label{Ex2-different-nu}            
	\end{figure}

We first carry out some simulations of following five cases of parameters $a$ and $b$ \cite{Guo2001}: 
\begin{subequations}\label{poise-five-cses}
	\begin{align}
		&{\rm Case}\:1:\: a=1,b=1, {\rm the\: MRT-LB\: model},\\
		&{\rm Case}\:2:\: a=0.2,b=a^2,{\rm the\:LW\:scheme },\\
		&{\rm Case}\:3:\: a=0.2,b=a,{\rm  the\:FP\:scheme },\\
		&{\rm Case}\:4:\: a=0.2,b=a^2(2-a),(a^2<b<a),\\
		&{\rm Case}\:5:\: a=0.2,b=a(2-a),(a<b<1), 
	\end{align} 
\end{subequations}
and present the results in Fig. \ref{Ex2-different-nu}  where $\Delta x=1/100$, $\Delta t=1/200$, $\upsilon=0.02$ and $\upsilon=0.10$. From this figure, one can observe that  the numerical solutions of the GPMRT-LB model and GPMFD scheme are in good agreement with the corresponding analytical solutions.

	Next, we take the viscosity coefficient $\upsilon=0.06$ to test the CR at the diffusive scaling, and show the   results  in Tables \ref{Tab-diff-NSE-LB-slope} and \ref{Tab-diff-NSE-FD-slope}, where a second-order accuracy is clearly observed.

	\begin{table} 
			\caption{The RMSEs and CRs of  GPMRT-LB model at the diffusive scaling for five cases of parameters $a$ and $b$ (\ref{poise-five-cses}). }
			\label{Tab-diff-NSE-LB-slope}  
		\begin{tabular*}{\textwidth}{@{\extracolsep\fill}lccccccccc}\toprule%
				$\Delta x$&
				$\Delta t$ &$(a,b)$&RMSE$_{\Delta x}$&RMSE$_{\Delta x/2}$&RMSE$_{\Delta x/4}$&RMSE$_{\Delta x/8}$&CR\\
				\midrule
				$\frac{1}{16}$&$\frac{1}{50}$&(1,1)&5.6927$\times 10^{-4}$&1.4234$\times 10^{-4}$&3.5586$\times 10^{-5}$&8.8965$\times 10^{-6}$&$\sim$2.0000 \\
				$\frac{1}{16}$&$\frac{1}{50}$&(0.2,0.04)&6.3294$\times 10^{-3}$&1.5823$\times 10^{-3}$&3.9559$\times 10^{-4}$&9.8896$\times 10^{-5}$&$\sim$2.0000 \\
				$\frac{1}{16}$&$\frac{1}{50}$&(0.2,0.2)&3.8197$\times 10^{-3}$&7.6164$\times 10^{-4}$& 1.6075$\times 10^{-4}$& 3.5918$\times 10^{-5}$&$\sim$2.2442 \\
				$\frac{1}{16}$&$\frac{1}{50}$&(0.2,0.072)&3.8478$\times 10^{-3}$&8.6419$\times 10^{-4}$&2.0272$\times 10^{-4}$&4.8929$\times 10^{-5}$&$\sim$2.0990\\
				$\frac{1}{16}$&$\frac{1}{50}$&(0.2,0.36)&8.4321$\times 10^{-3}$&2.1739$\times 10^{-3}$&5.5153$\times 10^{-4}$&1.3888$\times 10^{-4}$&$\sim$1.9747\\
				\botrule
			\end{tabular*}
	\end{table}
	
	\begin{table} 
			\caption{The RMSEs and CRs of  GPMFD model  at the diffusive scaling for five cases of parameters $a$ and $b$ (\ref{poise-five-cses}). }
			\label{Tab-diff-NSE-FD-slope}  
			\begin{tabular*}{\textwidth}{@{\extracolsep\fill}lccccccccc}\toprule%
				$\Delta x$&
				$\Delta t$ &$(a,b)$&RMSE$_{\Delta x}$&RMSE$_{\Delta x/2}$&RMSE$_{\Delta x/4}$&RMSE$_{\Delta x/8}$&CR\\
				\midrule
				$\frac{1}{16}$&$\frac{1}{50}$&(1,1)& 1.0618$\times 10^{-4}$ &2.6223$\times 10^{-5}$&5.7352$\times 10^{-6}$&1.4285$\times 10^{-6}$&$\sim$2.0719\\
				$\frac{1}{16}$&$\frac{1}{50}$&(0.2,0.04)& 1.2866$\times 10^{-4}$&3.1207$\times 10^{-5}$&7.1306$\times 10^{-6}$&1.7756$\times 10^{-6}$&$\sim$2.0823\\
				$\frac{1}{16}$&$\frac{1}{50}$&(0.2,0.2)& 1.4639$\times 10^{-4}$&3.1882$\times 10^{-5}$&6.8274$\times 10^{-6}$&1.5054$\times 10^{-6}$&$\sim$2.2011\\
				$\frac{1}{16}$&$\frac{1}{50}$&(0.2,0.072)&1.8749$\times 10^{-4}$&4.4785$\times 10^{-5}$&1.0404$\times 10^{-5}$&2.5998$\times 10^{-6}$&$\sim$2.1091 \\
				$\frac{1}{16}$&$\frac{1}{50}$&(0.2,0.36)&7.3677$\times 10^{-4}$&1.6241$\times 10^{-4}$&3.9908$\times 10^{-5}$&9.8000$\times 10^{-6}$&$\sim$2.0774 \\
				\botrule
			\end{tabular*}
	\end{table}

\noindent  \textbf{Example 3} We further consider the CDE (\ref{CDE}) problem with the  periodic boundary condition and  following initial condition: 
\begin{align}
	&\phi(x,0)=\sin(\pi x),-1\leq x\leq1, 
\end{align} and obtain the analytical solution of this problem as
\begin{align}
	\phi(x,t)=\sin[\pi(x-ut)]\exp{(-\kappa\pi^2t)}.
\end{align}
In this test, we consider the diffusion coefficient $\kappa=0.08$, velocity $u=1$,  lattice spacing $\Delta x=1/10,1/20,1/40,1/80$, time step $\Delta t=1/50,1/200,1/800,1/3200$, and measure the RMSEs between the numerical and analytical solutions at the time $t=2$. As seen from  Tables \ref{Tab-diff-CDE-LB-slope} and \ref{Tab-diff-CDE-FD-slope},  both the GPMRT-LB model and GPMFD scheme  at the diffusive scaling have a fourth-order CR in space.

\begin{table} 
	\caption{The RMSEs and CRs of F-GPMRT-LB model at the diffusive scaling for different value of $a$ and $b$.  }
	\label{Tab-diff-CDE-LB-slope}  
	\begin{tabular*}{\textwidth}{@{\extracolsep\fill}lccccccccc}\toprule%
		$\Delta x$&
		$\Delta t$ &$(a,b)$&RMSE$_{\Delta x}$&RMSE$_{\Delta x/2}$&RMSE$_{\Delta x/4}$&RMSE$_{\Delta x/8}$&CR\\
		\midrule
		$\frac{1}{10}$&$\frac{1}{50}$&(1,1)& 6.1216$\times 10^{-4}$ &3.7760$\times 10^{-5}$&2.3466$\times 10^{-6}$&1.4628$\times 10^{-7}$&$\sim$4.0103\\
		$\frac{1}{10}$&$\frac{1}{50}$&(0.6,0.9)& 1.6485$\times 10^{-3}$&1.0545$\times 10^{-4}$&6.6188$\times 10^{-6}$&4.1442$\times 10^{-7}$&$\sim$3.9859\\
		$\frac{1}{10}$&$\frac{1}{50}$&(0.9,0.8)& 5.8918$\times 10^{-4}$&3.6349$\times 10^{-5}$&1.1590$\times 10^{-6}$&1.4082$\times 10^{-7}$&$\sim$4.0110\\
		$\frac{1}{10}$&$\frac{1}{50}$&(0.6,0.6)&5.4684$\times 10^{-4}$&3.3716$\times 10^{-5}$&2.0952$\times 10^{-6}$&1.3062$\times 10^{-7}$&$\sim$4.0105 \\ 
		\botrule
	\end{tabular*}
\end{table}
\begin{table}
	\caption{The RMSEs and CRs of F-GPMFD model at the diffusive scaling for different values of $a$ and $b$.  }
	\label{Tab-diff-CDE-FD-slope}  
	\begin{tabular*}{\textwidth}{@{\extracolsep\fill}lccccccccc}\toprule%
		$\Delta x$&
		$\Delta t$ &$(a,b)$&RMSE$_{\Delta x}$&RMSE$_{\Delta x/2}$&RMSE$_{\Delta x/4}$&RMSE$_{\Delta x/8}$&CR\\
		\midrule%
		$\frac{1}{10}$&$\frac{1}{50}$&(1,1)&5.2505$\times 10^{-4}$ &3.7716$\times 10^{-5}$&2.5180$\times 10^{-6}$&1.6268$\times 10^{-7}$&$\sim$3.8774\\
		$\frac{1}{10}$&$\frac{1}{50}$&(0.6,0.9)& 7.2309$\times 10^{-4}$&5.1108$\times 10^{-5}$&3.3636$\times 10^{-6}$&2.1540$\times 10^{-7}$&$\sim$3.9043\\
		$\frac{1}{10}$&$\frac{1}{50}$&(0.9,0.8)& 2.0699$\times 10^{-4}$&1.4997$\times 10^{-5}$&1.0076$\times 10^{-6}$&6.5284$\times 10^{-8}$&$\sim$3.8768\\
		$\frac{1}{10}$&$\frac{1}{50}$&(0.6,0.6)&1.8924$\times 10^{-4}$&1.3172$\times 10^{-5}$&9.2169$\times 10^{-7}$&5.9723$\times 10^{-8}$&$\sim$3.8766 \\ 
		\botrule%
	\end{tabular*}
\end{table}

	\section{Conclusions}\label{Conclusion}

In this paper, we first derived the multiple-level GPMFD scheme  of the GPMRT-LB model on conservative moments, and then conducted the accuracy analysis for the GPMRT-LB model and GPMFD scheme through the Maxwell iteration method at the diffusive and acoustic scalings. Furthermore, for the NACDE and NSEs, we presented the first- and second-order modified equations of the GPMRT-LB model and GPMFD scheme at both diffusive and acoustic scalings. In particular, based on our previous work \cite{Chen2023}, we also developed the F-GPMRT-LB model and F-GPMFD scheme at the diffusive scaling for the one-dimensional CDE, which can be more stable than the MRT-LB model and the corresponding macroscopic finite-difference scheme \cite{Chen2023} through adjusting parameters $a$ and $b$ properly [see Fig. (\ref{stability-region})]. Finally, some numerical simulations of Gauss hill problem and Poiseuille flow were conducted to test the GPMRT-LB model and GPMFD scheme, and the results show that both of them have a second-order convergence rates in space. We also performed  a numerical test on  the F-GPMRT-LB model and F-GPMFD scheme for the one-dimensional CDE, and it is found that they are of fourth-order accuracy in space, which is consistent with our theoretical analysis.
\section*{Acknowledgements}
The computation is completed in the HPC Platform of Huazhong University of Science and Technology. This work was financially supported by the National Natural Science Foundation of China (Grants No. 12072127 and No. 51836003) and Interdiciplinary Research Program of Hust (2023JCJY002).

\begin{appendices}

	\setcounter{equation}{0}
\renewcommand\theequation{A.\arabic{equation}}
\section{ Derivation of Eqs. (\ref{second-diffusive-NACDEs}), (\ref{diff-con-second})  and (\ref{mon-second-diff})}
For the first-order ME (\ref{modified-equation-first-order-diffusive}), one can obtain 
\begin{align}\label{zero-order-diffuisve-appendix}
	&\big[\partial_t\bm{m}^{eq}\big]_1=\partial_t\sum_if_i^{eq},\Big[c\bm{\mathcal{W}}_0\bm{m}^{eq}\Big]_1=c\partial_{\alpha}\sum_i\bm{e}_{i\alpha}f_i^{eq},\bm{\tilde{F}}_1=\sum_{i=1}^q\Big(F_i+G_i+\frac{\Delta t}{2}\overline{D}_iF_i\Big),\notag\\
	&\Bigg[\Delta tc^2\bm{\mathcal{W}}_0\bigg(\bm{\hat{S}}_N^{-1}+\Big(\frac{b}{2a^2}-1\Big)\bm{I}\bigg)\bm{\mathcal{W}}_0\bm{m}^{eq}\Bigg]_{1} =\Delta tc^2 \partial_{\beta}\sum_{k=1}^{q}\sum_{j=1}^{q}\bm{e}_{j\beta}
	\bigg[\bm{\overline{\Lambda}}_{jk}+\Big(\frac{b}{2a^2}-1\Big)\delta_{jk}\bigg]\notag\\
	&\times \partial_{\theta}\bm{e}_{k\theta}f_k^{eq}=\Delta tc^2\partial_{\beta}\sum_{k=1}^{q}\big(S^{10}_{\beta}+S^{1}_{\beta\gamma}\bm{e}_{k\gamma}\big)\partial_{\theta}\bm{e}_{k\theta}f_k^{eq}+\Delta t\Big(\frac{b}{2a^2}-\frac{1}{2}\Big)\partial_{\beta}\partial_{\theta}\sum_{k=1}^q\bm{c}_{k\beta}\bm{c}_{k\theta}f_k^{eq}. 
\end{align} 
For the second-order ME (\ref{modified-equation-second-order-diffusive}), however, apart from   above equalities, we also need the following results:
\begin{align}\label{first-order-diffuisve-appendix}
	&\big[c \bm{\mathcal{W}}_0\bm{\hat{S}}_N^{-1}\bm{\tilde{F}}\big]_1=c\partial_{\alpha}\sum_{k=1}^{q }\bm{e}_{k\alpha}(F_k+G_k+\frac{\Delta t}{2}\overline{D}_kF_k)\notag\\
	&=c\partial_{\alpha}\sum_{k=1}^{q }\sum_{j=1}^{q }\bm{e}_{i\alpha}\bm{\overline{\Lambda}}_{jk}(F_k+G_k+\frac{\Delta t}{2} \overline{D}_kF_k)\notag\\
	&\qquad\qquad\qquad=\partial_{\alpha}\sum_{k=1}^{q }(cS^{10}_{\alpha}+S^1_{\alpha\beta}\bm{c}_{k\beta})(F_k+G_k+\frac{\Delta t}{2}\overline{D}_kF_k),\notag\\
	&\big[c\bm{\mathcal{W}}_0\bm{\hat{S}}_N^{-1}\partial_t\bm{m}^{eq}\big]_1=c\partial_{\alpha}\sum_{k=1}^{q }\sum_{j=1}^{q}\bm{e}_{i\alpha}\bm{\overline{\Lambda}}_{jk}\partial_tf_k^{eq}=\partial_{\alpha}\sum_{k=1}^{q}(cS^{10}_{\alpha}+S^1_{\alpha\beta}\bm{c}_{k\beta})\partial_tf^{eq}_k ,\notag\\
	&\big[c\partial_t\bm{\hat{S}}_N^{-1} \bm{\mathcal{W}}_0\bm{m}^{eq}\big]_1=s^{0}\partial_t\partial_{\alpha}\sum_{k=1}^{q}\bm{c}_{k\alpha}f^{eq}_k  ,   \big(c \bm{\mathcal{W}}_0\partial_t\bm{m}^{eq}\big)_1=\partial_{\alpha}\partial_t\sum_{k=1}^q\bm{c}_{k\alpha}f_k^{eq},
\end{align}
and
\begin{align}\label{second-order-diffuisve-appendix}
	& \bigg[ \Delta t^2 c ^3\bm{\mathcal{W}}_0^3\bm{m}^{eq} \bigg]_1 =\Delta t^2 \partial_{\alpha}\partial_{\beta}\partial_{\gamma}\sum_i\bm{c}_{i\alpha}\bm{c}_{i\beta}\bm{c}_{i\gamma}f_i^{eq} ,\notag\\
	& \bigg[ \Delta t^2 c ^3\bm{\mathcal{W}}_0^2\bm{\hat{S}}_N^{-1}\bm{\mathcal{W}}_0\bm{m}^{eq} \bigg]_1 =c^3\Delta t^2\partial_{\alpha}\partial_{\beta}\sum_{k=1}^{q}\bigg( \sum_{j=1}^{q }\bm{e}_{j\alpha}\bm{e}_{j\beta}\bm{\overline{\Lambda}}_{jk} \bigg)\partial_{\gamma}\bm{e}_{k\gamma}f_k^{eq} \notag\\
	&\qquad=\Delta t^2\partial_{\alpha}\partial_{\beta}\sum_{k=1}^{q }\bigg( c^2S^{20}_{\alpha\beta} +S^{2}_{\alpha\beta\theta\zeta}\bm{c}_{k\theta}\bm{c}_{k\zeta} \bigg)\partial_{\gamma}\bm{c}_{k\gamma}f_k^{eq}, \notag\\
	& \bigg[ \Delta t^2 c ^3\bm{\mathcal{W}}_0\bm{\hat{S}}_N^{-1}\bm{\mathcal{W}}_0^2\bm{m}^{eq} \bigg]_1 =c^3\Delta t^2\partial_{\alpha} \sum_{k=1}^{q }\bigg(\sum_{j=1}^{q}\bm{e}_{j\alpha}\bm{\overline{\Lambda}}_{jk}\bigg)\partial_{\beta}\bm{e}_{k\beta}\partial_{\gamma}\bm{e}_{k\gamma}f_k^{eq} \notag\\
	&\qquad=\Delta t^2\partial_{\alpha}\sum_{k=1}^{q }\bigg(cS^{10}_{\alpha}+S^{1}_{\alpha\theta}\bm{c}_{k\theta}\bigg)\partial_{\beta}\bm{c}_{k\beta}\partial_{\gamma}\bm{c}_{k\gamma}f_k^{eq},\notag\\
	& \bigg[ \Delta t^2 c ^3\bm{\mathcal{W}}_0\bm{\hat{S}}_N^{-1}\bm{\mathcal{W}}_0\bm{\hat{S}}_N^{-1}\bm{\mathcal{W}}_0\bm{m}^{eq} \bigg]_1 \notag\\
	&=c^3\Delta t^2\partial_{\alpha}\sum_{i=1}^{q }\Bigg[\sum_{k=1}^{q }\bigg(\sum_{j=1}^{q }\bm{e}_{j\alpha}\bm{\overline{\Lambda}}_{jk}\bigg)\partial_{\beta}\bm{e}_{k\beta}\bm{\overline{\Lambda}}_{ki}\Bigg]\partial_{\gamma}\bm{e}_{i\gamma}f_i^{eq} \notag\\
	&\qquad=\Delta t^2\partial_{\alpha}\sum_{i=1}^{q}\bigg[\sum_{k=1}^{q }\bigg(cS^{10}_{\alpha}+S^{1}_{\alpha\theta}\bm{c}_{k\theta}\bigg)\partial_{\beta}\bm{c}_{k\beta}\bm{\overline{\Lambda}}_{ki}\bigg]\partial_{\gamma}\bm{c}_{i\gamma}f_i^{eq}\notag\\
	&\qquad=\Delta t^2\partial_{\alpha}\sum_{i=1}^{q}\bigg[cS^{10}_{\alpha}\partial_{\beta}\Big(S^{10}_{\beta}+S^1_{\beta\eta}\bm{c}_{i\eta}\Big)\notag\\
	&\qquad\qquad\qquad\qquad\qquad+S^1_{\alpha\theta}\partial_{\beta}\Big(S^{20}_{\theta\beta}+S^{21}_{\theta\beta\eta}\bm{c}_{i\eta}+S^{2}_{\theta\beta\zeta\mu}\bm{c}_{i\zeta}\bm{c}_{i\mu}\Big)\bigg]\partial_{\gamma}\bm{c}_{i\gamma}f_i^{eq}\notag\\
	&\qquad=\Delta t^2\partial_{\alpha}\sum_{i=1}^{q }
	\bigg[\Big(cS^{10}_{\alpha} \partial_{\beta} S^{10}_{\beta}+S^1_{\alpha\theta}\partial_{\beta}S^{20}_{\theta\beta}\Big)
	+\Big(cS^{10}_{\alpha}\partial_{\beta}S^1_{\beta\eta} +S^1_{\alpha\theta}\partial_{\beta}S^{21}_{\theta\beta\eta}\Big)\bm{c}_{i\eta}\notag\\ &\qquad\qquad\qquad\qquad\qquad+S^1_{\alpha\theta}\partial_{\beta}S^{2}_{\theta\beta\zeta\mu}\bm{c}_{i\zeta}\bm{c}_{i\mu} \bigg] \partial_{\gamma}\bm{c}_{i\gamma}f_i^{eq}. 
\end{align} 
\subsection{NACDE: Derivation of Eq. (\ref{second-diffusive-NACDEs})}\label{appendix-NACDEs-diffusive}
Substituting the moment conditions (\ref{NACDE-moment-conditions}) for NACDE into above Eq. (\ref{zero-order-diffuisve-appendix}), one can derive
\begin{align}
	&\Big[c\bm{\mathcal{W}}_0\bm{m}^{eq}\Big]_1=\partial_{\alpha}B_{\alpha}, \big[\partial_t\bm{m}^{eq}\big]_1=\partial_t\phi,\bm{\tilde{F}}_1=R, \notag\\
	&\Bigg[\Delta tc^2\bm{\mathcal{W}}_0\bigg(\bm{\hat{S}}_N^{-1}+\Big(\frac{b}{2a^2}-1\Big)\bm{I}\bigg)\bm{\mathcal{W}}_0\bm{m}^{eq}\Bigg]_{1} \notag\\
	&\qquad=\Delta t\frac{\partial}{\partial x_{\beta}}\Bigg( cS^{10}_{\beta}\frac{\partial  B_\theta }{ \partial x_{\theta}} 
	+\bigg[S^{1}_{\beta\gamma}+\Big(\frac{b}{2a^2}-1\Big)\bigg]\delta_{\beta\gamma}c_s^2\chi\frac{\partial D_{\gamma \theta}}{ \partial x_{\theta}}\Bigg)+O(\Delta x^2)\notag\\
	&\qquad=\Delta t\chi\frac{\partial}{\partial x_{\beta}}\Bigg( \bigg[S^{1}_{\beta\gamma}+\Big(\frac{b}{2a^2}-1\Big)\delta_{\beta\gamma}\bigg]c_s^2\frac{\partial D_{\gamma \theta}}{ \partial x_{\theta}}\Bigg)+O(\Delta x), 
\end{align} 
then the  first-order ME of the GPMRT-LB model and GPMFD scheme for the NACDE can be given by
\begin{align}\label{NACDEs-first-order-diffusive-aapendix}
	& \partial_t \phi +\partial_{\alpha}B_{\alpha} 
	-\Delta t\frac{\partial}{\partial x_{\beta}}\bigg(\chi c_s^2\bigg[S^{1}_{\beta\gamma}-\Big(\frac{b}{2a^2}-1\Big)\delta_{\beta\gamma}\bigg]\frac{\partial  D_{\beta\theta}  }{\partial  x_{\theta}}\bigg)-R=O(\Delta x).
\end{align} 
Similarly, substituting the moment conditions (\ref{NACDE-moment-conditions}) for NACDE (\ref{NACDE}) into above Eqs. (\ref{first-order-diffuisve-appendix})  and (\ref{second-order-diffuisve-appendix}) yields
\begin{align}
	&\big[c \bm{\mathcal{W}}_0\bm{\hat{S}}_N^{-1}\bm{\tilde{F}}\big]_1=  \partial_{\alpha}\Big[cS^{10}_{\alpha}R+ \big(S^1_{\alpha\beta}-\delta_{\alpha\beta}/2\big)\partial_tB_{\beta}\Big]+O(\Delta x) =c\partial_{\alpha}S^{10}_{\alpha}R+O(1) ,\notag\\
	&\big[c \bm{\mathcal{W}}_0\bm{\hat{S}}_N^{-1}\partial_t\bm{m}^{eq}\big]_1=c\partial_{\alpha}S^{10}_{\alpha}\partial_t\phi+O(1), \big(c\partial_t\bm{\hat{S}}_N^{-1} \bm{\mathcal{W}}_0\bm{m}^{eq}\big)_1 =O(1),\notag\\
	&   \big(c \bm{\mathcal{W}}_0\partial_t\bm{m}^{eq}\big)_1 =O(1),\bigg[ \Delta t^2 c ^3\bm{\mathcal{W}}_0^3\bm{m}^{eq} \bigg]_1   = \partial_{\alpha}\partial_{\beta}\partial_{\gamma}\sum_{i=1}^q\bm{e}_{i\alpha}\bm{e}_{i\beta}\bm{e}_{i\gamma}\notag\\
	&\qquad\times\biggg(w_i\Bigg[\phi+\frac{\Big[c_s^2(D_{\eta\xi}-\phi\delta_{\eta\xi})+C_{\eta\xi}\Big](\bm{c}_{i\eta}\bm{c}_{i\xi}-c_s^2\delta_{\eta\xi})}{2c_s^4}\Bigg]+O(\Delta x)\biggg)\Delta  t^2c^3 =O(\Delta x^2),\notag\\
	& \bigg[ \Delta t^2 c ^3\bm{\mathcal{W}}_0^2\bm{\hat{S}}_N^{-1}\bm{\mathcal{W}}_0\bm{m}^{eq} \bigg]_1 = \Delta t^2\partial_{\alpha}\partial_{\beta} \big(c^2S^{20}_{\alpha\beta} \partial_{\gamma}B_{\gamma} +S^{2}_{\alpha\beta\theta\zeta}\partial_{\gamma}\Delta_{\theta\zeta\gamma\xi}B_{\xi}c_s^2\big)=  O(\Delta x^2),\notag\\
	& \bigg[ \Delta t^2 c ^3\bm{\mathcal{W}}_0\bm{\hat{S}}_N^{-1}\bm{\mathcal{W}}_0^2\bm{m}^{eq} \bigg]_1 =\Delta t^2\partial_{\alpha}\Big(cS^{10}_{\alpha }  \partial_{\beta}\partial_{\gamma} \big(\chi c_s^2 D_{\beta\gamma}+C_{\beta\gamma}\big)+S^{1}_{\alpha \theta}\partial_{\beta}\partial_{\gamma}  \Delta_{\theta\beta\gamma\xi}B_{\xi}c_s^2\Big) \notag\\
	&\quad\:\:\:\:\qquad\qquad\qquad\qquad\quad=\Delta t^2\chi\partial_{\alpha}\big(cS^{10}_{\alpha } c_s^2\partial_{\beta}\partial_{\gamma} D_{\beta\gamma} \big)+O(\Delta x^2),\notag\\
	& \bigg[ \Delta t^2 c ^3\bm{\mathcal{W}}_0\bm{\hat{S}}_N^{-1}\bm{\mathcal{W}}_0\bm{\hat{S}}_N^{-1}\bm{\mathcal{W}}_0\bm{m}^{eq} \bigg]_1= \Delta t^2\partial_{\alpha} 
	\bigg[\Big(cS^{10}_{\alpha} \partial_{\beta}S^{10}_{\beta}+S^1_{\alpha\theta}\partial_{\beta}S^{20}_{\theta\beta}\Big)B_{\gamma} \notag\\
	&\qquad
	+\Big(cS^{10}_{\alpha}\partial_{\beta}S^1_{\beta\eta}  +S^1_{\alpha\theta}\partial_{\beta}S^{21}_{\theta\beta\eta}\Big) \partial_{\gamma}\big(c_s^2\chi D_{\eta\gamma}+C_{\eta\gamma}\big)
	+S^1_{\alpha\theta}\partial_{\beta}S^{2}_{\theta\beta\zeta\mu}c_s^2\partial_{\gamma}\Delta_{\zeta\mu\gamma\delta}B_{\delta}\bigg] \notag\\
	&\qquad=\Delta t^2\chi \partial_{\alpha} 
	\bigg[  
	cS^{10}_{\alpha}\partial_{\beta}S^1_{\beta\eta}   c_s^2\partial_{\gamma}D_{\eta\gamma} \bigg] +O(\Delta x^2), 
\end{align}
thus one can derive the following second-order ME of the GPMRT-LB model and GPMFD scheme for NACDE (\ref{NACDE}):
\begin{align}
	&\partial_t \phi +\partial_{\alpha}B_{\alpha} 
	-\Delta t\frac{\partial}{\partial x_{\beta}}\bigg(\chi c_s^2\bigg[S^{1}_{\beta\gamma}-\Big(\frac{b}{2a^2}-1\Big)\delta_{\beta\gamma}\bigg]\frac{\partial  D_{\beta\theta}  }{\partial  x_{\theta}}\bigg)-R\notag\\
	&+\underbrace{\Delta tc\biggg( \frac{\partial}{\partial x_{\beta}}\Bigg[S^{10}_{\beta}\bigg(\frac{\partial \phi}{\partial t} + \frac{\partial  B_\theta }{ \partial x_{\theta}}-\frac{\partial }{\partial x_{\theta}}\bigg[\chi c_s^2\Delta t\bigg[S^1_{\theta\eta}+\Big(\frac{b}{2a^2}-1\Big)\delta_{\theta\eta}\bigg]\frac{\partial D_{\eta\gamma}}{\partial x_{\gamma}}\bigg] - R  \bigg)\Bigg]\biggg)}_{O(\Delta x^2)} \notag\\
	&=O(\Delta x^2) ,
\end{align} 
where the first-order ME (\ref{NACDEs-first-order-diffusive-aapendix}) has been used.

\subsection{NSEs:  Derivation of Eqs.   (\ref{diff-con-second}), and (\ref{mon-second-diff})}\label{appendix-NSEs-diffusive}	
Substituting the moment conditions (\ref{NSEs-moment-conditions}) for NSEs (\ref{NSEs}) into above Eq. (\ref{zero-order-diffuisve-appendix}), we have
\begin{align}\label{NSE-first-order-diffusive-1}
	&\big[\partial_t\bm{m}^{eq}\big]_1=\partial_t\rho,\Big[c\bm{\mathcal{W}}_0\bm{m}^{eq}\Big]_1=\partial_{\alpha}(\rho u_{\alpha}),\bm{\tilde{F}}_1=0, \notag\\
	&\Bigg[\Delta tc^2\bm{\mathcal{W}}_0\bigg(\bm{\hat{S}}_N^{-1}+\Big(\frac{b}{2a^2}-1\Big)\bm{I}\bigg)\bm{\mathcal{W}}_0\bm{m}^{eq}\Bigg]_{1} \notag\\
	&\qquad=\Delta t\frac{\partial}{\partial x_{\beta}}\Bigg( cS^{10}_{\beta}\frac{\partial  (\rho u_\theta) }{ \partial x_{\theta}} 
	+\bigg[S^{1}_{\beta\gamma}+\Big(\frac{b}{2a^2}-1\Big)\delta_{\beta\gamma}\bigg]\frac{\partial (\rho u_{\gamma}u_{\theta}+\rho c_s^2\delta_{\gamma\theta})}{ \partial x_{\theta}}\Bigg), 
\end{align} 
and for any $\alpha\in\{1\sim d\}$,  
\begin{align}\label{NSE-first-order-diffusive-2-d+1}
	&\Big[c\bm{\mathcal{W}}_0\bm{m}^{eq}\Big]_{\alpha+1}=\partial_{\beta}(\rho u_{\alpha}u_{\beta}+\rho c_s^2\delta_{\alpha\beta})/c,\big( \bm{\tilde{F}}\big)_{\alpha+1}=O(\Delta x),\notag\\
	& \Bigg[\Delta tc^2\bm{\mathcal{W}}_0\bigg(\bm{\hat{S}}_N^{-1}+\Big(\frac{b}{2a^2}-1\Big)\bm{I}\bigg)\bm{\mathcal{W}}_0\bm{m}^{eq}\Bigg]_{\alpha+1} \notag\\
	&=\Delta t\partial_{\beta}\bigg[c S^{20}_{\alpha\beta}\partial_{\gamma}(\rho u_{\gamma})  +  S^{21}_{\alpha\beta\xi_1}\partial_{\gamma}(\rho u_{\xi_1}u_{\gamma}+\rho c_s^2\delta_{\xi_1\gamma })\notag\\
	&\qquad\qquad +\frac{1}{c}\bigg(S^{2}_{\alpha\beta\xi_1\xi_2}+\Big(\frac{b}{2a^2}-1\Big)\delta_{\alpha\xi_1}\delta_{\beta\xi_2}\bigg)\partial_{\gamma}(\rho c_s^2\Delta_{\xi_1\xi_2\gamma\zeta}u_{\zeta})\bigg]\notag\\
	& =\Delta t\partial_{\beta}\bigg[   S^{21}_{\alpha\beta\xi_1}\partial_{\gamma}( \rho c_s^2\delta_{\xi_1\gamma })  \bigg]+O(\Delta x)=O(\Delta x). 
\end{align} 
From above Eqs. (\ref{NSE-first-order-diffusive-1}) and (\ref{NSE-first-order-diffusive-2-d+1}), the first-order MEs of the GPMRT-LB model and GPMFD scheme on conservative moments $m_1=\rho$ and $m_{\alpha+1}=(\rho u_{\alpha})/c$ can be expressed as
\begin{subequations}
	\begin{align}
		& \bigg[   c\bm{\mathcal{W}}_0 \bm{m}^{eq}   \bigg]_1  = \partial_{\beta}(\rho u_{\beta})	=O(\Delta x),\label{continuous-equation-first-order-diffusive-aapendix}\\
		& \bigg[   c\bm{\mathcal{W}}_0 \bm{m}^{eq}   \bigg]_{\alpha+1}  = \frac{1}{c}\partial_{\beta}(\rho c_s^2\delta_{\alpha\beta})  +O(\Delta x)	=O(\Delta x),\label{mometumn-equation-first-order-diffusive-aapendix}
	\end{align}
\end{subequations}
this illustrates that at the diffusive scaling,  $\nabla \rho =O(\Delta x^2)$, which will be used below. Subsequently, substituting the moment conditions (\ref{NSEs-moment-conditions}) for NSEs (\ref{NSEs}) into above Eqs. (\ref{first-order-diffuisve-appendix})  and (\ref{second-order-diffuisve-appendix}) yields
\begin{align}
	&\big[\Delta tc \bm{\mathcal{W}}_0\partial_t\bm{m}^{eq}\big]_1=\Delta t\partial_{\alpha}\partial_t(\rho u_{\alpha})=O(\Delta x^2),\notag\\
	&\big[\Delta tc \bm{\mathcal{W}}_0\bm{\hat{S}}_N^{-1}\partial_t\bm{m}^{eq}\big]_1=\Delta t\Big[c\partial_{\alpha}\big(S^{10}_{\alpha}\partial_t\rho\big]+\partial_{\alpha}\big[S^1_{\alpha\beta}\partial_t(\rho u_{\beta})\big]\Big]= c\partial_{\alpha}\big(S^{10}_{\alpha}\partial_t\rho\big] +O(\Delta x^2),\notag\\
	&\big[\Delta tc\partial_t\bm{\hat{S}}_N^{-1} \bm{\mathcal{W}}_0\bm{m}^{eq}\big]_1  =\partial_ts_0\partial_{\alpha}(\rho u_{\alpha})=O(\Delta x^2),\notag\\
	&\big[\Delta tc \bm{\mathcal{W}}_1\bm{\hat{S}}_N^{-1}\bm{\tilde{F}}\big]_1= \Delta t\partial_{\alpha}\Big[cS^{10}_{\alpha}\frac{\Delta t}{2}\partial_{\gamma}(\rho \bm{\hat{F}}_{x_{\gamma}})+ \ S^1_{\alpha\beta}(\rho \bm{\hat{F}}_{x_{\beta}})  +O(\Delta x)\Big]= O(\Delta x^2) ,\notag\\
	& \bigg[ \Delta t^2 c ^3\bm{\mathcal{W}}_0^3\bm{m}^{eq} \bigg]_1  =\Delta  t^2c^3 \partial_{\alpha}\partial_{\beta}\partial_{\gamma}\sum_{i=1}^q\bm{e}_{i\alpha}\bm{e}_{i\beta}\bm{e}_{i\gamma}\big(w_i\rho+O(\Delta x)\big) =O(\Delta x^2),\notag\\
	& \bigg[ \Delta t^2 c ^3\bm{\mathcal{W}}_0^2\bm{\hat{S}}_N^{-1}\bm{\mathcal{W}}_0\bm{m}^{eq} \bigg]_1=\Delta t^2\partial_{\alpha}\partial_{\beta} \big(c^2S^{20}_{\alpha\beta} \partial_{\gamma}(\rho u_{\gamma})+cS^{21}_{\alpha\beta\xi_1} \partial_{\gamma}(\rho u_{\xi_1}u_{\gamma}+c_s^2\rho \delta_{\xi_1\gamma})\notag\\
	&\qquad\qquad\qquad\qquad\qquad\qquad+S^{2}_{\alpha\beta\xi_1\xi_2}\partial_{\gamma}(\rho c_s^2\Delta_{\xi_1\xi_2\gamma\zeta}u_{\zeta})\big)= O(\Delta x^2),\notag\\ 
	& \bigg[ \Delta t^2 c ^3\bm{\mathcal{W}}_0\bm{\hat{S}}_N^{-1}\bm{\mathcal{W}}_0^2\bm{m}^{eq} \bigg]_1  =\Delta t^2\partial_{\alpha}\sum_{k=1}^{q }\bigg(cS^{10}_{\alpha}+S^{1}_{\alpha\xi_1}\bm{c}_{k\xi_1}\bigg)\partial_{\beta}\bm{c}_{k\beta}\partial_{\gamma}\bm{c}_{k\gamma}f_k^{eq}\notag\\
	&\qquad =\Delta t^2\partial_{\alpha}\big(cS^{10}_{\alpha } \partial_{\beta}\partial_{\gamma} (\rho u_{\beta}u_{\gamma}+c_s^2\rho \delta_{\beta\gamma})+S^{1}_{\alpha \xi_1}\partial_{\beta}\partial_{\gamma} (\rho c_s^2 \Delta_{\xi_1\beta\gamma\zeta}u_{\zeta})\big) = O(\Delta x^2),\notag\\
	& \bigg[ \Delta t^2 c ^3\bm{\mathcal{W}}_0\bm{\hat{S}}_N^{-1}\bm{\mathcal{W}}_0\bm{\hat{S}}_N^{-1}\bm{\mathcal{W}}_0\bm{m}^{eq} \bigg]_1 =  \Delta t^2\partial_{\alpha} 
	\bigg[\Big(cS^{10}_{\alpha} \partial_{\beta}S^{10}_{\beta}+S^1_{\alpha\xi_1}\partial_{\beta}S^{20}_{\xi_1\beta}\Big)\partial_{\gamma}(\rho u_{\gamma})\notag\\
	&\qquad\qquad\qquad\qquad\qquad\qquad
	+\Big(cS^{10}_{\alpha}\partial_{\beta}S^1_{\beta\zeta_1}  +S^1_{\alpha\xi_1}\partial_{\beta}S^{21}_{\xi_1\beta\zeta_1}\Big)\partial_{\gamma}(\rho u_{\zeta_1}u_{\gamma}+c_s^2\rho \delta_{\zeta_1\gamma})\notag\\
	&\qquad\qquad\qquad\qquad\qquad\qquad\qquad\qquad 
	+S^1_{\alpha\xi_1}\partial_{\beta}S^{2}_{\xi_1\beta\zeta_1\zeta_2}\partial_{\gamma}(c_s^2\rho\Delta_{\zeta_1\zeta_2\gamma\zeta}u_{\zeta})\bigg]   = O(\Delta x^2). 
\end{align} 
Additionally, one can also  derive the following second-order ME of the GPMRT-LB model and GPMFD scheme for the continuous equation (\ref{continuous-equation}): 
\begin{align}\label{second-order-contonuous-appendix}
	&  \partial_t \rho +\partial_{\alpha}(\rho u_{\alpha})   +\underbrace{\Delta tc\biggg( \frac{\partial}{\partial x_{\beta}}\Bigg[S^{10}_{\beta}\bigg(\frac{\partial \rho}{\partial t} + \frac{\partial  (\rho u_\gamma) }{ \partial x_{\gamma}}  \bigg)\Bigg]\biggg)}_{O(\Delta x^2)} =O(\Delta x^2),  
\end{align}  
where the first-order ME (\ref{continuous-equation-first-order-diffusive-aapendix}) has been applied.

For the GPMRT-LB model and the GPMFD  scheme for the momentum equation (\ref{momentum-equation}), one can obtain   
\begin{align}\label{second-order-remain}
	&\big[-\partial_t\bm{m}^{eq}\big]_{\alpha+1}=-\frac{1}{c}\partial_t(\rho u_{\alpha} ),\big[-c\bm{\mathcal{W}}_0\bm{m}^{eq}\big]_{\alpha+1}=-\frac{1}{c}\partial_{\beta}(\rho u_{\alpha}u_{\beta}+\rho c_s^2\delta_{\alpha\beta}),\notag\\
	&\big[\bm{\tilde{F}}\big]_{\alpha+1}=\hat{F}_{x_{\alpha}},\big[\Delta tc^2 \bm{\mathcal{W}}_0^2\bm{m}^{eq}\big]_{\alpha+1}= -\Delta t \frac{1}{c} \partial_{\beta}\partial_{\gamma}(\rho c_s^2\Delta_{\alpha\beta\gamma\zeta}u_{\zeta}),\notag\\
	&\big[\Delta tc\bm{\mathcal{W}}_0\bm{\hat{S}}_N^{-1}\partial_t\bm{m}^{eq}\big]_{\alpha+1}=\Delta tc\sum_{k=1}^{q}\Big[\sum_{j=0}^{q-1}\bm{e}_{j\alpha}\bm{e}_{j\beta}\partial_{\beta}\bm{\overline{\Lambda}}_{jk}\Big]\partial_tf_k^{eq}\notag\\
	&\qquad\qquad\qquad\qquad\qquad	=\Delta t  \partial_{\beta}\Big[cS^{20}_{\alpha\beta}\partial_t\rho+S^{21}_{\alpha\beta\xi_1}\partial_t(\rho u_{\xi_1})+\frac{1}{c}S^{2}_{\alpha\beta\xi_1\xi_2}\partial_t(\rho c_s^2\delta_{\xi_1\xi_2}) \Big],\notag\\
	&\big[\Delta tc^2\bm{\mathcal{W}}_0\bm{\hat{S}}_N^{-1}\bm{\mathcal{W}}_0\bm{m}^{eq}\big]_{\alpha+1}=\Delta tc^2\sum_{k=1}^{q}\Big[\sum_{j=1}^{q}\bm{e}_{j\alpha}\bm{e}_{j\beta}\partial_{\beta}\bm{\overline{\Lambda}}_{jk}\Big]\bm{e}_{k\gamma}\partial_{\gamma}f_k^{eq}\notag\\
	&	=\Delta t  \partial_{\beta}\Big[cS^{20}_{\alpha\beta}\partial_{\gamma}(\rho u_{\gamma})+S^{21}_{\alpha\beta\xi_1}\partial_{\gamma}(\rho  u_{\xi_1}u_{\gamma}+\rho c_s^2\delta_{\xi_1\gamma})+\frac{1}{c}S^{2}_{\alpha\beta\xi_1\xi_2}\partial_{\gamma}(\rho c_s^2\Delta_{\xi_1\xi_2\gamma\zeta}u_{\zeta}) \Big] , 
\end{align}  
and
\begin{align}\label{third-order-remain}
	&\big[-\Delta tc\bm{\mathcal{W}}_0\partial_t\bm{m}^{eq}\big]_{\alpha+1}\notag\\
	&\qquad\qquad\qquad=-\Delta t\frac{1}{c}\partial_{\beta}\partial_t(\rho u_{\alpha}u_{\beta}+\rho c_s^2\delta_{\alpha\beta}) =-\Delta t\frac{1}{c}\partial_{\beta}\partial_t(\rho c_s^2\delta_{\alpha\beta})+O(\Delta x^3),\notag\\
	&\big[-\Delta tc\bm{\mathcal{W}}_0\bm{\hat{S}}_N^{-1}\bm{\tilde{F}}\big]_{\alpha+1}=-c\Delta t\sum_{k=1}^{q }\Big[\sum_{j=1}^{q }\bm{e}_{j\alpha}\bm{e}_{j\beta}\partial_{\beta}\bm{\overline{\Lambda}}_{jk}\Big]\big(F_k+\frac{\Delta t}{2}D_kF_k\big)\notag\\
	&\qquad\qquad	=-c\Delta t\sum_{k=1}^{q }\partial_{\beta}\Big[S^{20}_{\alpha\beta}+S^{21}_{\alpha\beta\xi_1}\bm{e}_{k\xi_1}+S^{2}_{\alpha\beta\xi_1\xi_2}\bm{e}_{k\xi_1}\bm{e}_{k\xi_2}\Big]\Big(F_k+\frac{\Delta t}{2}D_k F_k\Big)\notag\\
	&\qquad\qquad	=- \Delta t \partial_{\beta}\Big[ S^{21}_{\alpha\beta\xi_1}(\rho \hat{F}_{x_{\alpha}})\Big]+O(\Delta x^3),\notag\\
	&\big[\Delta t c\partial_t\bm{\hat{S}}_N^{-1}\bm{\mathcal{W}}_0\bm{m}^{eq}\big]_{\alpha+1}=\Delta tc\partial_t\sum_{k=0}^{q-1}\Big[\sum_{j=0}^{q-1}\bm{e}_{j\alpha}\bm{\overline{\Lambda}}_{jk}\Big]\bm{e}_{k\beta}\partial_{\beta}f_k^{eq}\notag\\
	&\qquad\qquad=\Delta t\partial_tS^{10}_{\alpha}\partial_{\beta}(\rho u_{\beta})+\Delta t\frac{1}{c}\partial_tS^1_{\alpha\xi_1}\partial_{\beta}(\rho u_{\xi_1}u_{\beta}+\rho c_s^2\delta_{\xi_1\beta})\notag\\
	&\qquad\qquad=\Delta t\partial_t\Big[S^{10}_{\alpha}\partial_{\beta}(\rho u_{\beta})\Big]+O(\Delta x^3),\notag\\
	&\big[\Delta t^2c^3\bm{\mathcal{W}}_0^2\bm{\hat{S}}_N^{-1}\bm{\mathcal{W}}_0\bm{m}^{eq}\big]_{\alpha+1}=\Delta t^2c^3\sum_{k=1}^{q}\Big[\sum_{j=1}^{q }\bm{e}_{j\alpha}\bm{e}_{j\beta}\partial_{\beta}\bm{e}_{j\gamma}\partial_{\gamma}\bm{\overline{\Lambda}}_{jk}\Big]\bm{e}_{k\zeta}\partial_{k\zeta}f_k^{eq}\notag\\
	&=\Delta t^2 \partial_{\beta}\partial_{\gamma}\Big[c^2S^{30}_{\alpha\beta\gamma}\partial_{\zeta}(\rho u_{\zeta})  \Big]+O(\Delta x^3),\notag\\
	&\big[\Delta t^2c^3\bm{\mathcal{W}}_0\bm{\hat{S}}_N^{-1}\bm{\mathcal{W}}_0^2\bm{m}^{eq}\big]_{\alpha+1}=\Delta t^2c^3 \sum_{k=1}^{q}\Big[\sum_{j=1}^{q}\bm{e}_{j\alpha}\bm{e}_{j\beta}\partial_{\beta}\bm{\overline{\Lambda}}_{jk}\Big]\bm{e}_{k\gamma}\bm{e}_{k\zeta}\partial_{\gamma}\partial_{\zeta}f_k^{eq}\notag\\
	&=\Delta t^2 \partial_{\beta}\Big[cS^{20}_{\alpha\beta}\partial_{\gamma}\partial_{\zeta}\times O(1)+S^{21}_{\alpha\beta\xi_1}\partial_{\gamma}\partial_{\zeta}
	(\rho c_s^2\Delta_{\xi_1\gamma\zeta\chi}u_{\chi})+\frac{1}{c}S^{2}_{\alpha\beta\xi_1\xi_2}\partial_{\gamma}\partial_{\zeta}\notag\\
	&\qquad \times [\rho c_s^3\Delta_{\xi_1\xi_2\gamma\zeta}+O(1/\Delta x^2)]\Big]=\Delta t^2 \partial_{\beta}\Big[ S^{21}_{\alpha\beta\xi_1}\partial_{\gamma}\partial_{\zeta}
	(\rho c_s^2\Delta_{\xi_1\gamma\zeta\chi}u_{\chi}) \Big]+O(\Delta x^3),\notag\\
	&\big[-\Delta t^2c^3\bm{\mathcal{W}}_0\bm{\hat{S}}_N^{-1}\bm{\mathcal{W}}_0\bm{\hat{S}}_N^{-1}\bm{\mathcal{W}}_0\bm{m}^{eq}\big]_{\alpha+1} \notag\\
	&\qquad\qquad=-\sum_{i=1}^{q }\bigg[\sum_{k=1}^{q }\Big[\sum_{j=1}^{q }\bm{e}_{j\alpha}\bm{e}_{j\beta}\partial_{\beta}\bm{\overline{\Lambda}}_{jk}\Big]\bm{e}_{k\gamma}\partial_{\gamma}\bm{\overline{\Lambda}}_{ki}\bigg]\bm{e}_{i\theta}\partial_{\theta}f_i^{eq}\Delta t^2c^3\notag\\
	&\qquad\qquad=-\Delta t^2 \partial_{\beta} \bigg[ 
	S^{20}_{\alpha\beta}\partial_{\gamma}
	\big(
	c^2S^{10}_{\gamma}\partial_{\theta}(\rho u_{\theta})\big) +S^{21}_{\alpha\beta\xi_1}\partial_{\gamma}
	\Big(
	c^2S^{20}_{\xi_1\gamma}\partial_{\theta}(\rho u_{\theta})\notag\\
	&\qquad\qquad
	+S^{2}_{\xi_1\gamma\zeta_1\zeta_2}\partial_{\theta}(\rho c_s^2\Delta_{\zeta_1\zeta_2\theta\eta}u_{\eta})\Big) +S^{2}_{\alpha\beta\xi_1\xi_2}\partial_{\gamma}\Big(
	c^2S^{30}_{\xi_1\xi_2\gamma}\partial_{\theta}(\rho u_{\theta})
	\Big)
	\bigg]+O(\Delta x^3),
\end{align} 
where $\alpha\in\{1\sim d\}$. Together with Eqs. (\ref{second-order-remain}) and (\ref{third-order-remain}), the second-order ME of GPMRT-LB model and GPMFD scheme on conservative moment $m_{\alpha+1}=(\rho u_{\alpha})/c$ can be obtained,
\begin{align}\label{first-order-momentum-equation-appendix}
	&\Big[  \big(\partial_t\bm{I}+c \bm{\mathcal{W}}_0 \big)\bm{m}^{eq}\Big]_{\alpha+1}\notag\\
	&=\bigg[\bm{\tilde{F}}+\Delta t  c \bm{\mathcal{W}}_0 \Big(\bm{\hat{S}}_N^{-1}+\Big(\frac{b}{2a^2}-1\Big)\bm{I}\Big)  c \bm{\mathcal{W}}_0   \bm{m}^{eq}  +\Delta t    c \bm{\mathcal{W}}_0\partial_t \bm{\hat{S}}_N^{-1}  \bm{m}^{eq}\bigg]_{\alpha+1}\notag\\
	&=\bigg[ \partial_t(\rho u_{\alpha})+ \partial_{\beta}\big(\rho u_{\alpha}u_{\beta}+\rho c_s^2\delta_{\alpha\beta}) -\Delta t \partial_{\beta}S^{2}_{\alpha\beta\xi_1\xi_2}\partial_t(\rho c_s^2\delta_{\xi_1\xi_2}) \notag\\
	& - \rho  {\hat{F}}_{x_{\alpha}}-\Delta t  \partial_{\beta} \bigg[S^{2}_{\alpha\beta\xi_1\xi_1}+\Big(\frac{b}{2a^2}-1\Big)\delta_{\alpha\xi_1}\delta_{\beta\xi_2} \bigg] \partial_{\gamma}(\rho c_s^2\Delta_{\xi_1\xi_2\gamma\zeta}u_{\zeta}) \bigg]\frac{1}{c}=O(\Delta x^2), 
\end{align} 
where the second-order ME of the GPMRT-LB model and GPMFD scheme for the continuous equation (\ref{second-order-contonuous-appendix}) has been used. Actually, at the diffusive
scaling, 
$c=O(1/\Delta x)$, and above ME of the GPMRT-LB model and GPMFD scheme for the momentum equation (\ref{momentum-equation}) is first-order accurate. This also explains why we further conduct the expansion of $\bm{\Xi}$ up to $\bm{\Xi}^{(4)}$ in Eq. (\ref{five-order-diffusive}). To obtain a second-order ME of the GPMRT-LB model and GPMFD scheme for the momentum equation (\ref{momentum-equation}), it is necessary to compute Eq. (\ref{modified-equation-third-order-diffusive}), and after some manipulations, one can derive 
\begin{align}\label{third-order-remain-2}
	&\Big[-\frac{\Delta t}{2}\partial_{tt}\bm{m}^{eq}\Big]_{\alpha +1}=-\frac{1}{c}\frac{\Delta t}{2}\partial_{tt}(\rho u_{\alpha}),\notag\\
	&\big[\Delta t\partial_t\bm{\hat{S}}_N^{-1}\partial_t\bm{m}^{eq}\big]_{\alpha +1}=\Delta t\partial_t\sum_{k=0}^{q-1}\Big[\sum_{j=0}^{q-1}\bm{e}_{j\alpha}\bm{\overline{\Lambda}}_{jk}\Big]\partial_tf_k^{eq}\notag\\
	&\qquad\qquad=\Delta t\partial_t\sum_{k=0}^{q-1}\Big[S^{10}_{\alpha}+S^1_{\alpha\xi_1}\bm{e}_{k\xi_1}\Big]\partial_tf_k^{eq}=\Delta t\partial_t\Big[S^{10}_{\alpha}\partial_t\rho\Big]+O(\Delta x^3),\notag\\
	&\big[ \Delta t^2c^2\bm{\mathcal{W}}_0^2\bm{\hat{S}}_N^{-1}\partial_t\big]_{\alpha+1}= \Delta t^2c^2\sum_{k=1}^{q}\Big[\sum_{j=1}^{q }\bm{e}_{j\alpha}\bm{e}_{j\beta}\partial_{\beta}\bm{e}_{j\gamma}\partial_{\gamma}\bm{\overline{\Lambda}}_{jk}\Big]\partial_tf_k^{eq}\notag\\
	&\qquad\qquad =\Delta t^2c^2\partial_{\beta}\partial_{\gamma} S^{30}_{\alpha\beta\gamma}\partial_t\rho +O(\Delta x^3),\notag\\
	&\big[\Delta t^2c^2\bm{\mathcal{W}}_0\bm{\hat{S}}_N^{-1} \bm{\mathcal{W}}_0\partial_t\bm{m}^{eq}\big]_{\alpha+1}=\Delta t^2c^2\sum_{k=1}^{q}\Big[\sum_{j=1}^{q}\bm{e}_{j\alpha}\bm{e}_{j\beta}\partial_{\beta}\bm{\overline{\Lambda}}_{jk}\Big]\bm{e}_{k\gamma}\partial_{\gamma}\partial_tf_k^{eq}\notag\\
	&\qquad=\Delta t^2 \partial_{\beta} \Big[ cS^{20}_{\alpha\beta}\partial_{\gamma}\partial_t\times O(1)+S^{21}_{\alpha\beta\xi_1}\partial_{\gamma}\partial_t(\rho c_s^2\delta_{\xi_1\gamma})+\frac{1}{c}S^{2}_{\alpha\beta\xi_1\xi_2}\partial_{\gamma}\partial_t\times O(1/\Delta x^2) \Big]\notag\\
	&\qquad=\Delta t^2\partial_{\beta}\Big[S^{21}_{\alpha\beta\xi_1}\partial_{\gamma}\partial_t(\rho c_s^2\delta_{\xi_1\gamma})\Big]+O(\Delta x^3) ,\notag\\
	&\big[-\Delta t^2c^2\bm{\mathcal{W}}_0\bm{\hat{S}}_N^{-1}\bm{\mathcal{W}}_0\bm{\hat{S}}_N^{-1}\partial_t\bm{m}^{eq}\big]_{\alpha+1}=\sum_{i=1}^{q }\bigg[\sum_{k=1}^{q }\Big[\sum_{j=1}^{q }\bm{e}_{j\alpha}\bm{e}_{j\beta}\partial_{\beta}\bm{\overline{\Lambda}}_{jk}\Big]\bm{e}_{k\gamma}\partial_{\gamma}\bm{\overline{\Lambda}}_{ki}\partial_tf_i^{eq}\notag\\
	& =-\Delta t^2 \partial_{\beta}\bigg[ S^{20}_{\alpha\beta}
	\partial_{\gamma}(c^2S^{10}_{\gamma}\partial_t\rho )  
	+S^{21}_{\alpha\beta\xi_1}
	\partial_{\gamma}	(c^2S^{20}_{\xi_1\gamma}\partial_t\rho ) +S^{21}_{\alpha\beta\xi_1}
	\partial_{\gamma}\Big[S^{2}_{\xi_1\gamma\zeta_1\zeta_2}\partial_t(\rho c_s^2\delta_{\zeta_1\zeta_2})\Big]\notag\\
	&\qquad\qquad\qquad\qquad\qquad +S^{2}_{\alpha\beta\xi_1\xi_2}
	\partial_{\gamma}	(c^2S^{30}_{\xi_1\xi_2\gamma}\partial_t\rho)
	\bigg] +O(\Delta x^3),
\end{align} 

\begin{align}\label{w3}
	\big[\Delta t^2c^3\bm{\mathcal{W}}_0^3\bm{m}^{eq}\big]_{\alpha+1}&=\frac{1}{c}\Delta t^2 \partial_{\beta}\partial_{\gamma}\partial_{\zeta}\sum_{i=1}^q\bm{c}_{i\alpha}\bm{c}_{i\beta}\bm{c}_{i\gamma}\bm{c}_{i\zeta}\Big[w_i\rho+\bm{c}_{i\eta}\partial_{\eta}+O(\Delta x^2)\Big]\notag\\
	&=O(\Delta x^3),
\end{align} 
\begin{align}\label{w4}
	&\big[\Delta t^3c^4\bm{\mathcal{W}}_0^4\bm{m}^{eq}\big]_{\alpha+1}=\frac{1}{c}\Delta t^3 \partial_{\beta}\partial_{\gamma}\partial_{\zeta}\partial_{\eta}\sum_{i=1}^q\bm{c}_{i\alpha}\bm{c}_{i\beta}\bm{c}_{i\gamma}\bm{c}_{i\zeta}\bm{c}_{i\eta}\Big[w_i\rho+O(\Delta x)\Big]=O(\Delta x^3),
\end{align} 
\begin{align}\label{w2t}
	&\big[-\Delta t^2c^2 \bm{\mathcal{W}}_0^2\partial_t\bm{m}^{eq}\big]_{\alpha+1} =-\Delta t^2\frac{1}{c} \partial_{\beta}\partial_{\gamma}\partial_t(\rho c_s^2\Delta_{\alpha\beta\gamma\zeta}u_{\zeta})=O(\Delta x^3),
\end{align} 
\begin{align}\label{tsf}
	&\big[-\Delta t\partial_t\bm{\hat{S}}_N^{-1}\bm{\tilde{F}}\big]_{\alpha+1} =-\Delta t\sum_{k=1}^{q}\Big[\sum_{j=1}^{q }\bm{e}_{j\alpha} \partial_t\bm{\overline{\Lambda}}_{jk}\Big]\big(F_k+\frac{\Delta t}{2}D_kF_k\big)  =O(\Delta x^3),
\end{align}  
\begin{align}\label{w2sf}
	&\big[-\Delta t^2c^2\bm{\mathcal{W}}_0^2\bm{\hat{S}}_N^{-1}\bm{\tilde{F}}\big]_{\alpha+1} =-\Delta t^2c^2 \sum_{k=1}^{q }\partial_{\beta}\partial_{\gamma}\Big[\sum_{j=1}^{q}\bm{e}_{j\alpha}\bm{e}_{j\beta}\bm{e}_{j\gamma} \bm{\overline{\Lambda}}_{jk}\Big]\big(F_k+\frac{\Delta t}{2}D_kF_k\big)\notag\\
	&\qquad=-c\Delta t^2\sum_{k=1}^{q }\partial_{\beta}\partial_{\gamma}\Big[S^{30}_{\alpha\beta\gamma}+S^{31}_{\alpha\beta\gamma\xi_1}\bm{e}_{k\xi_1} +S^{33}_{\alpha\beta\gamma\xi_1\xi_2\xi_3}\bm{e}_{k\xi_1}\bm{e}_{k\xi_2}\bm{e}_{k\xi_3}\Big]\notag\\
	&\qquad\times\big(F_k+\frac{\Delta t}{2}D_kF_k\big) =O(\Delta x^3),
\end{align}  
\begin{align}\label{wswsf}
	&\big[-\Delta t^2c^2 \bm{\mathcal{W}}_0\bm{\hat{S}}_N^{-1}\bm{\mathcal{W}}_0\bm{\hat{S}}_N^{-1}\bm{\tilde{F}}\big]_{\alpha+1}\notag\\
	&\qquad=-\Delta t^2c^2 \sum_{i=1}^{q }\bigg[\sum_{k=1}^{q }\partial_{\beta}\Big[\sum_{j=1}^{q }\bm{e}_{j\alpha}\bm{e}_{j\beta}  \bm{\overline{\Lambda}}_{jk}\Big]\bm{e}_{j\gamma}\partial_{\gamma}\bm{\overline{\Lambda}}_{ki}\bigg]\big(F_i+\frac{\Delta t}{2}D_iF_i\big)\notag\\
	&\qquad=- \Delta t^2 \partial_{\beta}\bigg[  S^{20}_{\alpha\beta}\partial_{\gamma}S^1_{\gamma\zeta_1}(\rho \hat{F}_{x_{\zeta_1}})+S^{21}_{\alpha\beta\xi_1}\partial_{\gamma}\Big( S^{21}_{\xi_1\gamma\zeta_1}(\rho \hat{F}_{x_{\zeta_1}})+\frac{1}{c}S^{2}_{\xi_1\gamma\zeta_1\zeta_2}\times O(1)\Big)\notag\\
	&\qquad+S^{2}_{\alpha\beta\xi_1\xi_2}\partial_{\gamma}\Big(S^{31}_{\xi_1\xi_2\gamma\zeta_1}(\rho \hat{F}_{x_{\zeta_1}}) +\frac{1}{c^2}S^{33}_{\xi_1\xi_2\gamma\zeta_1\zeta_2\zeta_3}\times(\rho c_s^2\Delta_{\zeta\zeta_2\zeta_3\eta}\hat{F}_{x_{\eta}})\Big)\bigg]+O(\Delta x^6)\notag\\
	&\qquad =O(\Delta x^4),
\end{align}  
\begin{align}\label{wtsw}
	&\big[\Delta t^2c^2\bm{\mathcal{W}}_0\partial_t\bm{\hat{S}}_N^{-1} \bm{\mathcal{W}}_0\bm{m}^{eq}\big]_{\alpha+1}
	=\Delta t^2c^2\sum_{k=1}^{q }\Big[\sum_{j=1}^{q }\bm{e}_{j\alpha}\bm{e}_{j\beta}\partial_{\beta}\partial_t\bm{\overline{\Lambda}}_{jk}\Big]\bm{e}_{k\gamma}\partial_{\gamma}f_k^{eq}\notag\\
	&\qquad =\Delta t^2 \partial_{\beta}\partial_t\Big[ cS^{20}_{\alpha\beta}\partial_{\gamma}\times O(1)+S^{21}_{\alpha\beta\xi_1}\partial_{\gamma}\times O(1)+\frac{1}{c}S^{2}_{\alpha\beta\xi_1\xi_2}\partial_{\gamma}\times O(1/\Delta x^2) \Big] \notag\\
	&\qquad=O(\Delta x^3),   
\end{align}  
\begin{align}\label{tsw2}
	&\big[\Delta t^2c^2\partial_t\bm{\hat{S}}_N^{-1}\bm{\mathcal{W}}_0^2\bm{m}^{eq}\big]_{\alpha+1}=\Delta t^2c^2 \partial_t\sum_{k=1}^{q}\Big[\sum_{j=0}^{q-1}\bm{e}_{j\alpha}\bm{\overline{\Lambda}}^{-1}_{jk}\Big]\bm{e}_{k\beta}\partial_{\beta}\bm{e}_{k\gamma}\partial_{\gamma}f_k^{eq}\notag\\
	&\qquad=\Delta t^2 \partial_t\sum_{k=1}^{q }\Big[S^{10}_{\alpha}\partial_{\beta}\partial_{\gamma}\times O(1) +\frac{1}{c}S^1_{\alpha\xi_1}\partial_{\beta}\partial_{\gamma}\times O(1/\Delta x^2)\Big]=O(\Delta x^3) ,
\end{align}  
\begin{align}\label{w2sw2}
	& \big[\Delta t^3c^4\bm{\mathcal{W}}_0^2\bm{\hat{S}}_N^{-1}\bm{\mathcal{W}}_0^2\big]_{\alpha+1}= \Delta t^3c^4\sum_{i=1}^{q }\Big[\sum_{j=1}^{q}\bm{e}_{j\alpha}\bm{e}_{j\beta}\bm{e}_{j\gamma}\partial_{\beta}\partial_{\gamma}\bm{\overline{\Lambda}}_{jk}\Big]\partial_{\eta}\bm{e}_{k\eta}\partial_{\zeta}\bm{e}_{k\zeta}f_k^{eq}\notag\\ 
	&= \Delta t^3 \sum_{i=1}^{q}\Big[ c^2S^{30}_{\alpha\beta\gamma}\partial_{\eta}\partial_{\zeta}\times O(1)+cS^{31}_{\alpha\beta\gamma\xi_1}\partial_{\eta}\partial_{\zeta}\times O(1/\Delta x^2)  \notag\\
	&\qquad +\frac{1}{c}S^{33}_{\alpha\beta\gamma\xi_1\xi_2\xi_3}\partial_{\eta}\partial_{\zeta}\times O(1/\Delta x^4)\Big]=O(\Delta x^3),
\end{align} 

\begin{align}\label{wswsw2}
	&\big[-\Delta t^3c^4\bm{\mathcal{W}}_0\bm{\hat{S}}_N^{-1}\bm{\mathcal{W}}_0\bm{\hat{S}}_N^{-1}\bm{\mathcal{W}}_0^2\bm{m}^{eq}\big]_{\alpha+1}\notag\\
	&=-\Delta t^3c^4\sum_{i=1}^{q}\bigg[\sum_{k=1}^{q}\Big[\sum_{j=1}^{q}
	\bm{e}_{j\alpha}\bm{e}_{j\beta}\partial_{\beta}\bm{\overline{\Lambda}}_{jk}\Big]\bm{e}_{k\gamma}\partial_{\gamma}\bm{\overline{\Lambda}}_{ki}\bigg]\bm{e}_{i\zeta}\partial_{\zeta}\bm{e}_{i\eta}\partial_{\eta}f_i^{eq} \notag\\
	&=-\Delta t^3 \partial_{\beta} \bigg[ 
	S^{20}_{\alpha\beta}\partial_{\gamma}\Big(c^2S^{10}_{\gamma}\partial_{\zeta}\partial_{\eta}\times O(1)+cS^1_{\gamma\zeta_1}\partial_{\zeta}\partial_{\eta}\times O(1/\Delta x^2)\Big)\notag\\
	&\qquad+S^{21}_{\alpha\beta\xi_1}\partial_{\gamma}\Big(c^2S^{20}_{\xi_1\gamma}\partial_{\zeta}\partial_{\eta}\times O(1)+cS^{21}_{\xi_1\gamma\zeta_1}\partial_{\zeta}\partial_{\eta}\times O(1/\Delta x^2)\notag\\
	&\qquad 
	+S^{2}_{\xi_1\gamma\zeta_1\zeta_2}\partial_{\zeta}\partial_{\eta}\big[\rho c_s^4\Delta_{\zeta_1\zeta_2\zeta\eta}+O(1/\Delta x^2)\big]\Big)  +S^{2}_{\alpha\beta\xi_1\xi_2}\partial_{\gamma}\Big(c^2S^{30}_{\xi_1\xi_2\gamma}\partial_{\zeta}\partial_{\eta}\times O(1)\notag\\
	&\qquad+cS^{31}_{\xi_1\xi_2\gamma\zeta_1}\partial_{\zeta}\partial_{\eta}\times O(1/\Delta x^2) 	 
	+\frac{1}{c}S^{33}_{\xi_1\xi_2\gamma\zeta_1\zeta_2\zeta_3}\partial_{\zeta}\partial_{\eta}\times O(1/\Delta x^4)\Big) \bigg] =O(\Delta x^3),
\end{align}  
\begin{align}\label{wts2w}
	&\big[-\Delta t^2c^2\bm{\mathcal{W}}_0\partial_t\bm{\hat{S}}_N^{-1}\bm{\hat{S}}_N^{-1}\bm{\mathcal{W}}_0\bm{m}^{eq}\big]_{\alpha+1}\notag\\
	&\qquad=-\Delta t^2c^2\sum_{i=1}^{q }\bigg[\sum_{k=0}^{q-1}\Big[\sum_{j=1}^{q}\bm{e}_{j\alpha}\bm{e}_{j\beta}\partial_{\beta}\partial_t\bm{\overline{\Lambda}}_{jk}\Big]\bm{\overline{\Lambda}}_{ki}\bigg]\bm{e}_{i\gamma}\partial_{\gamma}f_i^{eq}\notag\\
	&\qquad=-\Delta t^2  \partial_{\beta}\partial_t\bigg[ 
	S^{20}_{\alpha\beta}cs_0 \partial_{\gamma}\times O(1)
	+S^{21}_{\alpha\beta\xi_1}(cS^{10}_{\xi_1} \partial_{\gamma}\times O(1)+S^1_{\xi_1\zeta_1} \partial_{\gamma}\times O(1))\notag\\
	&\qquad  +S^{2}_{\alpha\beta\xi_1\xi_2}\big(cS^{20}_{\xi_1\xi_2}\bm{c}_{i\gamma}\partial_{\gamma}+S^{21}_{\xi_1\xi_2\zeta_1} \partial_{\gamma}\times O(1)\notag\\
	&\qquad\qquad+\frac{1}{c}S^{ 2}_{\xi_1\xi_2\zeta_1\zeta_2} \partial_{\gamma}\times O(1/\Delta x^2)\big)\bigg] =O(\Delta x^3), 
\end{align}  
\begin{align}\label{wsw2sw}
	&\big[-\Delta t^3c^4\bm{\mathcal{W}}_0\bm{\hat{S}}_N^{-1}\bm{\mathcal{W}}_0^2\bm{\hat{S}}_N^{-1}\bm{\mathcal{W}}_0\bm{m}^{eq}\big]_{\alpha+1}\notag\\
	& =-\Delta t^3c^4 \sum_{i=1}^{q}\bigg[\sum_{k=1}^{q}\Big[\sum_{j=0}^{q-1}\bm{e}_{j\alpha}\bm{e}_{j\beta}\partial_{\beta}\bm{\overline{\Lambda}}_{jk}\Big]\bm{e}_{k\zeta}\partial_{\zeta}\bm{e}_{k\eta}\partial_{\eta}\bm{\overline{\Lambda}}_{ki}\bigg] \bm{e}_{i\theta}\partial_{\theta}f_i^{eq}\notag\\
	&=-\Delta t^3 \partial_{\beta}\bigg[ 
	S^{20}_{\alpha\beta}\partial_{\zeta}\partial_{\eta}\Big(c^3S^{20}_{\zeta\eta}\partial_{\theta}\times O(1)+c^2S^{21}_{\zeta\eta\zeta_1} \partial_{\theta}\times O(1)+cS^{2}_{\zeta\eta\zeta_1\zeta_2}\partial_{\theta}\times O(1/\Delta x^2)\Big)\notag\\
	&  +S^{21}_{\alpha\beta\xi_1} \partial_{\zeta} \partial_{\eta}
	\Big(c^3S^{30}_{\xi_1\zeta\eta}\partial_{\theta}\times O(1) +c^2S^{31}_{\xi_1\zeta\eta\zeta_1}\partial_{\theta}\times O(1) 
	\notag\\
	&\qquad\qquad\qquad+S^{33}_{\xi_1\zeta\eta\zeta_1\zeta_3}  \partial_{\theta}\big[\rho c_s^4\Delta_{\zeta_1\zeta_2\zeta_3\theta}+O(1/\Delta x^2)\big]
	\Big)\notag\\
	& +S^{2}_{\alpha\beta\xi_1\xi_2} \partial_{\zeta} \partial_{\eta} 
	\Big(c^3S^{40}_{\xi_1\xi_2\zeta\eta} \partial_{\theta}\times O(1)+c^2S^{31}_{\xi_1\zeta\eta\zeta_1} \partial_{\theta}\times O(1)
	+cS^{42}_{\xi_1\xi_2\zeta\eta\zeta_1\zeta_2}  \partial_{\theta}\times O(1/\Delta x^2)\notag\\
	& 
	+S^{43}_{\xi_1\xi_2\zeta\eta\zeta_1\zeta_3} \partial_{\theta}\big(\rho c_s^4\Delta_{\zeta_1\zeta_2\zeta_3\theta}+O(1/\Delta x^2))
	+\frac{1}{c}S^{44}_{\xi_1\xi_2\zeta\eta\zeta_1\zeta_3} \partial_{\theta}\times O(1/\Delta x^4)\Big) 
	\bigg] =O(\Delta x^3),  
\end{align}  
\begin{align}\label{swsw}
	&\big[-\Delta t^2c^2\partial_t\bm{\hat{S}}_N^{-1}\bm{\mathcal{W}}_0\bm{\hat{S}}_N^{-1}\bm{\mathcal{W}}_0\bm{m}^{eq}\big]_{\alpha+1}\notag\\
	&=-\Delta t^2c^2\partial_t\sum_{i=1}^{q }\bigg[\sum_{k=1}^{q }\Big[\sum_{j=1}^{q }\bm{e}_{j\alpha}\bm{\overline{\Lambda}}_{jk}\Big]\bm{e}_{k\beta}\partial_{\beta}\bm{\overline{\Lambda}}_{ki}\bigg]\bm{e}_{i\gamma}\partial_{\gamma}f_i^{eq}\notag\\
	&=-\Delta t^2 \partial_t \bigg[ S^{10}_{\alpha}\partial_{\beta}\big(cS^{10}_{\gamma} \partial_{\gamma}\times O(1)
	+S^1_{\gamma\zeta_1} \partial_{\gamma}\times O(1)\big)+S^1_{\alpha\xi_1}\partial_{\beta}\big(cS^{20}_{\xi_1\beta} \partial_{\gamma}\times O(1)\notag\\
	&\qquad\qquad
	+S^{21}_{\xi_1\beta\zeta_1} \partial_{\gamma}\times O(1)+\frac{1}{c}S^{2}_{\xi_1\beta\zeta_1\zeta_2}\partial_{\gamma}\times O(1/\Delta x^2)\big)\bigg]=O(\Delta x^3), 
\end{align}  
\begin{align}\label{w2swsw}
	&\big[-\Delta t^3c^4\bm{\mathcal{W}}_0^2\bm{\hat{S}}_N^{-1}\bm{\mathcal{W}}_0\bm{\hat{S}}_N^{-1}\bm{\mathcal{W}}_0\bm{m}^{eq}\big]_{\alpha+1}\notag\\
	&\qquad =-\Delta t^3c^4\sum_{i=1}^{q }\bigg[\sum_{k=1}^{q }\Big[\sum_{j=1}^{q }\bm{e}_{j\alpha}\bm{e}_{j\beta}\partial_{\beta}\bm{e}_{j\gamma}\partial_{\gamma}\bm{\overline{\Lambda}}_{jk}\Big]\bm{e}_{k\zeta}\partial_{\zeta}\bm{\overline{\Lambda}}_{ki}\bigg]\bm{e}_{i\eta}\partial_{\eta}f_i^{eq}\notag\\
	&\qquad =-\Delta t^3 \partial_{\gamma}\partial_{\beta}\bigg[  S^{30}_{\alpha\beta\gamma}\partial_{\zeta}\big(c^3S^{10}_{\zeta}\partial_{\eta}\times O(1)+c^2S^{31}_{\zeta\zeta_1} \partial_{\eta}\times O(1)\big)\notag\\
	&\qquad\qquad\qquad\qquad 
	+S^{31}_{\alpha\beta\gamma\xi_1}\partial_{\zeta}\big(c^3S^{20}_{\xi_1\zeta} \partial_{\eta}\times O(1)+c^2S^{21}_{\xi_1\zeta\zeta_1}\partial_{\eta}\times O(1)\notag\\
	&\qquad\qquad\qquad\qquad \qquad\qquad+cS^{2}_{\xi_1\zeta\zeta_1\zeta_2}\partial_{\eta}\times O(1/\Delta x^2)\big)\notag\\
	&\qquad\qquad\qquad\qquad 
	+S^{33}_{\alpha\beta\gamma\xi_1\xi_2\xi_3}\partial_{\zeta}\Big(c^3S^{40}_{\xi_1\xi_2\xi_3\zeta} \partial_{\eta}\times O(1)
	+c^2S^{41}_{\xi_1\xi_2\xi_3\zeta\zeta_1} \partial_{\eta}\times O(1)\notag\\
	&	\qquad\qquad\qquad\qquad\qquad\qquad\qquad +cS^{42}_{\xi_1\xi_2\xi_3\zeta\zeta_1\zeta_2} \partial_{\eta}\times O(1/\Delta x^2)\notag\\
	&\qquad\qquad\qquad\qquad\qquad\qquad\qquad 
	+S^{43}_{\alpha\beta\gamma\xi_1\xi_2\xi_3\zeta_1\zeta_2\zeta_3} \partial_{\eta}\times \big[\rho c_s^4\Delta_{\zeta_1\zeta_2\zeta_3\eta}+O(1/\Delta x^2)\big]\notag\\
	&\qquad\qquad\qquad\qquad\qquad\qquad\qquad 
	+\frac{1}{c}S^{44}_{\alpha\beta\gamma\xi_1\xi_2\xi_3\zeta_1\zeta_2\zeta_3} \partial_{\eta}\times O(1/\Delta x^4)\Big) 
	\bigg ]=O(\Delta x^3), 
\end{align}  
and
\begin{align}\label{wswswsw}
	&\big[\Delta t^3c^4\bm{\mathcal{W}}_0\bm{\hat{S}}_N^{-1}\bm{\mathcal{W}}_0\bm{\hat{S}}_N^{-1}\bm{\mathcal{W}}_0\bm{\hat{S}}_N^{-1}\bm{\mathcal{W}}_0\big]_{\alpha+1}\notag\\
	&=\Delta t^3c^4\sum_{l=1}^{q}\Bigg[\sum_{i=1}^{q }\bigg[\sum_{k=1}^{q}\partial_{\beta}\Big[\sum_{j=1}^{q}\bm{e}_{j\alpha}\bm{e}_{j\beta}\bm{\overline{\Lambda}}_{jk}\Big]\bm{e}_{k\gamma}\partial_{\gamma}\bm{\overline{\Lambda}}_{ki}\bigg]\bm{e}_{i\zeta}\partial_{\zeta}\bm{\overline{\Lambda}}_{il}\Bigg]\bm{e}_{l\eta}\partial_{\eta}f_l^{eq}\notag\\
	&=\Delta t^3  \partial_{\beta}\Bigg[  S^{20}_{\alpha\beta}\partial_{\gamma}S^{10}_{\gamma}\partial_{\zeta}\Big(c^3S^{10}_{\zeta} \partial_{\eta}\times O(1)+c^2S^{1}_{\zeta\chi_1}\partial_{\eta}\times O(1)\Big)\notag\\
	& 
	+S^{20}_{\alpha\beta}\partial_{\gamma}S^1_{\gamma\zeta_1} \partial_{\zeta}\Big(c^3S^{20}_{\zeta_1\zeta}\partial_{\eta}\times O(1)+c^2S^{21}_{\zeta_1\zeta\chi_1} \partial_{\eta}\times O(1)\notag\\
	&\qquad\qquad\qquad\qquad \qquad\qquad+cS^{2}_{\zeta_1\zeta\chi_1\chi_2} \partial_{\eta}\times O(1/\Delta x^2)\Big)\notag\\
	&
	+S^{21}_{\alpha\beta\xi_1}\partial_{\gamma}S^{20}_{\xi_1\gamma} \partial_{\zeta}\Big(c^3S^{10}_{\zeta}\partial_{\eta}\times O(1)+c^2S^1_{\zeta\chi_1} \partial_{\eta}\times O(1)\Big)\notag\\
	& 
	+S^{21}_{\alpha\beta\xi_1}\partial_{\gamma} S^{21}_{\xi_1\gamma\zeta_1}  \partial_{\zeta}\Big(c^3S^{20}_{\zeta_1\zeta}\partial_{\eta}\times O(1)+c^2S^{21}_{\zeta_1\zeta\chi_1} \partial_{\eta}\times O(1)\notag\\
	&\qquad\qquad\qquad\qquad \qquad\qquad+cS^{2}_{\zeta_1\zeta\chi_1\chi_2} \partial_{\eta}\times O(1/\Delta x^2)\Big)\notag\\
	& 
	+S^{21}_{\alpha\beta\xi_1}\partial_{\gamma}\partial_{\zeta} S^{2}_{\xi_1\gamma\zeta_1\zeta_2}  \partial_{\zeta}\Big(c^3S^{30}_{\zeta_1\zeta_2\zeta} \partial_{\eta}\times O(1)+c^2S^{31}_{\zeta_1\zeta_2\zeta\chi_1} \partial_{\eta}\times O(1)  \notag\\
	&\qquad\qquad\qquad\qquad\qquad +S^{33 }_{\zeta_1\zeta_2\zeta\chi_1\chi_2\chi_3} \partial_{\eta}\big[\rho c_s^4\Delta_{\chi_1\chi_2\chi_3\eta}+O(1/\Delta x^2)\big]\Big)\notag\\
	&\qquad\qquad\qquad\qquad\qquad+S^{ 2}_{\alpha\beta\xi_1\xi_2}\partial_{\gamma} S^{30}_{\xi_1\xi_2\gamma}\partial_{\zeta}\big(c^3S^{10}_{\zeta} \partial_{\eta}\times O(1)+c^2S^{1}_{\zeta\chi_1}\partial_{\eta}\times O(1)\big)\notag\\
	&  +S^{2 }_{\alpha\beta\xi_1\xi_2}\partial_{\gamma}S^{31}_{\xi_1\xi_2\gamma\zeta_1} \partial_{\zeta}\big(c^3S^{20}_{\zeta_1\zeta}\partial_{\eta}\times O(1)+c^2S^{21}_{\zeta_1\zeta\chi_1} \partial_{\eta}\times O(1)\notag\\
	&\qquad\qquad\qquad\qquad \qquad\qquad+cS^{2}_{\zeta_1\zeta\chi_1\chi_2} \partial_{\eta}\times O(1/\Delta x^2)\big)\notag\\
	& 
	+S^{2 }_{\alpha\beta\xi_1\xi_2}\partial_{\gamma} 
	S^{33 }_{\xi_1\xi_2\gamma\zeta_1\zeta_2\zeta_3} \partial_{\zeta}\Big(c^3S^{40}_{\zeta_1\zeta_2\zeta_3\zeta} \partial_{\eta}\times O(1)
	+c^2S^{41}_{\zeta_1\zeta_2\zeta_3\zeta\chi_1} \partial_{\eta}\times O(1)
	\notag\\
	&\qquad\qquad\qquad\qquad \qquad\qquad+cS^{42}_{\zeta_1\zeta_2\zeta_3\zeta\chi_1\chi_2} \partial_{\eta}\times O(1/\Delta x^2)\notag\\
	&\qquad\qquad\qquad\qquad\qquad +S^{43}_{\zeta_1\zeta_2\zeta_3\zeta\chi_1\chi_2\chi_3} \partial_{\eta}\times \big[\rho c_s^4\Delta_{\chi_1\chi_2\chi_3\eta}+O(1/\Delta x^2)\big]\notag\\
	&\qquad\qquad\qquad\qquad\qquad 
	+\frac{1}{c}S^{44}_{\zeta_1\zeta_2\zeta_3\zeta\chi_1\chi_2\chi_3\chi_4} \partial_{\eta}\times O(1/\Delta x^4)\Big)  \Bigg]=O(\Delta x^3). 
\end{align} 
With Eqs. (\ref{w3}$\sim$\ref{wswswsw}), the second-order ME of GPMRT-LB model and GPMFD scheme (\ref{modified-equation-third-order-diffusive}) for the momentum equation (\ref{momentum-equation}) becomes
\begin{align}\label{appen-above}
	&\Bigg[\bm{\tilde{F}}-\Big(c\bm{\mathcal{W}}_0+\partial_t\bm{I}+\Delta tc^2\Big(1-\frac{b}{2a^2}\Big)\bm{\mathcal{W}}_0^2  +\Delta t\bm{\mathcal{W}}_0\partial_t+\frac{\Delta t}{2}\partial_{tt}\bm{I}\Big)\bm{m}^{eq}-\Delta tc\bm{\mathcal{W}}_0\bm{\hat{S}}_N^{-1}\bm{\tilde{F}}\notag\\
	&+\bigg[\Delta tc^2\bm{\mathcal{W}}_0\bm{\hat{S}}_N^{-1}\bm{\mathcal{W}}_0+\Delta tc\bm{\mathcal{W}}_0\bm{\hat{S}}_N^{-1}\partial_t\bm{I}+\Delta t^2c^3\Big(1-\frac{b}{2a^2}\Big)\bm{\mathcal{W}}_0\bm{\hat{S}}_N^{-1}\bm{\mathcal{W}}_0^2\notag\\
	&\qquad\qquad +\Delta tc\partial_t\bm{I}\bm{\hat{S}}_N^{-1}\bm{\mathcal{W}}_0+\Delta t\partial_t\bm{\hat{S}}_N^{-1}\bm{I}\partial_t\bm{I}  +\Delta t^2c^3\Big(1-\frac{b}{2a^2}\Big)\bm{\mathcal{W}}_0^2\bm{\hat{S}}_N^{-1}\bm{\mathcal{W}}_0\notag\\
	&\qquad\qquad+\Delta t^2c^2\Big(1-\frac{b}{2a^2}\Big)\bm{\mathcal{W}}_0^2\bm{\hat{S}}_N^{-1}\partial_t\bm{I} +\Delta t^2c^2\bm{\mathcal{W}}_0\bm{\hat{S}}_N^{-1}\bm{\mathcal{W}}_0\partial_t\bm{I}\bigg]\bm{m}^{eq}\notag\\
	&-\bigg[\Delta t^2c^3\bm{\mathcal{W}}_0\bm{\hat{S}}_N^{-1}\bm{\mathcal{W}}_0\bm{\hat{S}}_N^{-1}\bm{\mathcal{W}}_0+\Delta   t^2c^2\bm{\mathcal{W}}_0\bm{\hat{S}}_N{-1}\bm{\mathcal{W}}_0\bm{\hat{S}}_N^{-1}\partial_t\bm{I}  \bigg]\bm{m}^{eq}\Bigg]_{\alpha+1} +O(\Delta x^3),
\end{align} 
then substituting Eqs. (\ref{second-order-remain}$\sim$\ref{third-order-remain-2}) into above Eq. (\ref{appen-above}), we have  
\begin{align}
	&\partial_t(\rho u_{\alpha})+\partial_{\beta}(\rho u_{\alpha}u_{\beta}+\rho c_s^2\delta_{\alpha\beta})\notag\\
	&= \Delta t  \partial_{\beta}\bigg[\Big[  S^{2}_{\alpha\beta\xi_1\xi_2}-\Big(1-\frac{b}{2a^2}\Big)\delta_{\alpha\xi_1}\delta_{\beta\xi_2} \Big]\partial_{\gamma}(\rho c_s^2\Delta_{\xi_1\xi_2\gamma\zeta}u_{\zeta})\bigg] \notag\\
	&  +\Delta t  \partial_{\beta}\Big[   S^{2}_{\alpha\beta\xi_1\xi_2}\partial_t(\rho c_s^2\delta_{\xi_1\xi_2}) \Big]+ \rho F_{x_\alpha}     =O(\Delta x^2) ,
\end{align} 
where the first-order ME (\ref{first-order-momentum-equation-appendix}) has been used.
\section{ Derivation of Eqs. (\ref{second-acoustic-NACDEs}), (\ref{acou-con-second}), and (\ref{mon-second-acou})}
\subsection{NACDE:  Derivation of Eq. (\ref{second-acoustic-NACDEs})}\label{appendix-NCDEs-acoustic}
Due to the fact that 
\begin{align}\label{first-acoustic-ncde-aapendix}
	&m^{eq}_1=\phi,     c\bm{\mathcal{W}}_0m^{eq}_1=\sum_{i=1}^{q}\bm{c}_{i\alpha}f_i^{eq}=B_{\alpha},  \notag\\
	&\bm{\tilde{F}}_{1}=\Big[\bm{M}\big( \bm{F} +\bm{G}+\frac{\Delta t}{2}\bm{D}\bm{F}\big)\Big]_{1}=S+\frac{\Delta t}{2 }\sum_{i=1}^{q} \big(\partial_t+\bm{c}_{i\gamma}\partial_{\gamma}\big) F_i =S+\frac{\Delta t}{2}\partial_tS ,  
\end{align} 
one can obtain the first-order ME: 
\begin{align}\label{first-order-NCDE-appendix-acoustic}
	&\frac{\partial \phi}{\partial t}+\nabla\cdot\bm{B} =S+O(\Delta x), 
\end{align}  
For the second-order ME, we also need to further consider the following equalities: 
\begin{align}\label{second-acoustic-ncde-aapendix}
	&\Big[-\frac{\Delta t}{2}\partial_t\bm{\tilde{F}}\Big]_{1}=-\frac{\Delta t}{2}\partial_tS+O(\Delta t^2),\notag\\
	&\Big[-\Delta tc\bm{\mathcal{W}}_0\bm{\hat{S}}_N^{-1}\bm{\tilde{F}}\Big]_{1}=-c\Delta t\partial_{\beta} \sum_{k=1}^{q }\sum_{j=1}^{q}\bigg[ \bm{e}_{j\beta}\bm{\overline{\Lambda}}_{jk}(F_k+G_k) \bigg]+O(\Delta t^2)\notag\\
	&\qquad\qquad  =- \Delta t\partial_{\beta}\sum_{k=1}^{q}\big(cS^{10}_{\beta}+S^{1}_{\beta\gamma}\bm{c}_{k\gamma}\big)(F_k+G_k)+O(\Delta t^2)\notag\\
	&\qquad\qquad =-c\Delta t\partial_{\beta}S^{10}_{\beta}S-\Delta t\partial_{\beta}\Big(S^{1}_{\beta\gamma}-\frac{\delta_{\beta\gamma}}{2}\Big) \partial_tB_{\gamma}\notag\\
	&\qquad\qquad= -\Delta t\partial_{\beta}\bigg[S^{1}_{\beta\gamma}-\Big(1-\frac{b}{2a^2}\Big)\delta_{\beta\gamma}\bigg] \partial_{\theta} C_{\gamma\theta}  +O(\Delta x^2),\notag\\
	&\bigg[c\bm{\mathcal{W}}_0\big(\bm{\hat{S}}_N^{-1}-\bm{I}/2\big)\partial_t\bm{m}^{eq}\bigg]_{1} =  \partial_{\beta}\sum_{k=1}^{q }\big(cS^{10}_{\beta}+S^{1}_{\beta\xi_1}\bm{c}_{k\xi_1}\big)\partial_tf_k^{eq}-\frac{1}{2}\partial_t\partial_{\beta}B_{\beta}\notag\\
	&\qquad\qquad=\partial_{\beta}\bigg[cS^{10}_{\beta}\partial_t\phi+\Big(S^1_{\beta\xi_1}-\frac{1}{2}\delta_{\beta\xi_1}\Big)\partial_tB_{\xi_1} \bigg],\notag\\
	&\bigg[c^2\bm{\mathcal{W}}_0 \bm{\hat{S}}_N^{-1} \bm{\mathcal{W}}_0\bm{m}^{eq}\bigg]_{1} =  \partial_{\beta}\sum_{k=1}^{q }\big(cS^{10}_{\beta}+S^{1}_{\beta\gamma}\bm{c}_{k\gamma}\big)\partial_{\theta}\bm{c}_{k\theta}f_k^{eq} \notag\\
	&\qquad\qquad =\frac{\partial}{\partial x_{\beta}}\bigg[cS^{10}_{\beta}\frac{\partial B_{\theta}}{\partial x_{\theta}}+S^{1}_{\beta\gamma}\frac{\partial \big(D_{\gamma \theta}c_s^2  \chi+\zeta C_{\gamma\theta}\big)}{ \partial x_{\theta}} \bigg],\notag\\
	& \bigg[c^2\bm{\mathcal{W}}_0\bm{\mathcal{W}}_0\bm{m}^{eq}\bigg]_{1}=  \partial_{\beta}\partial_{\theta} \sum_{j=1}^{q} 
	\bm{c}_{j\beta}  \bm{c}_{j\theta}f_j^{eq}=  \frac{\partial ^2 \big(\chi c_s^2D_{\beta\theta} +\zeta C_{\beta\theta}\big)}{\partial x_{\beta}x_{\theta}},
\end{align}
substituting Eqs. (\ref{first-acoustic-ncde-aapendix} ) and (\ref{second-acoustic-ncde-aapendix}) into Eq. (\ref{first-modified-equation-acoustic}), one can obtain
\begin{align}
	& \bigg[  \big(\partial_t\bm{I}+c\bm{\mathcal{W}}_0\big)\bm{m}^{eq} \bigg]_{1}+\frac{\Delta t}{2}\Big(1-\frac{b}{a^2}\Big)\Big[ c \bm{\mathcal{W}}_0^2\bm{m}^{eq}\Big]_{1}\notag\\
	&-\Delta t\bigg[c\bm{\mathcal{W}}_0\big(\bm{\hat{S}}_N^{-1}-\bm{I}/2\big)\big(\partial_t+c\bm{\mathcal{W}}_0\big)\bm{m}^{eq}\bigg]_{1}
	\notag\\
	&-\bm{\tilde{F}}_{1}+\frac{\Delta t}{2}\partial_t \bm{\tilde{F}}_{1}+\Delta t\Big[c\bm{\mathcal{W}}_0\bm{\hat{S}}_N^{-1}\bm{\tilde{F}}\Big]_{1}= \partial_t \phi +\partial_{\alpha}B_{\alpha}-S\notag\\
	&\qquad  -\Delta tc_s^2\frac{\partial }{\partial x_{\beta}}\bigg(\Big[S^1_{\beta\gamma}+\Big(\frac{b}{2a^2}-1\Big)\delta_{\beta\gamma}\Big]\bigg)\frac{\partial  D_{\gamma\theta} }{\partial x_{\theta}} \notag\\
	&+\underbrace{c\Delta t\partial_{\beta}S^{10}_{\beta}\Big[\partial_t\phi+\frac{\partial B_{\theta}}{\partial x_{\theta}}-S\Big]}_{O(\Delta x^2)} =O(\Delta x^2),
\end{align} 
where the first-order ME (\ref{first-order-NCDE-appendix-acoustic}) has been used.

\subsection{NSEs:  Derivation of Eqs. (\ref{acou-con-second}), and (\ref{mon-second-acou})}\label{appendix-NSEs-acoustic}
For the first conservative moment, i.e., the density $\rho$, one can obtain 
\begin{align}\label{acoustic-first-order-con-appendix}
	&m^{eq}_1=\rho,  \bm{\tilde{F}}_{1}=\Big[\bm{M}\big( \bm{F} +\frac{\Delta t}{2}\bm{D}\bm{F}\big)\Big]_{1}=0+\frac{\Delta t}{2 }\sum_{i=1}^{q} \big(\partial_t+\bm{c}_{i\gamma}\partial_{\gamma}\big)F_i =\frac{\Delta t}{2}\partial_{\gamma}(\rho \hat{F}_{\gamma}),\notag\\
	& \big[c\bm{\mathcal{W}}_0\bm{m}^{eq}\big]_1=\sum_{i=1}^{q }\bm{c}_{i\alpha}f_i^{eq}=\rho u_{ {\alpha}},   
\end{align} 
and for the second to $(d+1)_{th}$ conservative moments, i.e., the momentum $\rho u_{\alpha}$,  we have
\begin{align}\label{acoustic-first-order-mon-appendix}
	&m^{eq}_{\alpha+1}=\frac{\rho u_{\alpha}}{c}, \big[c\bm{\mathcal{W}}_0\bm{m}^{eq}\big]_{ \alpha+1}=\frac{1}{c}\partial_{\beta}\sum_{i=1}^{q }\bm{c}_{i\alpha}\bm{c}_{i\beta}f_i^{eq}=\frac{\big(\rho u_{ {\alpha}}u_{ {\beta}}+\rho c_s^2\delta_{\alpha\beta}\big)}{c},\notag \\
	&\big[\bm{\tilde{F}}\big]_{\alpha+1}=\Big[\bm{M}\big( \bm{F} +\frac{\Delta t}{2}\bm{D}\bm{F}\big)\Big]_{\alpha+1}=\frac{\rho  \hat{F}_{ x_{\alpha}}}{c}+\frac{\Delta t}{2c}\sum_{i=1}^{q }\bm{c}_{i\beta}\big(\partial_t+\bm{c}_{i\gamma}\partial_{\gamma}\big)F_i\notag \\
	&\qquad\quad\: =\frac{\rho  \hat{F}_{x_{\alpha}}}{c}+\frac{\Delta t}{2c}\partial_t(\rho \hat{F}_{ x_{\alpha}})
	+\frac{\Delta t}{2c}\partial_{\gamma}\big(\rho \hat{F}_{\beta}u_{\gamma}+\rho \hat{F}_{\gamma}u_{\beta}\big), 
\end{align} 
then the corrseponding first-order MEs of the GPMRT-LB model and GPMFD scheme for the NSEs (\ref{NSEs}) are given by
\begin{subequations}
	\begin{align}
		&\partial_t\rho +\partial_{\beta}(\rho u_{\beta})=O(\Delta x),\label{con-appendix-first}\\
		&\partial_t(\rho u_{\alpha})+\partial_{\alpha}\partial_{\beta}\big(\rho u_{\alpha}u_{\beta}+c_s^2\rho\delta_{\alpha\beta}\big)-\rho \hat{F}_{x_{\alpha}}=O(\Delta x).\label{mon-appendix-first}
	\end{align}
\end{subequations}
To obtain the second-order ME, the following equalities are needed:
\begin{align}\label{acoustic-second-order-con-appendix}
	&-\Delta t\Big[c\bm{\mathcal{W}}_0\bm{\hat{S}}_N^{-1}\bm{\tilde{F}}\Big]_{1}=-\Delta t\partial_{\beta} \sum_{k=1}^{q}\sum_{j=1}^{q }\bigg[ \bm{c}_{j\beta}\bm{\overline{\Lambda}}_{jk}F_k \bigg]+O(\Delta t^2)\notag \\
	&\qquad\qquad\qquad\qquad\quad=-\Delta t\partial_{\beta}\sum_{k=1}^{q }\big(S^{10}_{\beta}+S^{1}_{\beta\gamma}\bm{c}_{k\xi_1}\big)F_k=-\Delta t\partial_{\beta}S^{1}_{\beta\xi_1}\rho \hat{F}_{\xi_1}+O(\Delta x^2),\notag \\
	&\Big[-\frac{\Delta t}{2}\partial_t\bm{\tilde{F}}\Big]_{1}=\bm{0}, \bigg[c\bm{\mathcal{W}}_0 \partial_t\bm{m}^{eq}\bigg]_{1}=\partial_{\beta}\partial_t\sum_{i=1}^{q }  \bm{c}_{i\beta} f_i^{eq}= \partial_{\beta}\partial_t(\rho u_{\beta}),\notag \\
	&\bigg[c\bm{\mathcal{W}}_0\bm{\hat{S}}_N^{-1}\partial_t\bm{m}^{eq}\bigg]_{1}=\partial_{\beta}\sum_{k=1}^{q }\sum_{j=1}^{q } \bm{c}_{i\beta}\bm{{\Lambda}}_{jk}\partial_tf_k^{eq}= \partial_{\beta}\sum_{k=1}^{q}\big(cS^{10}_{\beta}+S^{1 }_{\beta\xi_1}\bm{c}_{k\xi_1}\big)\partial_tf_k^{eq}\notag \\
	&\qquad\qquad\qquad\qquad=\partial_{\beta}\bigg[cS^{10}_{\beta}\partial_t \rho + S^{1 }_{\beta\xi_1} \partial_t(\rho u_{\xi_1})\bigg],\notag \\
	&\bigg[c^2\bm{\mathcal{W}}_0\bm{\hat{S}}_N^{-1}\bm{\mathcal{W}}_0\bm{m}^{eq}\bigg]_{1} = \partial_{\beta}\sum_{k=1}^{q }\sum_{j=1}^{q }\bm{c}_{j\beta}
	\bm{{\Lambda}}_{jk}\bm{c}_{k\theta}\partial_{\theta}f_k^{eq}\notag \\
	&\qquad=\partial_{\beta}\partial_{\theta}\sum_{k=1}^{q}\big(cS^{10}_{\beta}+S^{1 }_{\beta\xi_1}\bm{c}_{k\xi_1}\big)\bm{c}_{k\theta}\partial_{\theta}f_k^{eq}\notag\\
	&\qquad=\frac{\partial}{\partial x_{\beta}}\Bigg[ cS^{10}_{\beta}\frac{\partial}{ \partial x_{\theta}}(\rho u_\theta)+S^{1}_{\beta\xi_1}\frac{\partial}{ \partial x_{\theta}}\big(\rho u_{\xi_1}u_{\theta}+\rho c_s^2\delta_{\xi_1\theta}\big) \Bigg],\notag \\
	& \bigg[c^2\bm{\mathcal{W}}_0\bm{\mathcal{W}}_0\bm{m}^{eq}\bigg]_{1}= \partial_{\beta}\partial_{\theta} \sum_{j=1}^{q } 
	\bm{c}_{j\beta}  \bm{c}_{j\theta}f_j^{eq}=\frac{\partial ^2\big(\rho u_{\beta}u_{\theta}+\rho c_s^2\delta_{\beta\theta}\big)}{\partial x_{\beta}x_{\theta}},
\end{align} 
combining Eqs. (\ref{acoustic-first-order-con-appendix}) and (\ref{acoustic-second-order-con-appendix}) yields 
\begin{align}
	& \bigg[  \big(\partial_t\bm{I}+c\bm{\mathcal{W}}_0\big)\bm{m}^{eq} \bigg]_{1}+\frac{\Delta t}{2}\Big(1-\frac{b}{a^2}\Big)\Big[ c \bm{\mathcal{W}}_0^2\bm{m}^{eq}\Big]_{1}\notag\\
	&-\Delta t\bigg[c\bm{\mathcal{W}}_0\big(\bm{\hat{S}}_N^{-1}-\bm{I}/2\big)\big(\partial_t+c\bm{\mathcal{W}}_0\big)\bm{m}^{eq}\bigg]_{1}-\bm{\tilde{F}}_{1}+\frac{\Delta t}{2}\partial_t\bm{\tilde{F}}_{1}+\Delta t\Big[c\bm{\mathcal{W}}_0\bm{\hat{S}}_N^{-1}\bm{\tilde{F}}\Big]_{1}\notag\\
	&\qquad= \partial_t \rho +\partial_{\alpha}\big(\rho u_{\alpha} ) -\Delta t\Big(\frac{b}{2a^2}-\frac{1}{2}\Big)\frac{\partial ^2\big(\rho u_{\beta}u_{\theta}+\rho c_s^2\delta_{\beta\theta}\big)}{\partial x_{\beta}x_{\theta}}\underbrace{+\partial_{\beta}cS^{10}_{\beta}\Big[\partial\rho+\partial_{\theta}(\rho u_{\theta})\Big]}_{O(\Delta x^2)}\notag\\
	& -\underbrace{\partial_{\beta}\Big(S^1_{\beta\xi_1}-\frac{1}{2}\delta_{\beta\xi_1}\Big)\Big[\partial_t(\rho u_{\xi_1})+\partial_{\theta}\big(\rho u_{\xi_1\theta}+\rho c_s^2\delta_{\xi_1\theta}\big)-\rho \hat{F}_{x_{\xi_1}}\Big]}_{O(\Delta x^2)}+O(\Delta x^2)\notag\\
	&\qquad=\partial_t \rho +\partial_{\alpha}\big(\rho u_{\alpha} )-\Delta t\Big(\frac{b}{2a^2}-\frac{1}{2}\Big)\frac{\partial ^2\big(\rho u_{\beta}u_{\theta}+\rho c_s^2\delta_{\beta\theta}\big)}{\partial x_{\beta}x_{\theta}}+O(\Delta x^2 ),
\end{align}
where the first-order MEs (\ref{con-appendix-first}) and (\ref{mon-appendix-first}) have been used. In addition, for any $\alpha\in\{1\sim d\}$, we can derive the following equations:
\begin{align}\label{acoustic-second-order-mon-appendix}
	&\Big[-\frac{\Delta t}{2}\partial_t \bm{\tilde{F}}\Big]_{\alpha+1}=-\frac{\Delta t}{2c}\partial_t(\rho \hat{F}_{ x_{\alpha}})+O(\Delta x^2),\notag \\
	&\Big[-\Delta tc\bm{\mathcal{W}}_0\bm{\hat{S}}_N^{-1}\bm{\tilde{F}}\Big]_{\alpha+1}=-\Delta t\frac{1}{c}\sum_{k=1}^{q}\sum_{j=1}^{q }\bigg[\bm{c}_{j\alpha}\bm{c}_{j\beta}\partial_{\beta}\bm{\overline{\Lambda}}_{jk}F_k \bigg]+O(\Delta t^2),\notag \\
	&\qquad\qquad\qquad=-\Delta t\frac{1}{c}\partial_{\beta}\sum_{k=1}^{q }\Bigg(S^{20}_{\alpha\beta}+S^{21}_{\alpha\beta \xi_1}\bm{c}_{k\xi_1}+S^{22}_{\alpha\beta\xi_1\xi_2}\bm{c}_{k\xi_1}\bm{c}_{k\xi_2}\Bigg)F_k+O(\Delta t^2)\notag \\
	&\qquad\qquad\qquad=-\Delta t\frac{1}{c}\partial_{\beta}\Big[S^{21}_{\alpha\beta \xi_1}(\rho \hat{F}_{\xi_1}) +S^{22}_{\alpha\beta\xi_1\xi_2}(\rho \hat{F}_{\xi_1}u_{\xi_2}+\rho \hat{F}_{\xi_2}u_{\xi_1})\Big]+O(\Delta t^2),\notag \\
	&\bigg[c\bm{\mathcal{W}}_0   \partial_t\bm{m}^{eq}\bigg]_{\alpha+1} =\frac{1}{c} \sum_{k=1}^{q} \bm{c}_{k\alpha}\bm{c}_{k\beta}\partial_{\beta} \partial_tf_k^{eq}=\frac{1}{c} \partial_{\beta}  \partial_t(\rho c_s^2\delta_{\alpha\beta}+\rho u_{\alpha}u_{\beta}) ,\notag \\
	&\bigg[c\bm{\mathcal{W}}_0 \bm{\hat{S}}_N^{-1} \partial_t\bm{m}^{eq}\bigg]_{\alpha+1}=\frac{1}{c} \partial_{\beta}\sum_{k=0}^{q-1}\big(S^{20}_{\alpha\beta}+S^{21}_{\alpha\beta\xi_1}\bm{c}_{k\xi_1}+S^{22}_{\alpha\beta\xi_1\xi_2}\bm{c}_{k\xi_1}\bm{c}_{k\xi_2}\big)\partial_tf_k^{eq}\notag \\
	&\qquad\qquad\qquad=\frac{1}{c} \partial_{\beta}\Big(S^{20}_{\alpha\beta}\partial_t\rho+S^{21}_{\alpha\beta\xi_1}\partial_t(\rho u_{\xi_1})+S^{22}_{\alpha\beta\xi_1\xi_2}\partial_t(\rho c_s^2\delta_{\xi_1\xi_2}+\rho u_{\xi_1}u_{\xi_2})\Big),\notag \\
	&\bigg[c^2\bm{\mathcal{W}}_0 \bm{\hat{S}}_N^{-1} \bm{\mathcal{W}}_0\bm{m}^{eq}\bigg]_{\alpha+1} =\frac{1}{c}\partial_{\beta}\sum_{k=1}^{q}\big(S^{20}_{\alpha\beta}+S^{21}_{\alpha\beta\xi_1}\bm{c}_{k\xi_1}+S^{22}_{\alpha\beta\xi_1\xi_2}\bm{c}_{k\xi_1}\bm{c}_{k\xi_2}\big)\bm{c}_{k\theta}\partial_{\theta}f_k^{eq}\notag \\
	&=\frac{1}{c}\partial_{\beta}\Big[S^{20}_{\alpha\beta}\partial_{\theta}(\rho u_{\theta})+S^{21}_{\alpha\beta\xi_1}\partial_{\theta}(\rho u_{\theta}u_{\xi_1}+\rho c_s^2\delta_{\xi_1\theta}) +S^{22}_{\alpha\beta\xi_1\xi_2}\partial_{\theta}(\rho c_s^2\Delta _{\xi_1\xi_2\theta\zeta u_{\zeta}})\Big],\notag \\
	& \bigg[c^2\bm{\mathcal{W}}_0\bm{\mathcal{W}}_0\bm{m}^{eq}\bigg]_{\alpha+1}=\frac{1}{c} \partial_{\beta}\partial_{\theta} \sum_{j=1}^{q }\bm{c}_{j\alpha}
	\bm{c}_{j\beta}  \bm{c}_{j\theta}f_j^{eq}=\frac{1}{c} \partial_{\beta}\partial_{\theta} (\rho c_s^2\Delta_{\alpha\beta\theta\zeta} u_{\zeta}) ,
\end{align} 
considering  Eqs. (\ref{acoustic-first-order-mon-appendix}) and (\ref{acoustic-second-order-mon-appendix}) yields
\begin{align}
	& \bigg[  \big(\partial_t\bm{I}_q+c\bm{\mathcal{W}}_0\big)\bm{m}^{eq} \bigg]_{\alpha+1}+\frac{\Delta t}{2}\Big(1-\frac{b}{a^2}\Big)\Big[ c^2 \bm{\mathcal{W}}_0^2\bm{m}^{eq}\Big]_{\alpha+1}\notag \\
	&-\Delta t\bigg[c\bm{\mathcal{W}}_0\big(\bm{\hat{S}}_N^{-1}-\bm{I}_q/2\big)\big(\partial_t+c\bm{\mathcal{W}}_0\big)\bm{m}^{eq}\bigg]_{\alpha+1}\notag\\
	&-\bm{\tilde{F}}_{\alpha+1}+\frac{\Delta t}{2}\partial_t\bm{\tilde{F}}_{\alpha+1}+\Delta t\Big[c\bm{\mathcal{W}}_0\bm{\hat{S}}_N^{-1}\bm{\tilde{F}}\Big]_{\alpha+1}=\frac{1}{c}\partial_t(\rho u_{\alpha})+\partial_{\beta}\big(\rho u_{\alpha}u_{\beta}+\rho c_s^2\delta_{\alpha\beta})\notag \\
	& +\frac{\Delta t}{2c}\Big(1-\frac{b}{2a^2}\Big)\partial_{\beta}\partial_{\theta}\big(\rho c_s^2\Delta_{\alpha\beta\theta\zeta}u_{\zeta}\big) -\frac{\rho  \hat{F}_{x_{\alpha}}}{c}-\underbrace{\Delta t\frac{1}{c}\partial_{\beta}S^{20}_{\alpha\beta}\Big(\partial_t\rho +\partial_{\theta}(\rho u_{\theta})\Big)}_{O(\Delta x^2)} \notag \\
	&-\underbrace{\Delta t\frac{1}{c}\partial_{\beta} S^{21}_{\alpha\beta\xi_1}\Big(\partial_t(\rho u_{\xi_1}) +\partial_{\theta}(\rho u_{\theta}u_{\xi_1}+\rho c_s^2\delta_{\xi_1\theta})-\rho\hat{F}_{x_{\xi_1}}\Big)}_{O(\Delta x^2)}\notag \\
	&-\Delta t\frac{1}{c}\partial_{\beta}\Big(S^{22}_{\alpha\beta\xi_1\xi_2}-\frac{1}{2}\delta_{\xi_1\alpha}\delta_{\xi_2\beta}\Big)\Big(\partial_t(\rho c_s^2\delta_{\xi_1\xi_2}+\rho u_{\xi_1}u_{\xi_2})+\partial_{\theta}(\rho c_s^2\Delta_{\xi_1\xi_2\theta\zeta}u_{\zeta})\notag \\
	&-\rho \hat{F}_{x_{\xi_1}}u_{\xi_2}-\rho \hat{F}_{x_{\xi_2}}u_{\xi_1}\Big) =O(\Delta x^2 ),
\end{align}
where the results of first-order MEs (\ref{con-appendix-first}) and (\ref{mon-appendix-first}) have been adopted.
\section{The parameters of GPMFD scheme  (\ref{fourth-scheme})}\label{coe-fourth-order-fd-scheme}
The parameters $\alpha_i(i\in\{1\sim 3\}),\beta_k,\gamma_k(k\in\{1\sim 5\})$  of the GPMFD scheme  (\ref{fourth-scheme}) are given by
\begin{align}
	& \alpha_1=\frac{1}{c^2\theta}\big(6a^2c^2 - 4a^2bc^2 - 8a^2c^2s_1 - 2a^2c^2s_2 - bc^2s_1^2 - 2b^2c^2s_1 + 2a^2c^2s_1^2 + b^2c^2s_1^2 \notag\\
	&\qquad+ 3bc^2s_1 - 2a^2bc^2s_1^2 - bc^2s_1s_2 + 6a^2bc^2s_1 + 2a^2c^2s_1s_2 - 2a^2bs_2u^2 + 2a^2bc^2s_2w_0 \notag\\
	&\qquad + a^2bs_1s_2u^2 + b^2c^2s_1s_2w_0 - 2a^2bc^2s_1s_2w_0\big), \notag\\
	&\alpha_2=\frac{1}{2c^2\theta}\big(4a^2bc^2 + 2b^2c^2s_1 - b^2c^2s_1^2 + 2a^2bc^2s_1^2 + 2a^3cs_1u - 6a^2bc^2s_1 + 2a^2bs_2u^2  \notag \\
	&\qquad- 2a^3cs_1^2u- 2a^2bc^2s_2w_0 - a^2bs_1s_2u^2 - b^2c^2s_1s_2w_0 + abcs_1^2u + 2a^2bc^2s_1s_2w_0\big), \notag\\         
	&\alpha_3=-\frac{1}{2c^2\theta}\big(b^2c^2s_1^2 - 2b^2c^2s_1 - 4a^2bc^2 - 2a^2bc^2s_1^2 + 2a^3cs_1u + abcs_1^2u - 2a^2bc^2s_1s_2w_0 \notag\\
	&\qquad+ 6a^2bc^2s_1 - 2a^2bs_2u^2 - 2a^3cs_1^2u + 2a^2bc^2s_2w_0 + a^2bs_1s_2u^2 + b^2c^2s_1s_2w_0\big), \notag\\
	&\beta_1=-\frac{1}{2c^2\theta}\big(12a^2c^2 + 2a^4c^2 - 16a^2bc^2 - 20a^2c^2s_1 - 8a^2c^2s_2 - 4a^4c^2s_1- 4bc^2s_1^2   \notag\\
	&\qquad - 8b^2c^2s_1+ 3b^3c^2s_1 + 2a^4s_2u^2 + 6a^2b^2c^2 + 8a^2c^2s_1^2 + 2a^4c^2s_1^2+ 6b^2c^2s_1^2  \notag \\
	&\qquad - 3b^3c^2s_1^2 + 6bc^2s_1- 13a^2bc^2s_1^2 - 12a^2b^2c^2s_1 - 4a^2c^2s_1^2s_2+ 6a^2b^2s_2u^2  \notag \\
	&\qquad- 2b^2c^2s_1^2s_2 + a^4s_1^2s_2u^2 - 4bc^2s_1s_2 + 6a^2b^2c^2s_1^2 + 29a^2bc^2s_1+ 8a^2bc^2s_2   \notag \\
	&\qquad + 12a^2c^2s_1s_2 - 4a^2bs_2u^2 + 2bc^2s_1^2s_2+ 4b^2c^2s_1s_2 - 2a^4c^2s_2w_0- 3a^4s_1s_2u^2   \notag \\
	&\qquad + 4a^2bc^2s_2w_0 + 6a^2bs_1s_2u^2+ 4a^4c^2s_1s_2w_0  + 2b^2c^2s_1s_2w_0 - 3b^3c^2s_1s_2w_0 \notag\\
	&\qquad - 6a^2b^2c^2s_2w_0 - 2a^2bs_1^2s_2u^2- 9a^2b^2s_1s_2u^2 - 2a^4c^2s_1^2s_2w_0 - 2b^2c^2s_1^2s_2w_0  \notag \\
	&\qquad+ 3a^2b^2s_1^2s_2u^2 + 5a^2bc^2s_1^2s_2w_0 + 12a^2b^2c^2s_1s_2w_0  - 6a^2b^2c^2s_1^2s_2w_0,\notag \\
	&\qquad- 9a^2bc^2s_1s_2w_0+ 4a^2bc^2s_1^2s_2 - 12a^2bc^2s_1s_2+ 3b^3c^2s_1^2s_2w_0 \big)\notag\\
	&\beta_2=\frac{1}{2c^2\theta}\big(2b^3c^2s_1 - 4b^2c^2s_1 - 8a^2bc^2 + 4a^2b^2c^2- 4a^2b^2c^2s_1^2s_2w_0 - 4a^2bc^2s_1s_2w_0 \notag \\
	&\qquad+ 3b^2c^2s_1^2 - 2b^3c^2s_1^2 - 6a^2bc^2s_1^2 - 8a^2b^2c^2s_1  + 4a^2b^2s_2u^2 - b^2c^2s_1^2s_2 \notag\\
	&\qquad + 4a^2b^2c^2s_1^2 + 14a^2bc^2s_1 + 4a^2bc^2s_2 - 2a^2bs_2u^2 + 2b^2c^2s_1s_2 + 2a^3cs_1^2u \notag\\
	&\qquad + 2a^2bc^2s_2w_0 + 3a^2bs_1s_2u^2 - 2a^3cs_1^2s_2u + b^2c^2s_1s_2w_0 - 2b^3c^2s_1s_2w_0 \notag\\
	&\qquad - 4a^2b^2c^2s_2w_0 - a^2bs_1^2s_2u^2 - 6a^2b^2s_1s_2u^2 - b^2c^2s_1^2s_2w_0 + 2b^3c^2s_1^2s_2w_0\notag \\
	&\qquad- 6a^2bc^2s_1s_2 + 2a^3cs_1s_2u + 2a^2b^2s_1^2s_2u^2 + 2a^2bc^2s_1^2s_2w_0  \notag\\
	&\qquad + 8a^2b^2c^2s_1s_2w_0- 2a^3cs_1u + 2a^2bc^2s_1^2s_2 - abcs_1^2u + abcs_1^2s_2u  \big), \notag\\
	&\beta_3=-\frac{(a^2 - b^2)(s_1 - 1)}{4c^2\theta}\notag\\
	&\qquad\times\big(2a^2c^2(1-s_1) +  a^2s_2(2-s_1)(u^2-2w_0c^2) + bc^2s_1 (1-s_2w_0)  \big) , \notag\\
	& \beta_3=\beta_2, \beta_5=\beta_4, \gamma_1=\frac{(s_1 - 1)(s_2 - 1)(- a^2  + b^2 )}{4}( 2a^2+ 6b^2  - 8b+ 4), \notag\\
	&\gamma_2=(s_1 - 1)(s_2 - 1)(b- b^2  ), \gamma_3=\gamma_2,\gamma_4=\frac{(s_1 - 1)(s_2 - 1)(- a^2  + b^2 )}{4}, \gamma_5=\gamma_4.
\end{align} 
where $\theta=1/(bs_1-2a^2s_1+2a^2)$.

\end{appendices}

\end{document}